%
%
%
%
%
\def\spaces{\space\space\space\space\space\space\space\space\space\space}
\def\spacess{\message{\spaces\spaces\spaces\spaces\spaces\spaces\spaces}}
\spacess
\spacess
\message{Annals of Mathematics Style: Current Version: 1.0. March 25, 1992}
\spacess
\spacess
\catcode`@=11



\font\elevensc=cmcsc10 scaled\magstephalf
\font\tensc=cmcsc10

\font\eightsc=cmcsc10 scaled800

\font\elevenrm=cmr10 scaled \magstephalf
\font\ninerm=cmr9
\font\eightrm=cmr8
\font\sixrm=cmr6
\font\fiverm=cmr5

\font\eleveni=cmmi10 scaled\magstephalf
\font\ninei=cmmi9
\font\eighti=cmmi8
\font\sixi=cmmi6
\font\fivei=cmmi5
\skewchar\ninei='177 \skewchar\eighti='177 \skewchar\sixi='177
\skewchar\eleveni= '177

\font\elevensy=cmsy10 scaled\magstephalf
\font\ninesy=cmsy9
\font\eightsy=cmsy8
\font\sixsy=cmsy6
\font\fivesy=cmsy5
\skewchar\ninesy='60 \skewchar\eightsy='60 \skewchar\sixsy='60
\skewchar\elevensy='60

\font\eighteenbf=cmbx10 scaled\magstep3

\font\twelvebf=cmbx10 scaled \magstep1
\font\elevenbf=cmbx10 scaled \magstephalf
\font\tenbf=cmbx10
\font\ninebf=cmbx9
\font\eightbf=cmbx8
\font\sixbf=cmbx6
\font\fivebf=cmbx5

\font\elevenit=cmti10 scaled\magstephalf
\font\nineit=cmti9
\font\eightit=cmti8

\font\eighteenmib=cmmib10 scaled \magstep3
\font\twelvemib=cmmib10 scaled \magstep1
\font\elevenmib=cmmib10 scaled\magstephalf
\font\tenmib=cmmib10
\font\eightmib=cmmib10 scaled 800
\font\sixmib=cmmib10 scaled 600

\font\eighteensyb=cmbsy10 scaled \magstep3
\font\twelvesyb=cmbsy10 scaled \magstep1
\font\elevensyb=cmbsy10 scaled \magstephalf
\font\tensyb=cmbsy10
\font\eightsyb=cmbsy10 scaled 800
\font\sixsyb=cmbsy10 scaled 600

\font\elevenex=cmex10 scaled \magstephalf
\font\tenex=cmex10
\font\eighteenex=cmex10 scaled \magstep3


\def\elevenpoint{\def\rm{\fam0\elevenrm}%
  \textfont0=\elevenrm \scriptfont0=\eightrm \scriptscriptfont0=\sixrm
  \textfont1=\eleveni \scriptfont1=\eighti \scriptscriptfont1=\sixi
  \textfont2=\elevensy \scriptfont2=\eightsy \scriptscriptfont2=\sixsy
  \textfont3=\elevenex \scriptfont3=\elevenex \scriptscriptfont3=\elevenex
  \def\bf{\fam\bffam\elevenbf}%
  \def\it{\fam\itfam\elevenit}%
  \textfont\bffam=\elevenbf \scriptfont\bffam=\eightbf
   \scriptscriptfont\bffam=\sixbf
\normalbaselineskip=13.9pt
  \setbox\strutbox=\hbox{\vrule height9.5pt depth4.4pt width0pt\relax}%
  \normalbaselines\rm}

\elevenpoint 

\def\ninepoint{\def\rm{\fam0\ninerm}%
  \textfont0=\ninerm \scriptfont0=\sixrm \scriptscriptfont0=\fiverm
  \textfont1=\ninei \scriptfont1=\sixi \scriptscriptfont1=\fivei
  \textfont2=\ninesy \scriptfont2=\sixsy \scriptscriptfont2=\fivesy
  \textfont3=\tenex \scriptfont3=\tenex \scriptscriptfont3=\tenex
  \def\it{\fam\itfam\nineit}%
  \textfont\itfam=\nineit
  \def\bf{\fam\bffam\ninebf}%
  \textfont\bffam=\ninebf \scriptfont\bffam=\sixbf
   \scriptscriptfont\bffam=\fivebf
\normalbaselineskip=11pt
  \setbox\strutbox=\hbox{\vrule height8pt depth3pt width0pt\relax}%
  \normalbaselines\rm}

\def\eightpoint{\def\rm{\fam0\eightrm}%
  \textfont0=\eightrm \scriptfont0=\sixrm \scriptscriptfont0=\fiverm
  \textfont1=\eighti \scriptfont1=\sixi \scriptscriptfont1=\fivei
  \textfont2=\eightsy \scriptfont2=\sixsy \scriptscriptfont2=\fivesy
  \textfont3=\tenex \scriptfont3=\tenex \scriptscriptfont3=\tenex
  \def\it{\fam\itfam\eightit}%
  \textfont\itfam=\eightit
  \def\bf{\fam\bffam\eightbf}%
  \textfont\bffam=\eightbf \scriptfont\bffam=\sixbf
   \scriptscriptfont\bffam=\fivebf
\normalbaselineskip=12pt
  \setbox\strutbox=\hbox{\vrule height8.5pt depth3.5pt width0pt\relax}%
  \normalbaselines\rm}


\def\eighteenbold{\def\rm{\fam0\eighteenrm}%
  \textfont0=\eighteenbf \scriptfont0=\twelvebf \scriptscriptfont0=\tenbf
  \textfont1=\eighteenmib \scriptfont1=\twelvemib\scriptscriptfont1=\tenmib
  \textfont2=\eighteensyb \scriptfont2=\twelvesyb\scriptscriptfont2=\tensyb
  \textfont3=\eighteenex \scriptfont3=\tenex \scriptscriptfont3=\tenex
  \def\bf{\fam\bffam\eighteenbf}%
  \textfont\bffam=\eighteenbf \scriptfont\bffam=\twelvebf
   \scriptscriptfont\bffam=\tenbf
\normalbaselineskip=20pt
  \setbox\strutbox=\hbox{\vrule height13.5pt depth6.5pt width0pt\relax}%
\everymath {\fam0 }
\everydisplay {\fam0 }
  \normalbaselines\bf}

\def\elevenbold{\def\rm{\fam0\elevenrm}%
  \textfont0=\elevenbf \scriptfont0=\eightbf \scriptscriptfont0=\sixbf%
  \textfont1=\elevenmib \scriptfont1=\eightmib \scriptscriptfont1=\sixmib%
  \textfont2=\elevensyb \scriptfont2=\eightsyb \scriptscriptfont2=\sixsyb%
  \textfont3=\elevenex \scriptfont3=\elevenex \scriptscriptfont3=\elevenex%
  \def\bf{\fam\bffam\elevenbf}%
  \textfont\bffam=\elevenbf \scriptfont\bffam=\eightbf%
   \scriptscriptfont\bffam=\sixbf%
\normalbaselineskip=14pt%
  \setbox\strutbox=\hbox{\vrule height10pt depth4pt width0pt\relax}%
\everymath {\fam0 }%
\everydisplay {\fam0 }%
  \normalbaselines\bf}


\hsize=31pc
\vsize=48pc

\parindent=22pt
\parskip=0pt

\widowpenalty=10000
\clubpenalty=10000

\topskip=12pt

\skip\footins=20pt
\dimen\footins=3in 

\abovedisplayskip=6.95pt plus3.5pt minus 3pt
\belowdisplayskip=\abovedisplayskip


\voffset=7pt
\def\bindingoffset{\ifodd\pageno \hoffset=-18pt\else\hoffset9pc\fi}

\newif\iftitle
%

\def\amheadline{\iftitle
\hbox to\hsize{\hss\currannalsline\hss}\else\line{\ifodd\pageno
{\fiverm\DATUM}\hfill\thetitle\hfill\llap{\elevenrm\folio}\else\rlap
{\elevenrm\folio}\hfill\theauthors\hfill{\fiverm\DATUM}\fi}\fi}

\headline={\amheadline}
\footline={\global\titlefalse}
\output={\bindingoffset\plainoutput}


\def\annalsline#1#2{\vfill\eject
\ifodd\pageno\else 
\line{\hfill}
\vfill\eject\fi
\global\titletrue
%
%
%
%
%
\def\currannalsline{\eightrm Manuscript, variant of
     #1, #2;  \  \thepages}}

\def\titleheadline#1{\def\one{#1}\ifx\one\empty\else
\gdef\thetitle{{\frenchspacing%
\let\\ \relax\eightsc\uppercase{#1}}}\fi}

\newif\ifshort

\let\shorttitle\titleheadline

\def\onpages#1#2{\def\thepages{#1--#2}}

\def\thismuchskip[#1]{\vskip#1pt}
\def\ilook{\ifx\next[ \let\go\thismuchskip\else
\let\go\relax\vskip1pt\fi\go}

\def\institution#1{\def\theinstitutions{\vbox{\baselineskip10pt
\def\and{{\eightrm and }}
\def\\{\futurelet\next\ilook}\eightsc #1}}}
\let\institutions\institution

\newwrite\auxfile

\def\startingpage#1{\def\one{#1}\ifx\one\empty\global\pageno=1\else
\global\pageno=#1\fi
\theoremcount=0 \eqcount=0 \sectioncount=0
\openin1 \jobname.aux \ifeof1
\onpages{\ifx\one\empty\else#1\fi}{}
\else\closein1 \relax\input \jobname.aux
\onpages{\ifx\one\empty\else#1\fi}{\lastpage}
\fi\immediate\openout\auxfile=\jobname.aux
}

\def\endarticle{\ifRefsUsed\global\RefsUsedfalse%
\else\vskip21pt\theinstitutions%
\nobreak\vskip8pt\vbox{\thereceived\therevised}\fi%
\write\auxfile{\string\def\string\lastpage{\the\pageno}}}

\outer\def\bye{\endarticle\par \vfill \supereject \end}

\newif\ifacks
\long\def\acknowledgements#1{\def\one{#1}\ifx\one\empty\else
\global\ackstrue\vskip-\baselineskip%
\footnote{\ \unskip}{*#1}\fi}

\def\title#1{\vbox to76pt{\vfill
\baselineskip=18pt
\parindent=0pt
\overfullrule=0pt
\hyphenpenalty=10000
\everypar={\hskip\parfillskip\relax}
\hbadness=10000
\def\\ {\vskip1sp}
\eighteenbold#1\vskip1sp}
\titleheadline{#1}}

\newif\ifauthor
\def\author#1{\vskip15pt
\hbox to\hsize{\hss\tenrm By \tensc#1\ifacks\global\acksfalse*\fi\hss}
\ifshort\else\xdef\theauthors{{\eightsc\uppercase{#1}}}\fi%
\vskip7pt
\vskip\baselineskip
\global\authortrue\everypar={\global\authorfalse\everypar={}}}

\def\twoauthors#1#2{\vskip15pt
\hbox to\hsize{\hss%
\tenrm By \tensc#1 {\tenrm and} #2\ifacks\global\acksfalse*\fi\hss}
\ifshort\else\xdef\theauthors{{\eightsc\uppercase{#1 and #2}}}\fi
\vskip7pt
\vskip\baselineskip
\global\authortrue\everypar={\global\authorfalse\everypar={}}}


\newcount\theoremcount
\newcount\sectioncount
\newcount\eqcount

\newif\ifspecialnumon

\def\eqnumber=#1 {\global\eqcount=#1 \global\advance\eqcount by-1\relax}
\def\sectionnumber=#1 {\global\sectioncount=#1
\global\advance\sectioncount by-1\relax}
\def\proclaimnumber=#1 {\global\theoremcount=#1
\global\advance\theoremcount by-1\relax}

\newif\ifsection
\newif\ifsubsection

\def\intro{\centerline{\bf Introduction}\global\everypar={}
\global\authorfalse\vskip6pt}

\def\elevenboldmath#1{$#1$\egroup}
\def\mathbold{\hbox\bgroup\elevenbold\elevenboldmath}

\def\section#1{\global\theoremcount=0
\global\eqcount=0
\ifauthor\global\authorfalse\else%
\vskip18pt plus 18pt minus 6pt\fi%
{\parindent=0pt\everypar={\hskip\parfillskip}%
\def\\ {\vskip1sp}\elevenpoint\bf%
\ifspecialnumon\global\specialnumonfalse$\rm\spnum$%
\gdef\sectnum{$\rm\spnum$}%
\else\interlinepenalty=10000%
\global\advance\sectioncount by1\relax\the\sectioncount%
\gdef\sectnum{\the\sectioncount}%
\fi.\hskip6pt#1\vrule width 0pt depth12pt}\hskip\parfillskip\break%
\global\sectiontrue%
\everypar={\global\sectionfalse\global\interlinepenalty=0\everypar={}}%
\ignorespaces}


\newif\ifspequation

\let\eqno\leqno 

\newif\ifineqalignno
\let\saveleqalignno\leqalignno
\def\leqalignno{\let\eqnu\Eeqnu\saveleqalignno}

\def\sectandeqnum{%
\ifspecialnumon\global\specialnumonfalse
$\rm\spnum$\gdef\eqnum{$\rm\spnum$}\else\global\firstlettertrue
\global\advance\eqcount by1
\ifappend\applett\else\the\sectioncount\fi.%
\the\eqcount
\xdef\eqnum{\ifappend\applett\else\the\sectioncount\fi.\the\eqcount}\fi}

\def\eqnu{\leqno{\hbox{\elevenrm\ifspequation\else(\fi\sectandeqnum
\ifspequation\global\spequationfalse\else)\fi}}}

\def\Speqnu{\global\setbox\leqnobox=\hbox{\elevenrm
\ifspequation\else%
(\fi\sectandeqnum\ifspequation\global\spequationfalse\else)\fi}}

\def\Eeqnu{\hbox{\elevenrm
\ifspequation\else%
(\fi\sectandeqnum\ifspequation\global\spequationfalse\else)\fi}}

\newif\iffirstletter
\global\firstlettertrue
\def\eqletter#1{\global\specialnumontrue\iffirstletter\global\firstletterfalse
\global\advance\eqcount by1\fi
\gdef\spnum{\ifappend\applett\else\the\sectioncount\fi.\the\eqcount#1}\eqnu}

\newbox\leqnobox
\def\outsideeqnu#1{\global\setbox\leqnobox=\hbox{#1}}

\def\eatone#1{}

\def\dosplit#1#2{\vskip-.5\abovedisplayskip
\setbox0=\hbox{$\let\eqno\outsideeqnu%
\let\eqnu\Speqnu\let\leqno\outsideeqnu#2$}%
\setbox1\vbox{\noindent\hskip\wd\leqnobox\ifdim\wd\leqnobox>0pt\hskip1em\fi%
$\displaystyle#1\mathstrut$\hskip0pt plus1fill\relax
\vskip1pt
\line{\hfill$\let\eqnu\eatone\let\leqno\eatone%
\displaystyle#2\mathstrut$\ifmathqed~~\qed\fi}}%
\copy1
\ifvoid\leqnobox
\else\dimen0=\ht1 \advance\dimen0 by\dp1
\vskip-\dimen0
\vbox to\dimen0{\vfill
\hbox{\unhbox\leqnobox}
\vfill}
\fi}

\everydisplay{\lookforbreak}

\long\def\lookforbreak #1$${\def\mathone{#1}%
\expandafter\testforbreak\mathone\splitmath @}

\def\testforbreak#1\splitmath #2@{\def\mathtwo{#2}\ifx\mathtwo\empty%
#1$$%
\ifmathqed\vskip-\belowdisplayskip
\setbox0=\vbox{\let\eqno\relax\let\eqnu\relax$\displaystyle#1$}%
\vskip-\ht0\vskip-3.5pt\hbox to\hsize{\hfill\qed}
\vskip\ht0\vskip3.5pt\fi
\else$$\vskip-\belowdisplayskip
\vbox{\dosplit{#1}{\let\eqno\eatone
\let\splitmath\relax#2}}%
\nobreak\vskip.5\belowdisplayskip
\noindent\ignorespaces\fi}


\newif\ifmathqed



\newcount\linenum
\newcount\colnum
\def\spline{\omit&\multispan{\the\colnum}{\hrulefill}\cr}
\def\colcounter{\ifnum\linenum=1\global\advance\colnum by1\fi}

\def\everyline{\noalign{\global\advance\linenum by1\relax}%
\ifnum\linenum=2\spline\fi}

\def\mtable{\bgroup\offinterlineskip
\everycr={\everyline}\global\linenum=0
\halign\bgroup\vrule height 10pt depth 4pt width0pt
\hfill$##$\hfill\hskip6pt\ifnum\linenum>1
\vrule\fi&&\colcounter\hskip12pt\hfill$##$\hfill\hskip12pt\cr}

\def\endmtable{\crcr\egroup\egroup}


\def\xast{*}
\newcount\intable
\newcount\mathcol
\newcount\savemathcol
\newcount\topmathcol
\newdimen\arrayhspace
\newdimen\arrayvspace

\arrayhspace=8pt 
\arrayvspace=12pt 

\newif\ifdollaron

\def\mathalign#1{\def\arg{#1}\ifx\arg\xast%
\let\go\relax\else\let\go\mathalign%
\global\advance\mathcol by1 %
\global\advance\topmathcol by1 %
\expandafter\def\csname  mathcol\the\mathcol\endcsname{#1}%
\fi\go}

\def\arraypickapart#1]#2*{\if#1c \ifmmode\vcenter\else
\global\dollarontrue$\vcenter\fi\else%
\if#1t\vtop\else\if#1b\vbox\fi\fi\fi\bgroup%
\def\one{#2}}

\def\arraystrut{\vrule height .7\arrayvspace depth .3\arrayvspace width 0pt}

\def\array#1#2*{\def\firstarg{#1}%
\if\firstarg[ \def\two{#2} \expandafter\arraypickapart\two*\else%
\ifmmode\vcenter\else\vbox\fi\bgroup \def\one{#1#2}\fi%
\global\everycr={\noalign{\global\mathcol=\savemathcol\relax}}%
\def\\ {\cr}%
\global\advance\intable by1 %
\ifnum\intable=1 \global\mathcol=0 \savemathcol=0 %
\else \global\advance\mathcol by1 \savemathcol=\mathcol\fi%
\expandafter\mathalign\one*%
\mathcol=\savemathcol %
\halign\bgroup&\hskip.5\arrayhspace\arraystrut%
\global\advance\mathcol by1 \relax%
\expandafter\if\csname mathcol\the\mathcol\endcsname r\hfill\else%
\expandafter\if\csname mathcol\the\mathcol\endcsname c\hfill\fi\fi%
$\displaystyle##$%
\expandafter\if\csname mathcol\the\mathcol\endcsname r\else\hfill\fi\relax%
\hskip.5\arrayhspace\cr}

\def\endarray{\crcr\egroup\egroup%
\global\mathcol=\savemathcol %
\global\advance\intable by -1\relax%
\ifnum\intable=0 %
\ifdollaron\global\dollaronfalse $\fi
\loop\ifnum\topmathcol>0 %
\expandafter\def\csname  mathcol\the\topmathcol\endcsname{}%
\global\advance\topmathcol by-1 \repeat%
\global\everycr={}\fi%
}

\def\big#1{{\hbox{$\left#1\vbox to 10pt{}\right.\n@space$}}}
\def\Big#1{{\hbox{$\left#1\vbox to 13pt{}\right.\n@space$}}}
\def\bigg#1{{\hbox{$\left#1\vbox to 16pt{}\right.\n@space$}}}
\def\Bigg#1{{\hbox{$\left#1\vbox to 19pt{}\right.\n@space$}}}


\newif\ifappend

\def\appendix#1#2{\def\applett{#1}\def\two{#2}%
\global\appendtrue
\global\theoremcount=0
\global\eqcount=0
\vskip18pt plus 18pt
\vbox{\parindent=0pt
\everypar={\hskip\parfillskip}
\def\\ {\vskip1sp}\elevenpoint\bf Appendix%
\ifx\applett\empty\gdef\applett{A}\ifx\two\empty\else.\fi%
\else\ #1.\fi\hskip6pt#2\vskip12pt}%
\global\sectiontrue%
\everypar={\global\sectionfalse\everypar={}}\nobreak\ignorespaces}

\newif\ifRefsUsed
\long\def\references{\global\RefsUsedtrue\vskip21pt
\theinstitutions
\everypar={}\global\bibnum=0
\vskip20pt\goodbreak\bgroup
\vbox{\centerline{\eightsc References}\vskip6pt}%
\ifdim\maxbibwidth>0pt\leftskip=\maxbibwidth\else\leftskip=18pt\fi%
\ifdim\maxbibwidth>0pt\parindent=-\maxbibwidth\else\parindent=-18pt\fi%
\everypar={\amref}%
\ninepoint
\frenchspacing
\nobreak\ignorespaces}

\def\endreferences{\vskip1sp\egroup\everypar={}%
\nobreak\vskip8pt\vbox{\thereceived\therevised}}

\newcount\bibnum

\def\amref#1 {\def\one{#1}\global\advance\bibnum by1%
\immediate\write\auxfile{\string\expandafter\string\def\string\csname
\space #1\string\endcsname{[\the\bibnum]}}%
\leavevmode\hbox to18pt{\hbox to13.2pt{\hss[\the\bibnum]}\hfill}}

\def\bibline{\hbox to30pt{\hrulefill}\/\/}

\def\name#1{{\eightsc#1}}

\newdimen\maxbibwidth
\def\AuthorRefNames [#1] {\def\amref{\spamref}
\setbox0=\hbox{[#1] }\global\maxbibwidth=\wd0\relax}

\def\spamref[#1] {\leavevmode\hbox to\maxbibwidth{\hss[#1]\hfill}}


\def\footnoterule{\kern-3pt\hrule width1in height.5pt\kern2.5pt}

\let\savefootnote\footnote
\def\footnote#1#2{%
\savefootnote{#1}{{\eightpoint\normalbaselineskip11pt
\normalbaselines#2}}}

\def\vfootnote#1{%
\insert \footins \bgroup \eightpoint\normalbaselineskip11pt\normalbaselines
\interlinepenalty \interfootnotelinepenalty
\splittopskip \ht \strutbox \splitmaxdepth \dp \strutbox \floatingpenalty
\@MM \leftskip \z@skip \rightskip \z@skip \spaceskip \z@skip
\xspaceskip \z@skip
{#1}$\,$\footstrut \futurelet \next \fo@t}


\newif\iffirstadded
\newif\ifadded

\def\addedlett{}

\def\alltheoremnums{%
\ifspecialnumon\global\specialnumonfalse
\ifadded\global\addedfalse
\iffirstadded\global\firstaddedfalse
\global\advance\theoremcount by1 \fi
\ifappend\applett\else\the\sectioncount\fi.\the\theoremcount\addedlett%
\xdef\theoremnum{\ifappend\applett\else\the\sectioncount\fi.%
\the\theoremcount\addedlett}%
\else$\rm\spnum$\def\theoremnum{$\rm\spnum$}\fi%
\else\global\firstaddedtrue
\global\advance\theoremcount by1
\ifappend\applett\else\the\sectioncount\fi.\the\theoremcount%
\xdef\theoremnum{\ifappend\applett\else\the\sectioncount\fi.%
\the\theoremcount}\fi}

\def\allcorolnums{%
\ifspecialnumon\global\specialnumonfalse
\ifadded\global\addedfalse
\iffirstadded\global\firstaddedfalse
\global\advance\corolcount by1 \fi
\the\corolcount\addedlett%
\else$\rm\spnum$\def\corolnum{$\rm\spnum$}\fi%
\else\global\advance\corolcount by1
\the\corolcount\fi}


\newcount\corolcount
\def\xcorol{Corollary}
\def\xtheorem{Theorem}
\def\xmaintheorem{Main Theorem}

\newif\ifthtitle

\let\saverparen)
\let\savelparen(
\def\rmparenl{{\rm(}}
\def\rmparenr{{\rm\/)}}
{
\catcode`(=13
\catcode`)=13
\gdef\makeparensRM{
\catcode`(=13
\catcode`)=13
\let(=\rmparenl
\let)=\rmparenr
\everymath{\let(\savelparen
\let)\saverparen}
\everydisplay{\let(\savelparen
\let)\saverparen}}
}

\medskipamount=8pt plus.1\baselineskip minus.05\baselineskip
\def\proclaim#1{\vskip-\lastskip
\def\one{#1}\ifx\one\xtheorem\global\corolcount=0\fi
\ifsection\global\sectionfalse\vskip-6pt\fi
\medskip
{\elevensc#1}%
\ifx\one\xmaintheorem\global\corolcount=0
\gdef\theoremnum{Main Theorem}\else%
\ifx\one\xcorol\ \allcorolnums\else\ \alltheoremnums\fi\fi%
\ifthtitle\ \global\thtitlefalse{\rm(\thethtitle)}\fi.%
\hskip.5pc\bgroup\makeparensRM
\it\ignorespaces}

\def\nonumproclaim#1{\vskip-\lastskip
\def\one{#1}\ifx\one\xtheorem\global\corolcount=0\fi
\ifsection\global\sectionfalse\vskip-6pt\fi
\medskip
{\elevensc#1}.\ifx\one\xmaintheorem\global\corolcount=0
\gdef\theoremnum{Main Theorem}\fi\hskip.5pc%
\bgroup\makeparensRM\it\ignorespaces}

\def\endproclaim{\egroup\medskip}


\def\xproof{Proof}
\def\xremark{Remark}
\def\xcase{Case}
\def\xsubcase{Subcase}
\def\xconjecture{Conjecture}
\def\xstep{Step}
\def\xof{of}

\def\deconstruct#1 #2 #3 #4 #5 @{\def\one{#1}\def\two{#2}\def\three{#3}%
\def\four{#4}%
\ifx\two\empty #1\else%
\ifx\one\xproof%
\ifx\two\xof%
  \ifx\three\xcorol Proof of Corollary \rm#4\else%
     \ifx\three\xtheorem Proof of Theorem \rm#4\else\xone\fi%
  \fi\fi%
\else\xone\fi\fi.}

\def\pickup#1 {\def\this{#1}%
\ifx\this\xproof\global\let\go\demoproof
\global\let\enddemo\endproof\else
\ifx\this\xremark\global\let\go\demoremark\else
\ifx\this\xcase\global\let\go\demostep\else
\ifx\this\xsubcase\global\let\go\demostep\else
\ifx\this\xconjecture\global\let\go\demostep\else
\ifx\this\xstep\global\let\go\demostep\else
\global\let\go\demoproof\fi\fi\fi\fi\fi\fi}

\def\demo#1{\vskip-\lastskip
\ifsection\global\sectionfalse\vskip-6pt\fi
\def\one{#1 }\def\two{#1*}%
\setbox0=\hbox{\expandafter\pickup\one}\expandafter\go\two}

\def\numbereddemo#1{\vskip-\lastskip
\ifsection\global\sectionfalse\vskip-6pt\fi
\def\two{#1*}%
\expandafter\demoremark\two}

\def\demoproof#1*{\medskip\def\xone{#1}
{\ignorespaces\it\expandafter\deconstruct\xone {} {} {} {} {} @%
\unskip\hskip6pt}\rm\ignorespaces}

\def\demoremark#1*{\medskip
{\it\ignorespaces#1\/} \alltheoremnums\unskip.\hskip1pc\rm\ignorespaces}

\def\demostep#1 #2*{\vskip4pt
{\it\ignorespaces#1\/} #2.\hskip1pc\rm\ignorespaces}

\def\enddemo{\medskip}

\def\endproof{\ifmathqed\global\mathqedfalse\medskip\else
\parfillskip=0pt~~\hfill\qed\medskip
\fi\global\parfillskip0pt plus 1fil\relax
\gdef\enddemo{\medskip}}

\def\qed{\vbox{\hrule\hbox{\vrule height6pt\hskip6pt\vrule}\hrule}}








\def\stripbs#1#2*{\def\one{#2}}
\def\newdef#1#2#{\expandafter\stripbs\string#1*
\defining{#1}{#2}}

\def\defining#1#2#3{\ifdefined{\one}{\message{<<<Sorry,
\string#1\space
is already defined. Please use a new name.>>>}}
{\expandafter\gdef#1#2{#3}}}

\def\emptyspace{ }
\def\nextthing{}
\def\newline{***}
\def\eatone#1{ }

\def\lookatspace#1{\ifcat\noexpand#1\ \else%
\gdef\nextthing{}\xdef\next{#1}%
\ifx\next\emptyspace%
\let\nextthing\emptyspace\else\ifx\next\newline%
\gdef\nextthing{\eatone}\fi\fi\fi\egroup\nextthing#1}

{\catcode`\^^M=\active%
\gdef\spacer{\bgroup\catcode`\^^M=\active%
\let^^M=\newline\obeyspaces\lookatspace}}

\def\ref#1{\seeifdefined{#1}#1\spacer}


\def\seeifdefined#1{\expandafter\stripbs\string#1*%
\crorefdefining{#1}}

\newif\ifcromessage
\global\cromessagetrue

\def\crorefdefining#1{\ifdefined{\one}{}
{\ifcromessage\global\cromessagefalse%
\message{\spaces\spaces\spaces\spaces\spaces\spaces\spaces}%
\message{<Undefined reference.}%
\message{Please TeX file once more to have accurate cross-references.>}%
\message{\spaces\spaces\spaces\spaces\spaces\spaces\spaces}\fi%
\expandafter\gdef#1{??}}}

\def\label#1#2*{\gdef\ctest{#2}\seeifdefined{#1}%
\ifx\empty\ctest\write\auxfile{\noexpand\def\noexpand#1{\tocstuff}}%
\else\def\ctemp{#2}%
\immediate\write\auxfile{\noexpand\def\noexpand#1{\ctemp}}\fi}

\def\ifdefined#1#2#3{\expandafter\ifx\csname#1\endcsname\relax#3\else#2\fi}




\def\articlecontents{
\vskip20pt\centerline{\bf Table of Contents}\everypar={}\vskip6pt
\bgroup \leftskip=3pc \parindent=-2pc
\def\item##1{\vskip1sp\indent\hbox to2pc{##1.\hfill}}}

\def\endcontents{\vskip1sp\leftskip=0pt\egroup}

\def\journalcontents{\vfill\eject
\def\currannalsline{\hfill}
\global\titletrue
\vglue3.5pc
\centerline{\tensc\hskip12pt TABLE OF CONTENTS}\everypar={}\vskip30pt
\bgroup \leftskip=34pt \rightskip=-12pt \parindent=-22pt
  \def\\ {\vskip1sp\noindent}
\def\pagenum##1{\unskip\parfillskip=0pt\dotfill##1\vskip1sp
\parfillskip=0pt plus 1fil\relax}
\def\name##1{{\tensc##1}}}


\institution{}
\onpages{0}{0}
\def\lastpage{}
\def\thetitle{Title ???}
\def\theauthors{Authors ???}
\def\thereceived{}
\def\therevised{}

\catcode`\@=12

\annalsline{June}{2003}

\startingpage{1}

%
%
%
%


\font\teneu=eufm10\font\seveneu=eufm7\font\fiveeu=eufm5
\font\tenlv=msbm11\font\sevenlv=msbm7\font\fivelv=msbm5

\newfam\eufam  \def\eu{\fam\eufam\teneu}
\textfont\eufam=\teneu \scriptfont\eufam=\seveneu 
                              \scriptscriptfont\eufam=\fiveeu

\newfam\lvfam  \def\lv{\fam\lvfam\tenlv}
\textfont\lvfam=\tenlv \scriptfont\lvfam=\sevenlv 
                              \scriptscriptfont\lvfam=\fivelv


\def\bfit#1{{\setbox0=\hbox{#1}%
   \setbox1=\hbox{\hbox to.1pt{}\copy0}\copy1\kern-\wd0%
   \copy1\kern-\wd0\copy1\kern-\wd0\copy1\kern-\wd0\copy1\kern-\wd0\copy1}}

\def\birmap{\hhb2\hbox{-\hhb{1}-\hhb{1}-}\hhb{-1}%
\lower-1pt\hbox{$\scriptscriptstyle>$}\hhb2}

\def\({\big(}  \def\){\big)}  \def\nix{{\phantom{|}}}
\def\aussag#1.{\copy#1}  \def\AUSSAG\nmc#1.{\copy#1}

 \def\bsn{\bigskip\noindent}


\def\horrdwndiag#1#2#3#4#5#6#7#8{
              $$\def\normalbaselines{
                    \baselineskip20pt\lineskip3pt\lineskiplimit3pt}
              \matrix{#1&\horr{#5}&#2&\cr
                      \dwn{#6}&&\dwn{#7}&\cr
                      #3&\horr{#8}&#4&\cr}
                              $$}


\def\matrixbaselines{
   \def\normalbaselines{\baselineskip20pt\lineskip3pt\lineskiplimit3pt}}

    \def\dwn#1{\Big\downarrow\rlap{$\vcenter{\hbox{$\scriptstyle{#1}$}}$}}

     \def\equnum#1.{\global\advance\NO by 1%
          {\global\setbox#1=\hbox{\rm\number\NOCAP.\number\NO}}}
 

  
  \def\hhb#1{\hbox to#1pt{}}
  \def\hor#1{\smash
     {\mathop{{\lgrghtar}}\limits^{\lower2pt\hbox{$\scriptstyle{#1}$}}}}
  \def\horr#1{\smash
     {\mathop{{\lglgrghtar}}\limits^{\lower2pt\hbox{$\scriptstyle{#1}$}}}}
  \def\horrr#1{\smash
     {\mathop{{\lglglgrghtar}}\limits^{\lower2pt\hbox{$\scriptstyle{#1}$}}}}

\def\ilim#1{\hbox to14pt{lim\kern-14pt\lower4.5pt\hbox{$\scriptstyle
  \longrightarrow$}\kern-8pt\lower8.5pt\hbox{$\scriptstyle{#1}$}}\hhb{3}}

   

\def\lgrghtar{{\hhb2{\relbar\joinrel\rightarrow}\hhb2}}
\def\lglgrghtar{{\hhb1{\relbar\joinrel\relbar\joinrel\rightarrow}\hhb1}}
\def\lglglgrghtar{{\hhb1{\relbar\joinrel\relbar\joinrel%
\relbar\joinrel\rightarrow}\hhb1}}

    \def\msn{\medskip\noindent}

\def\nmnm#1{{\tensc #1}}
\def\nix{{\phantom{|}}} 



\def\plim#1{\hbox to14pt{lim\kern-14pt\lower4.5pt\hbox{$\scriptstyle
  \longleftarrow$}\kern-8pt\lower8.5pt\hbox{$\scriptstyle{#1}$}}\hhb{3}}


\def\ratmap{\hhb2\hbox{-\hhb{1}-\hhb{1}-}\hhb{-1}%
\lower-1pt\hbox{$\scriptscriptstyle>$}\hhb2}

\def\sag#1{\hhb1{\hbox to#1mm{\rightarrowfill}}\hhb1}

    \def\ssn{\smallskip\noindent}
    
    \def\sdp{\hbox to10pt{\hss\hbox{\mathsurround=0pt$\times$\kern-1.6pt
       \hbox{\vrule height5.2pt width.6pt}\hbox to1.6pt{}}\hss}}


 \def\td{{\rm td}}


\def\vid{{\hhb{-3.25}\not\hhb{-2.5}\lower-1.25pt\hbox{\mathsurround=0pt$
                              \scriptscriptstyle\bigcirc$}}}

\def\vid{{\lower-1pt\hbox{/}\kern-7pt{\hbox{O}}}}


\def\abstr{pre-divisorial}

\def\cycellprefrm{Galois format\-ion}


\def\abab#1{#1^{\ell,{\rm ab}}}

\def\bfv{{\bf v}}  
  
\def\ker{{\rm ker}}

\def\calz{{\cal Z}}
\def\cd{{\rm cd}}
\def\chr{{\rm char}}
\def\clg{{\cal G}}

\def\divdiv#1#2{\oplus\, v_a\LL}
\def\divinrt#1{{\eu div.inr}(#1)}

\def\fell#1{{#1}^\ell}

\def\Gell#1{G^\ell_{#1}}
\def\Gellab#1{G^{\ell,\rm ab}_{#1}}
\def\ggg{G} 
\def\ggv{G_v}

\def\Gprm#1{G^{\,'}_{\!#1}}

\def\hhhom#1#2{{\rm Hom}\big(#1,#2\big)}
\def\hhmell#1{{\rm Hom}\big(#1,{\lv Z}_\ell\big)}
\def\hhone#1#2{{\rm H}^1\Big(#1,#2\Big)}

\def\inrt#1{{\eu inr}(#1)}

\def\ker{{\rm ker}}
\def\kpel{K^\times_\pel}
\def\kpkp{{\hhb1{\eu k}\hhb1}}
\def\KK{{{\eu K}\hhb1}}

\def\kxzk{{\cal K}}

\def\LX{{\cal L}}
\def\lxzl{\LX_{(\ell)}}
\def\lxlxzl#1{\LX_{#1,(\ell)}}

\def\mardv{{\rm Div}}

\def\pel{{(\ell)}}

\def\piellab{{\pi_1^{\ell,ab}}}
\def\prm#1{#1'}

\def\rszv{{\rm res}_{\zzv}}

\def\shra#1{\!\!\!\hookrightarrow\!\!\!}
\def\shor#1{\!\!\!\!\!\hor{#1}\!\!\!\!\!}

\def\shorrr#1{\!\!\!\horrr{#1}\!\!\!}
\def\shto{\!\!\!\!\!\to\!\!\!\!\!}

\def\Tellab#1{T^{\ell,\rm ab}_{#1}}
\def\Tell#1{T^{\ell}_{#1}}
\def\ttab{{T^{\rm ab}}}

\def\ttt{T}
\def\ttv{T_v}

\def\ttw{T_w}

\def\vurl{very unruly}

\def\whlx{\widehat\LX}
\def\whph{{\hat\varphi}}

\def\Zell{{\lv Z}_\ell}

\def\zetell{{\lv Z}_{(\ell)}}
\def\zzv{Z_v}



\def\DATUM{June, 2003}

\title{Pro-$\ell$ birational anabelian geometry \\
                     over \\
          algebraically closed fields I}

\shorttitle{Birational anabelian geometry}

\acknowledgements{ \ Ideas used in this work originate
in part from my visit to the IAS Princeton in 1994 and 1996.
I would like to thank Pierre Deligne for several fruitful 
discussions and suggestions, see also [P2], [P3]. The main 
technical details of the present manuscript were solved 
during my visit to the IAS Princeton in 2002. During this 
last visit I was partially supported by the NSF grant DMS 
9729992.}

\author{Florian Pop}

\institutions{Mathematisches Institut,
              Universit\"at Bonn}


\intro

Let $\ell$ be a fixed rational prime number. For every  
field $K$ of positive characteristic $\neq\ell$ containing 
the $\ell^{\rm th}$ roots of unity, let 
$\fell K|K$ be the maximal pro-$\ell$ Galois extension, 
and $\Gell K={\rm Gal}(\fell K|K)$ its Galois group. 

\ssn

The aim of this paper is to prove a  
{\it geometric pro-$\ell$ version \/} of Grothen\-dieck's 
birational anabelian conjecture over algebraic closures
of finite fields as follows: 

\nonumproclaim{Theorem} Let ${\cal F}$ be the category 
of all function fields $K|k$, with $\td(K|k)>1$ and $k$ 
an algebraic closure of some finite field.  Then there
exists a group theoretic recipe by which we can recover 
every field $K\in{\cal F}$ from $\Gell K$, up to a pure
inseparable extensions. This recipe is invariant under
profinite group isomorphisms. In particular, if $K$ and 
$L$ are in ${\cal F}$, there exists a canonical bijection
$$
{\rm Isom}^{\rm i}\,(L^{\rm i},K^{\rm i})\hor{}
      {\rm Out}\,(\Gell{K},\Gell{L})\,,
$$
where $(\hhb3)^{\rm i}$ denotes pure inseparable closure, 
and ${\rm Isom}^{\rm i}$ means up to Frobenius twists, and
${\rm Out}$ denotes outer isomorphisms of profinite groups.
\ssn
\indent
Equivalently, if $\Phi:\Gell{K}\to\Gell{L}$ is an isomorphism
of profinite groups, then up to a Frobenius twist, there 
exists a unique field isomorphism $\phi:\fell{L^{\rm i}}\to
\fell{K^{\rm i}}$ such that $\Phi(g)=\phi^{-1}\,g\,\phi$ for 
all $g\in\Gell{K}$. In particular, $\phi(L^{\rm i})=K^{\rm i}$.

\endproclaim

\nonumproclaim{Remarks} {\rm

\ssn

I first want to mention that this manuscript is motivated
by the {\it program initiated by\/} \nmnm{Bogomolov}~[Bo].
And the result above completes that program in the case 
of function fields over algebraic closures of finite fields.
At least at the level of rough ideas, the paper [Bo] was 
quite inspiring for me... and at a first glance there is 
a lot of similarity of what we are doing here with loc.cit..\  
Without going into details, we hope nevertheless that the 
advised reader will realize the essential differences. 
Compare also with \nmnm{Bogomolov--Tschinkel}~[B--T\hhb{1}2].

\ssn

1) The Theorem above implies corresponding ``full Galois
assertions'', i.e., the corresponding assertions for the
full Galois group $G_K$ in stead of the quotient $\Gell K$.

\ssn

Second, the corresponding questions concerning the Galois 
characterisation of finitely generated fields $K$ with 
$\td(K)>1$ in positive characteristic, see [P2], [P3], [P4]
can be reduced to the Main Theorem above. Thus the above 
result {\it generalises\/} in a non-trivial way the 
celebrated results by \nmnm{Neukirch, Ikeda, Iwasawa, Uchida} 
concerning the Galois characterisation global fields. 

\ssn

2) The above result {\it cannot\/} be true in the case 
$\td(K|k)=1$, as in this case $\Gell K$ is a pro-$\ell$ 
free group on countably many generators, thus its structure
does not depend on $K$ at all. 

\ssn

Further, because of the reason just mentioned above, one
cannot expect a {\it Hom-form\/} type result in the above
context, as proved by \nmnm{Mochizuki} over sub-$p$-adic 
fields. Indeed, for any given function field $K|k$ there
exist ``many'' surjective homomorphisms of profinite groups
$G_K\to G_{k(t)}$, given say by field $k$-embeddings 
$k(t)\hookrightarrow K$ such that $k(t)$ is relatively
closed in $K$; and further, since $G_{k(t)}$ is profinite
free, see \nmnm{Harbater}~[Ha] and \nmnm{Pop}~[Pop],
there are ``many'' surjective projections
$G_{k(t)}\to G_K$. Thus finally, there are many surjective
group homomorphisms $G_K\to G_K$ which do not arise in a
geometric way...  
  
\msn

3) Finally, we remark that in stead of working with the 
full pro-$\ell$ Galois group $\Gell K$, one 
could work as well with truncations of it, e.g., with ones 
coming from the central series. This is doable, but the 
resulting assertions are quite technical, and at the 
moment maybe too complicated in order to be interesting... 
}
\endproclaim


\ssn

{\it Rough idea of Proof\/} (Comp.\ with \nmnm{Bogomolov} [Bo])

\msn

Let $\KK|\kpkp$ be an extension of fields of characteristic 
$\neq\ell$, and suppose that $\kpkp$ is algebraically closed. 
We consider algebraic extensions $\prm\KK|\KK$ with the 
following properties:  

\ssn

{i)} $\prm\KK|\KK$ is a pro-$\ell$ Galois extension.

\ssn

{ii)} There exists a Galois sub-extension $\KK_1$
of $\prm\KK|\KK$ such that $\prm\KK|\KK_1$ is a maximal 
pro-$\ell$ Abelian extension of $\KK_1$. 

\ssn

We remark that the extension $\fell K|K$ as introduced before 
the Main Theorem as well as $K^{\ell,{\rm ab}}|K$ satisfies the 
conditions i), ii) above. 

\ssn

To fix notations, let $\Gprm\KK={\rm Gal}(\prm\KK|\KK)$,
and ${\lv T}_\ell\hor{\imath}\Zell$ a fixed 
identification as $\Gprm\KK$-modules. By Kummer Theory 
there is a functorial isomorphism 
$$
\widehat\delta_\KK:\widehat\KK\hor{}
      \hhone{\KK}{{\lv T}_\ell}
           \hor{\imath}\hhhom{\Gprm\KK}{\Zell}\,,
\leqno{\indent(*)}
$$
where $\widehat\KK$ denotes the $\ell$-adic completion 
of the multiplicative group $\KK^\times$ of $\KK$. Since
$\kpkp^\times$ is divisible, and $\KK^\times\!/\kpkp^\times$ 
is a free Abelian group, the $\ell$-adic completion 
homomorphism defines an exact sequence 
$$
1\to\kpkp^\times\hookrightarrow\KK^\times
                \hor{\jmath_\KK}\widehat\KK\,. 
$$

Denote ${\cal P}(\KK)=\KK^\times\!/\kpkp^\times$ inside 
$\widehat\KK$, and view ${\cal P}(\KK)$ as the projectivization 
of the (infinite) dimensional $\kpkp$-vector space $(\KK,+)$.
{\bf Suppose} that we have a Galois theoretic recipe in order 
to detect: First the {\it image of\/} ${\cal P}(\KK)$ inside 
the Galois theoretically ``known'' $\widehat\KK$. 
Second, the {\it projective lines in\/} ${\cal P}(\KK)$. Remark
that the multiplication $l_x$ by any non-zero element $x\in\KK$ 
defines an ``automorphisms'' of ${\cal P}(\KK)$ which respects 
``co-lineations''. This automorphism is noting but the
translation $l_{\jmath_\KK(x)}$ in ${\cal P}(\KK)$,
this time viewed again as the multiplicative group  
${\cal P}(\KK)=\KK^\times\!/\kpkp^\times$. Now {\bf using} 
the {\it Fundamental Theorem of projective geometry\/}, see 
e.g.\ \nmnm{Artin}~[A], it follows that the additive structure 
of $\KK$ can be deduced from the knowledge of all the 
projective lines in ${\cal P}(\KK)$. And finally, since 
the ``multiplications'' $l_{\jmath_\KK(x)}$ do respect 
this structure, we finally deduce from this the field 
structure of $\KK$, as well the field extension $\KK|\kpkp$. 

\msn

Now coming back to the case of a function fields $K|k$
as in the Theorem above, the problems we have to tackle 
are the following: Give {\it Galois theoretic recipes\/} 
in order to detect: 

\msn

I) \ $K_{(\ell)}:=K^\times\otimes{\lv Z}_{(\ell)}$ 
{\it inside\/} $\widehat K$. 

\msn

This almost answers the question about detecting 
${\cal P}(K)$: We know namely its ${\lv Z}_{(\ell)}$-version, 
but not ${\cal P}(K)$ itself. 

\msn

II) \ ${\cal P}(K)$ {\it and its projective lines.\/}

\msn

III) \ Using the functoriality of the construction, show
that the recipe for detecting $(K,+,\cdot)$ is invariant
under isomorphisms. 

\bsn
{\it Organization of the paper\/}

\msn
In the first part we put together the necessary tools for 
the proof as follows:

\ssn

a) \ First, a ``pro-$\ell$ Local theory'' similar to
the Local Theory form [P1], etc.. It relies on results 
by \nmnm{Ware}~[W], \nmnm{Koenigsmann}~[Ko], see also 
\nmnm{Bogomolov--Tschinkel} [B--T\hhb{1}1], 
\nmnm{Efrat}~[Ef], \nmnm{Engeler--Koenigsmann} [E--K]. 
The aim of this local theory is to recover in a 
{\it functorial way\/} from $\Gell K$ the set 
${\eu D}_K$ of all the Zariski prime divisors of 
$K|k$, where $K$ is any function field over an 
algebraic closure $k$ of a finite field such that
$\td(K|k)>1$.   

\ssn

b) \ Part of the local theory is to show that in the 
context from above, the whole inertia structure which 
is significant for us is encoded in $\Gell K$. 

\ssn

c) \ Third, a more technical result by which we recover 
the ``geometric sets of prime divisors'' of function 
fields $K|k$ as above. By definition, a geometric set 
of prime divisors of $K$ is the set of Zariski prime 
divisors ${\eu D}_X\subset{\eu D}_K$ defined by the 
Weil prime divisors of some quasi-projective normal model 
$X$ of $K$. This result itself relies on \nmnm{de~Jong}'s 
{\it theory of alterations\/}~[J].

\msn
In the second part of the manuscript, we give a 
simplified version of a part of the ``abstract non-sense'' 
from Part~II of [P4] (which follows a suggestion by 
\nmnm{Deligne} [D2]), reminding one in some sense 
of the abstract class field theory. We define so called 
{\it \abstr\ \cycellprefrm s\/}. The aim of this 
theory is to lay an axiomatic strategy for the proof of 
the main result. 

The interesting example of \abstr\ \cycellprefrm s are 
the {\it geometric\/} \cycellprefrm s, which arise from 
geometry (and arithmetic).

\ssn

A very basic results here is
Proposition~3.18, which shows that in the case $K|k$ is 
a function field with $\td(K|k)>1$ and $k$ an algebraic 
closure of a finite field, the geometric \cycellprefrm s\ 
on $G=G^{\ell,{\rm ab}}_K$ are group theoretically encoded 
in $\Gell K$.

\ssn

Finally in the last Section we prove the main Theorem
announced above. The main tool here is Proposition~4.1,
which gives a Galois characterization of the ``rational
projections''. 

\msn

{\bf Thanks:} I would like to thank several people who showed
interest in this work, both for criticism and suggestions,
and for careful reading. My special thanks go to David Harbater, 
Jochen Koenigsmann, Pierre Lochak, Hiroaki Nakamura, Frans
Oort, Tam\'as Szamuely, Jakob Stix, Akio Tamagawa, etc..
%
%
%
%
%
\section{Local theory and inertia elements}
In this section we first recall the main facts concerning 
the Local theory from [P1], but in a pro-$\ell$ setting. Our 
first aim is to give to give a group theoretic recipe for finding 
information about the {\it space of Zariski prime divisors\/}  
in a functorial way.
 
\ssn

Let $\KK|\kpkp$ be a function field over some base field 
$\kpkp$, and $\td(\KK|\kpkp)=d>0$. We consider the family 
of all models $X_i\to\kpkp$ of $\KK|\kpkp$ such that the 
structure sheaf of $X_i$ is a sheaf of sub-rings of $\KK$ 
with inclusions as structure morphisms. On the family of 
all the $X_i$'s there exists a naturally defined domination 
relation as follows: $X_j\geq X_i$ if  there exists a 
surjective $\kpkp$-morphism $\varphi_{ji}:X_j\to X_i$ which 
at the structure sheaf level is defined by inclusions. Let 
${\eu Proj}_\KK$ be the subfamily of all projective, normal 
models of the function field $\KK$. The following is well 
known, see e.g.\ \nmnm{Zariski--Samuel} [Z--S], Ch.VI, 
especially \S17: 

\ssn

{$\bullet$} Every complete model is dominated by some 
$X_i\in{\eu Proj}_\KK$ (Chow Lemma).
%
%

\ssn

{$\bullet$} The set ${\eu Proj}_\KK$ is increasingly 
filtered with respect to $\geq$, hence it is a surjective 
projective system.

\ssn

{$\bullet$} Denote ${\eu R}_\KK=\plim{i}X_i$ as  
topological spaces. We will call ${\eu R}_\KK$ the 
{\it Riemann space\/} of $\KK$. The points of 
${\eu R}_\KK$ are in bijection with the space of all 
valuation $\kpkp$-rings of $\KK$. For $v=(x_i)_i$ in 
${\eu R}_\KK$ one has: $x_i$ is the centre of $v$ on $X_i$ in 
the usual sense, and ${\cal O}_v=\cup_i{\cal O}_{X_i,x_i}$.
%
%

\proclaim{Fact/Definition} {\rm 
Using e.g.\ [BOU], Ch.IV, \S3, one shows that for 
a point $v=(x_i)_i$ in ${\eu R}_\KK$ the following 
conditions are equivalent:

\ssn

{i)} For $i$ sufficiently large, $x_i$ has co-dimension 1, 
or equivalently, $x_i$ is the generic point of a prime Weil 
divisor of $X_i$. Hence $v$ is the discrete $\kpkp$-valuation 
of $\KK$ with valuation ring ${\cal O}_{X_i,x_i}$. 

\ssn

{ii)} $\td(\KK v\,|\,\kpkp)=\td(\KK|\kpkp)-1$.

\ssn

We will say that a point $v=(x_i)_i$ in ${\eu R}_\KK$ 
satisfying the above equivalent conditions is a {\it Zariski 
prime divisor\/} of $\KK|\kpkp$. We denote the space 
of all Zariski prime divisors of $\KK|\kpkp$ by 
${\eu D}_\KK\subset{\eu R}_\KK$. One has: The space 
${\eu D}_\KK\subset{\eu R}_\KK$ of all Zariski prime 
divisors of $\KK$ is the union of the spaces ${\eu D}_i$ 
of prime Weil divisors of all the models $X_i$ of 
$\KK|\kpkp$ in ${\eu Proj}_\KK$ (if we identify every 
prime Weil divisor with the discrete valuation on $\KK$ 
it defines). 

\ssn

More generally, let $\KK'|\KK$ be an arbitrary algebraic 
extension, and $v'$ a valuation on $\KK'$. Then we will 
say that $v'$ is a Zariski prime divisor (of $\KK$ or of
$\KK'$), if its restriction $v$ to $\KK$ is a Zariski prime 
divisor of $\KK|\kpkp$.
}
\endproclaim   

A) \ {\it On the decomposition group\/}

\ssn

Let $\KK'|\KK$ be some Galois field extension as in the 
Introduction, {\it Rough idea of proof,\/} and $v$ is a 
valuation on $\KK'$ which is trivial on $\kpkp$. Let
$Z_v$, $T_v$, and $V_v$ be respectively the decomposition 
group, the inertia group, and the ramification group of 
$v$ in $\Gprm\KK:={\rm Gal}(\KK'|\KK)$. Since 
$\ell\neq{\rm char}(\kpkp)$,
and $v$ is trivial on $\kpkp$, $V_v$ is trivial. We denote 
by $\KK^T$, and $\KK^Z$ the corresponding fixed fields.

\ssn

\proclaim{Fact} \
{\rm The following are well known facts from Hilbert 
decomposition, and/or ramification theory for general 
valuations: $\KK'v\,|\,\KK v$ is a Galois field extension, 
which also satisfies the properties i), ii), from loc.cit..
Let $\Gprm{\KK v}={\rm Gal}\,(\prm\KK v\,|\,\KK v)$ be its
Galois group. One has a canonical exact sequence 
$$
1\to T_v\to Z_v\hor{\pi_v} G_v:=\Gprm{\KK v}\to1\>.
$$
Moreover, $v(\KK^T)=v(\KK^Z)=v\KK$, and $\KK v=\KK^{\!Z}v$; 
and further $v\prm\KK$ is an $\ell$-divisible hull of $v\KK$.

\ssn

Finally, there exists a pairing 
$\Psi_{\KK'}:T_v\times v\KK'\to\mu_{\KK'v}$,
$(g,\,vx)\mapsto(gx/x)v$, and the following hold: The left 
kernel of $\Psi_{\KK'}$ is exactly $V_v=\{1\}$, and the right 
kernel of $\Psi_{\KK'}$ is $v\KK$. In particular, $T_v$ is 
Abelian. Further, $\Psi$ is compatible with the action of 
$G_v$ on $T_v$ via $\pi_v$. In particular, this action 
is via the cyclotomic character, thus trivial, as $\KK v$ 
contains by hypothesis the algebraically closed field $\kpkp$.
Therefore, we finally have: $Z_v\cong T_v\times G_v\,$.

\ssn

1) Let ${\cal B}=(vx_i)_i$ be an ${\lv F}_\ell$-basis of 
$v\KK\,/\,\ell$. Then $T_v\cong{\lv Z}_\ell^{\cal B}$, and 
finally
$$ 
Z_v\,\cong\,T_v\times\Gprm{\KK v}\,\cong\,
           {\lv Z}_\ell^{\cal B}\times\Gprm{\KK v}\,
$$ 
as profinite groups. In particular, 
$\cd(Z_v)=\cd(\Gprm{\KK v})+|{\cal B}|$, 
where $\cd$ denotes the cohomological dimension.

\ssn

2) \ Thus finally there exist canonical isomorphisms
$\theta^v:\widehat{v\KK}\to\hhhom{T_v}{{\lv Z}_\ell}$
\ and \ $\theta_v:T_v\to\hhhom{v\KK}{{\lv Z}_\ell}\,$.

\ssn

3) \ Recalling the notations and remarks from 
Introduction, one gets a commutative diagram of the form:

\horrdwndiag{\;\KK^\times}{v\KK}
 {\hhhom{\Gprm\KK}{{\lv Z}_\ell}}{\hhhom{T_v}{{\lv Z}_\ell}}
   {v}{\jmath_\KK}{\theta^v}{\jmath^v}

\ssn

4) \ Let $U_v\subset\KK^\times$ denote the $v$-units, 
and $U_v^1=1+{\eu m}_v\subset U_v$ the principal $v$-units. 
Then $U_v^1$ becomes $\ell$-divisible in $\KK^Z$, and 
$(\KK v)^\times=U_v/U_v^1$. We have a commutative diagram 
of the form:

\horrdwndiag{U_v}{\KK v}{\hhhom{\Gprm\KK}{{\lv Z}_\ell}}
     {\hhhom{\Gprm{\KK v}}{{\lv Z}_\ell}}
   {}{\jmath_\KK}{\jmath_{\KK v}}{\jmath_v}

\ssn

B) \ {\it Recovering ${\eu D}_K$ from $\Gell K$\/}

\proclaim{Definition/Remark} 
{\rm In the notations from Introduction and the above ones, 
let $K|k$ be a function field with $\td(K|k)>0$ and $k$ an
algebraic closure of a finite field.
 
\ssn

1) The decomposition group $Z_v\subset\Gell K$ of some 
divisorial valuation $v\in{\eu D}_{\fell K}$ is called 
a {\it divisorial subgroup.\/} 

\ssn

2) A subgroup $Z\subset\Gell K$ which is isomorphic 
to a divisorial subgroup of some function field of 
transcendence degree equal to $\td(K|k)$ is called a 
{\it divisorial like subgroup.\/}
}
\endproclaim 
\proclaim{Proposition} 
Let $K|k$ be a function field with $\td(K|k)>1$ and 
$k$ an algebraic closure of a finite field. Let 
$Z_v\subset\Gell K$ be a divisorial subgroup, say 
defined by a divisorial valuation $v$ on $\fell K$. Let 
$T_v$ be the inertia group of $v$. Then the following 
hold:

\ssn

{\rm(1)} $Z_v$ is self-normalising in $\Gell{K}$. 
Further, if $Z_{v'}\neq Z_v$ is another divisorial 
subgroup, then one has: $Z_{v'}\cap Z_v=1$.

\ssn

{\rm(2)} $T_v\cong\Zell$ as a $\Gell{K v}$-module. 
Further $T_v$ is the unique maximal pro-$\ell$ Abelian 
normal subgroup of $Z_v$. And finally
$
Z_v\cong T_v\times\Gell{Kv}\cong\Zell\times\Gell{Kv}.
$
\endproclaim 
\label\propdecomgr\theoremnum*
\demo{Proof}
To (1): Both assertions follow using a result of 
\nmnm{F.\ K.\ Schmidt}, see e.g., \nmnm{Pop}~[P1], 
Proposition~1.3; see also the proof of Proposition~1.14 
from loc.cit.. Concerning (2), the only non-obvious 
part is the fact that given any function field $K'=Kv$ 
over $k$, the Galois group $\Gell{K'}$ has no non-trivial 
Abelian normal subgroups. This follows from the 
Hilbertianity of $K v$. 
\enddemo

\ssn

A first content of the local theory is that ``morally'' 
the converse of the above Proposition is also true, i.e., 
if $Z\subset\Gell K$ is a divisorial like subgroup, 
then it comes from a Zariski prime divisor.    
\proclaim{Proposition} Let $K|k$ be a function field 
with $\td(K|k)=d>1$ and $k$ an algebraic closure of 
a finite field. Then one has:

\ssn

{\rm(1)} For every divisorial like subgroup $Z$ of 
$\Gell K$ there exists a unique divisorial valuation 
$v$ of $\fell K$ such that $Z\subset Z_v$ and 
$\chr(K v)\neq\ell$. 

\ssn

{\rm(2)} Moreover, if $T\subset Z$ is the unique maximal
Abelian normal subgroup of $Z$, then $T=Z\cap T_v$, where 
$T_v$ is the inertia subgroup of $Z_v$.

\ssn

Therefore, the space ${\eu D}_{\fell K}$ of all 
Zariski prime divisors of $\fell K$ is in bijection 
with the divisorial subgroups of $\Gell K$. This 
bijection is given by $v\mapsto Z_v$. 
\endproclaim
\label\chardivsubg\theoremnum*

\demo{Proof} The main step in the proof is the following 
$\ell$-Lemma below, which replaces the $q$-Lemma from 
[P1], Local theory. After having the $\ell$-Lemma, the 
remaining steps in the proof are similar to (but easier
than) the ones from the proof of Theorem~1.16 from [P1],
Local theory. We will skip the remaining details.
\enddemo
\nonumproclaim{The $\ell$-Lemma {\rm(revisited)}} 
In the context of Proposition~\chardivsubg, let 
$Z_0\cong\Zell^d$ be a closed subgroup of $\Gell K$. 
Then there exists a valuation $v_0$ of $\fell K$ such
that $Z_0\subset Z_{v_0}$, and ${\rm char}(Kv)\neq\ell$.
\endproclaim

There are several ways to prove the $\ell$-Lemma above: 
First, one could develop the corresponding model theoretic 
machinery, and proceed as in [P1], Local theory. Second, 
one can apply the results from \nmnm{Ware}~[W], and 
\nmnm{Koenigsmann}~[Ko], more precisely 
\nmnm{Enge\-ler--Koenigsmann}~[E--K]; see also \nmnm{Efrat}~[Ef].
Or third, apply \nmnm{Bogomolov--Tschinkel} [B--T\hhb{1}1].

\msn

C) \ {\it Inertia elements of\/} $\Gell K$

\ssn

Let $K|k$ be a function field with $\td(K|k)>1$ and 
$k$ an algebraic closure of a finite field. In this 
subsection we give a description of an important class 
of inertia elements in $\Gell K$ via divisorial 
inertia elements. 

\proclaim{Definition.} \ 
{\rm Let $K'|K$ be a Galois extension, and let denote
$G':={\rm Gal}(K'|K)$. For a valuation $v$ of 
$K$ and a prolongation $v'$ of $v$ to $K'$, let 
$T_{v'}\subset Z_{v'}$ be its inertia, respectively 
decomposition group in $G'$. 

\ssn

1) An element $g$ of $G'$ is called a $v$-inertia element, 
if $g\in T_{v'}$, for some prolongation $v'|v$. In general, 
an element $g$ of $G'$ is called an inertia element, if 
there exists $v$ such that $g$ is a $v$-inertia element. 

\ssn

Denote by $\inrt{K'}$ is the set of all the inertia 
elements in $G'$.

\ssn

2) An inertia element $g$ of $G'$ is called divisorial 
inertia element, if $g$ is a $v$-inertia element to some 
divisorial valuation $v$ of $\KK$. 

\ssn

Denote by $\divinrt{K'}$ is the set of all divisorial 
inertia elements of $G'$. 
}
\endproclaim

\proclaim{Remarks} \ 
{\rm 
Let $K_0$ be a finitely generated infinite field, e.g.\ 
a global field. If ${\rm char}(K_0)>0$, then denote by
$K=K_0k$ a maximal constant extension of $K_0$. (Thus $k$
is an algebraic extension of the constants of $K_0$, 
and $K|k$ is a function field as above.) As above, let 
$K_0'|K_0$ be some Galois extension, and $G_0'$ its Galois
group. Thus if $K_0$ has positive characteristic, then 
$K':=K_0'k$ is a Galois extension of $K$, and let 
$G'\subset G_0'$ is a closed subgroup of $G_0'$.    

\ssn

{\rm a)} Let ${\cal X}\to{\lv Z}$ be a model of $K_0$.
Recall, that an element $\sigma$ of ${\rm Gal}(K_0'|K_0)$ 
is called a Frobenius element of $K_0'$ (over ${\cal X}$), 
if there exits a regular closed point $x\in{\cal X}$, 
and a decomposition groups $D_x\subset{\rm Gal}(K_0'|K_0)$ 
over $x$, such that $\sigma\in D_x$ is in the pre-image 
of the Frobenius at $x$ in $D_x$. Let 
${\eu Frob}(K_0')$ be the set of all the Frobenius elements 
in ${\rm Gal}(K_0'|K_0)$ (for the several models of $K_0$ 
as above). By the {\it Chebotarev Density Theorem,\/} 
${\eu Frob}(K_0')$ is dense in ${\rm Gal}(K'_0|K_0)$. 

\ssn

{\rm b)} In contrast to the situation above, the set 
$\divinrt{K_0'}$ is in general not dense in 
${\rm Gal}(K_0'|K_0)$. Namely, if $T'$ is the closed 
subgroup of ${\rm Gal}(K_0'|K_0)$ generated by all the divisorial
inertia elements, and if ${\cal X}$ is any {\it complete 
regular\/} model of $K_0$ (if any such models do exist), 
then ${\rm Gal}(K_0'|K_0)\,/\,T'\cong\pi_1({\cal X})$. 

\ssn

This is a totally different situation than that of the 
Frobenius elements of global fields. But a much more 
stronger assertion holds, as follows. First, recall that
for every valuation $v$ on $K$ one has the following 
fundamental inequality, see e.g.\ [BOU], Ch.IV, \S10, 3:
$$
\td(K|k)\geq {\rm rr}(vK)+\td(Kv|k)
$$
where ${\rm rr}(vK)$ the rational rank of the group $vK$.
A valuation of $K$ is called {\it defectless,\/} if the
above inequality is an equality. 
}
\endproclaim

\nonumproclaim{Fact} \ {\rm 
The most prominent defectless valuations are the 
``iterations'' of Zariski prime divisors as follows: 
Let $\delta\leq d=\td(K|k)$ be a fixed positive integer.
Define inductively a sequence of function fields $K_k$ 
endowed with Zariski prime divisors $v_k$ as follows:
$v_1$ is a Zariski prime divisor of $K_1:=K$, and set
$K_2:=K_1 v_1$, etc.. Thus inductively: $v_k$ is a 
Zariski prime divisor of $K_k:=K_{k-1} v_{k-1}$ for all 
$k\leq\delta$. Finally set $v=v_\delta\circ\dots\circ v_1$
as a valuation on $K$. Then $v$ is a ``generalised''
discrete valuation of $K$ with $v K\cong {\lv Z}^\delta$
lexicographically ordered, and 
$\td(Kv|k)=\td(K_\delta v_\delta|k)=\td(K|k)-\delta$.
Therefore $v$ is defectless.
}
\endproclaim

\proclaim{Theorem} \ 
In the notations from the Definition above, let ${K}|{k}$
be a function field with ${k}$ an algebraically closed 
base field, $\prm{K}|{K}$ be a Galois sub-extension of 
$\fell{K}\,|\,{K}$. Then the following hold: 

\ssn

{\rm(1)} \ $\,\inrt{\prm{K}}$ is closed in 
${\rm Gal}(\prm{K}|{K})$.

\ssn

{\rm(2)} \ The closure of $\,\divinrt{\prm{K}}$ in 
${\rm Gal}(\prm{K}|{K})$ contains all the inertia elements 
at defectless valuations $v$ which are trivial on ${k}$. 
\endproclaim
\label\inertiaelem\theoremnum*

\demo{Proof} In both cases we can suppose that 
${\prm{K}}=\fell{K}$, as both closed subsets and 
inertia elements are closed under the projection 
$\Gell{K}\to{\rm Gal}(\prm{K}|{K})$.

\ssn

To (1): \ Let $g\neq 1$ lie in the closure of 
$\inrt{\fell K}$. Equivalently, for every  finite Galois 
sub-extension ${K}_i|{K}$ of $\Gell{K}\,|\,{K}$, there do 
exist:

{i)} A quasi divisorial valuation $v_i$ on $\fell{K}$.

\ssn 

{ii)} Some $g_i$ in the inertia group $T_i:=T_{v_i}$ 
of $v_i$ in $\Gell{K}$, such that $g$ and $g_i$ have the 
same restriction to ${K}_i$.

\ssn
 
We will show that $g$ is a $v$-inertia element of $\Gell{K}$. 
First, in the notations from above, let $\ell^{n_i}$ be the 
common order of (the restriction of) $g$ and $g_i$ on 
${K}_i$. Further, let ${K}_0$ be the fixed field of $g$ 
in $\fell{K}$, thus $\fell{K}\,|\,{K}_0$ has Galois group 
generated by $g$. By Kummer Theory, there exits some 
$x\in{K}_0$ which is not an $\ell^{\rm th}$ power in 
${K}_0^\times$; and for every $\alpha_n\in\fell{K}$ such 
that $\alpha^{\ell^n}=x$, and ${K}_n:={K}_0[\alpha_n]$ 
is the unique extension of ${K}_0$ of degree $\ell^n$ 
inside $\fell{K}$. Moreover, we can suppose that 
$\alpha_n^\ell=\alpha_{n-1}$, where $\alpha_0:=x$. 

In a first approximation, we consider only those ${K}_i$
which contain $x$, as this is a co-final set of finite
Galois sub-extensions of $\fell{K}\,|\,{K}$. In a second
approximation, we consider the subset of those ${K}_i$
for which $\alpha_{n_i}\in{K}_i$. We claim that this is 
also co-final in the set of all finite Galois sub-extensions 
of $\fell{K}\,|\,{K}$. Indeed, it is sufficient to show
that the restriction of $g$ to ${K}_i[\alpha_{n_i}]$
has order $\ell^{n_i}$. But this follows from the fact
that the restriction of $g$ to both ${K}_i$ and 
${K}_{n_i}:={K}[x,\alpha_{n_i}]$ is $\ell^{n_i}$. 

Now for a co-final set ${K}_i|{K}$ as above, let $v_i$
and $g_i$ with the properties i), ii), above. We claim
that $v_i(x)$ is not divisible by $\ell$ in $v_i{K}$.
Indeed, otherwise let $u\in{K}$ such that $\ell\!\cdot\! 
vu=vx$. Then $y=x/u^\ell$ is a $v_i$-unit. Thus in the 
inertia field ${K}^{T_i}$ we have: $y$ is an $\ell^{\rm th}$ 
power, and hence $x$ is an $\ell^{\rm th}$ power in
${K}^{T_i}$. Therefore, $\alpha_1$ lies in ${K}^{T_i}$. 
Hence $g_i$ acts trivially on $\alpha_1$. This contradicts 
the fact that $g$ acts non-trivially on $\alpha_1\,$!
Therefore, we finally have: $x$ is not a $v_i$-unit, and 
$\fell{K}={K}^{g_i}[(\alpha_n)_n]$.

Now consider and ultrafilter ${\cal D}$ on the (index set 
$I=\{i\}$ of the) family ${K}_i$ such that for every $n$, 
${\cal D}$ contains all the ${K}_i$ with $n_i\geq n$. 
By the ``general non-standard non-sense'' we have the 
following facts:

\ssn

{a)} \ $^*{K}:={\fell{K}}^I/\,{\cal D}$ carries the valuation 
$^*v=\prod v_i\,/\,{\cal D}$ such that the residual field
$^*{K}\,^*v$ has characteristic $\neq\ell$.

\ssn

{b)} \ ${\eu L}:=\prod{K}^{g_i}/\,{\cal D}$ has a 
canonical embedding into $^*{K}$. Furthermore,
one has $\fell{{\eu L}}={\eu L}[(\alpha_n)_n]$, and
$^*v\,x$ is not divisible by $\ell$ in $^*v\,{\eu L}$.

\ssn

{c)} \ Since each $v_i$ is pro-$\ell$ Henselian on 
${K}^{g_i}$, it follows that the restriction ${\eu v}$ 
of $^*v$ to ${\eu L}$ is pro-$\ell$ Henselian on ${\eu L}$.  

\ssn

From this we deduce: ${\eu L}\cap\fell{K}={K}_0$ is the 
fixed field of $g$. Therefore, if $v$ is the restriction of 
${\eu v}$ to $\fell{K}\subset\fell{{\eu L}}$, then $v$ is 
pro-$\ell$ Henselian on ${K}_0$. And since $v x$ is not 
divisible by $\ell$ in $v {K}_0$, it follows that 
$T_v=\Gell{{K}_0}$. Thus $g$ lies in $\inrt{\fell{K}}\,$.
 
\ssn

To (2): \ Let $g\neq1$ be an inertia element at some $v$. 
Generally, in order to show that $g$ is an accumulation 
point of divisorial inertia elements, w.l.o.g.\ we can 
suppose that $g$ is not an $\ell^{\rm th}$-power in $T_v$. 

\ssn

As in the proof of (1), let ${K}_0$ be the fixed field 
of $g$ in $\fell{K}$, and $(\alpha_n)_n$ the compatible 
system of roots $\alpha_n^{\ell^n}=x$ of $x$ such that
${K}_n:={K}_0[\alpha_n]$ is the unique extension of
degree $\ell^n$ of ${K}_0$. By replacing $x$ by its inverse 
(if necessary), we can suppose that $vx>0$. Finally consider
a co-final set of finite Galois sub-extensions ${K}_i|{K}$
of $\Gell{K}\,|\,{K}$ as in the proof of assertion~(1) above,
i.e., such that: $x\in {K}_i$, and the restriction of $g$
to ${K}_i$ has order $\ell^{n_i}$, and $\alpha_{n_i}\in{K}_i$,
where $\alpha_{n_i}$ are as above.

\ssn
{\bf Hypothesis/Remark.} \ 
For every $i$, there exist projective models $X_i\to{k}$ 
of ${K}^g_i|{k}$ such that the centre ${\eu x}_i$ of $v$ 
is a regular point. 

\ssn
 
We remark that by Abhyankar's desingularization Theorem, 
\nmnm{Abhyan\-kar}~[A], this hypothesis is satisfied if 
$\td({K}|{k})\leq2$. Further, if $v$ is a defectless
valuation, the above hypothesis is satisfied by the main 
result of \nmnm{Knaf--Kuhlmann}, [K--K].

\msn

For every $i$ consider a projective model $X_i\to{k}$ of 
${K}^g_i|{k}$, such that the center ${\eu x}_i$ of $v$ on 
$X_i$ is a regular point. Equivalently, ${\cal O}_v$
dominates the regular local ring ${\cal O}_i:={\cal O}_{X_i,{\eu x}_i}$
of ${\eu x}_i$. In particular, if ${\cal M}_i={\cal O}_i\cap
{\cal M}_v$ is the maximal ideal of ${\cal O}_i$, then $vx>0$
implies $x\in{\cal M}_i$. On the other hand, since ${\cal O}_i$ 
is regular, it is a factorial ring. Let 
$x=\epsilon\pi_1^{\nu_1}\dots\pi_r^{\nu_r}$ 
be its representation as product of powers of non-associate 
prime elements in ${\cal O}_i$. 

\ssn

{\it Claim.\/} \ There is an exponent $\nu_m$ which is not 
divisible by $\ell$. 

\ssn
 
Indeed, by its choice, $x$ is not an $\ell^{\rm th}$ power 
in ${\cal O}_v^g$. Further, ${\cal O}_v^g$ is pro-$\ell$ 
Henselian. As it dominates ${\cal O}_i$, it contains a pro-$\ell$ 
Henselisation ${\cal O}^{\rm h}_i$ of ${\cal O}_i$. Further, 
since $\epsilon$ is unit in ${\cal O}_i$, and the residue 
field $\kappa_i={\cal O}_i/{\cal M}_i={k}$ is algebraically
closed, it follows that $\epsilon$ is an $\ell^{\rm th}$ 
powers in ${\cal O}^{\rm h}_i$. Therefore, if all the exponents 
$\nu_m$ were divisible by $\ell$, then $x$ would be an 
$\ell^{\rm th}$ power in ${\cal O}^{\rm h}_i\subset{\cal O}_v^g$,
contradiction!

\ssn

Now let $\nu_{m_i}$ be an exponent which is not divisible by 
$\ell$. Let $v^g_i$ be the divisorial valuation on ${K}^g_i$ 
defined by the prime element $\pi_{m_i}$, and let $v_i$ be 
some prolongation of $v^g_i$ to $\fell{K}$. Then we have: 

\ssn
 
a) \ The value $v_i(x)=v^g_i(x)=\nu_{m_i}$ is not 
divisible by $\ell$. 

\ssn

b) \ Since $\alpha_i\in{K}_i$ satisfies $\alpha_i^{\ell^{n_i}}=x$, 
we have $\ell^{n_i}\cdot v_i(\alpha_i)=v_i(x)$. Thus $v_i$ is
totally (tamely) ramified in ${K}_i|{K}^g_i$. 

\ssn

In particular, by general decomposition theory, it follows that
the restriction of $T_{v_i}$ to ${K}_i$ contains 
${\rm Gal}\,({K}|{K}^g_i)$. Thus, there exists $g_i\in T_{v_i}$ 
such that $g_i$ and $g$ coincide on ${K}_i$. 

\ssn

This concludes the proof of assertion~(2) of Theorem~2.10. 
\enddemo
\vskip-\baselineskip
\section{First consequences of the local theory}
Let $K|k$ be a function field over with $\td(K|k)>1$
and $k$ an algebraic closure of a finite field. Using the 
Theorem above one can recover from $\Gell K$ certain geometric 
invariants as follows. 

\msn

A) \ {\it Recovering $\td(K|k)$.\/}

\ssn

We first remark that from the generalized Milnor Conjectures
it immediately follows that $\cd(\Gell K)=\td(K|k)$. But this 
fact can be deduced from the local theory as follows:

\proclaim{Fact.} \
{\rm In the above context we have:

\ssn

1) The following conditions are equivalent, and 
they are satisfied if an only if $\td(K|k)=2$:

\ssn

(i) $\Gell{Kv}$ is pro-$\ell$ free for some 
$v\in{\eu D}_{Kv}$.

\ssn

(ii) $\Gell{Kv}$ is pro-$\ell$ free for all 
$v\in{\eu D}_{Kv}$. 

\ssn

2) By induction, we can characterize $d=\td(K|k)$ as 
follows: The following conditions are equivalent, and 
they are satisfied if an only if $\td(K|k)=d$:

\ssn

(j) For some $v\in{\eu D}_K$, the pro-$\ell$
group $\Gell{K v}$ has the property characterizing the 
fact that $\td(Kv|k)=d-1$.

\ssn

(jj) $\Gell{K v}$ has the property characterizing 
the fact that $\td(K v)=d-1$ for all $v\in{\eu D}_v$. 
}
\endproclaim

\demo{Proof} \ It is sufficient to prove 1). First
assume $\td(K|k)=2$. Then $K v$ is a function field 
in one variable over $k$. Thus the result follows by the 
fact that the absolute Galois group of $K v$ is profinite
free, see e.g.\ \nmnm{Harbater}~[Ha], or \nmnm{Pop}~[Po].
Conversely, suppose that $\Gell{K v}$ is pro-$\ell$ 
free for some $v\in{\eu D}_K$. By contradiction,
suppose that $\td(K|k>2$. Since $\td(K v|k)=\td(K|k)-1>1$, 
$\Gell{K v}$ contains divisorial subgroups of the form 
$Z_w\cong\Zell\times\Gell{(K v)w}$, which have $\cd>1$.
Thus $\Gell{K v}$ is not pro-$\ell$ free, contradiction!
\enddemo

\ssn

B) \ {\it Recovering the residually divisorial inertia\/}

\ssn 

In the usual hypothesis from above, let $v$ be a Zariski 
prime divisor, and $\Gell{Kv}$ the residual Galois extension. 
In order to have a name, we will say that the divisorial 
inertia elements $\,\divinrt{\fell{Kv}}$ 
are {\it $v$-residually divisorial inertia\/} elements of
$\Gell{K}$. 

\ssn

By the local theory, if $\td(K|k)>2$, then $\td(K v|k)>1$. 
Thus the Zariski prime divisors of $Kv$ are encoded in 
$\Gell{Kv}$ as shown in Proposition~\chardivsubg\ of the Local 
Theory. In particular, the $v$-residually divisorial inertia 
elements of $\Gell{Kv}$ are encoded in $\Gell{Kv}$ (and 
indirectly, in $\Gell K$ too). 

\ssn

{\it Nevertheless,\/} if $\td(K|k)=2$, a fact which by 
the subsection above is encoded in $\Gell K$, then the
residue Galois group $\Gell{K v}$ is pro-$\ell$ free. 
Thus there is no group theoretic way to recognise the
inertia elements of $\Gell{Kv}$ from $\Gell{Kv}$ itself,
as this group does not depend on $Kv$ at all...

\ssn

The solution to recovering the divisorial inertia elements 
in the residual Galois group $\Gell{Kv}$ in general in 
the following: Let $v_0$ be a Zariski divisor of $Kv$, 
and $T_{v_0}\subset Z_{v_0}$ its inertia, respectively 
decomposition group in $\Gell{Kv}$. Further let 
$w=v_0\circ v$ the composition of $v_0$ with $v$. Thus 
one has $Kw=(Kv)v_0$ and $0\to vK\to wK\to v_0(Kv)\to0$ in
a canonical way. Further, via the canonical exact sequence 
$
1\to T_v\to Z_v\hor{\pi_v}\Gell{Kv}\to1
$ 
we have: $T_w\subset Z_w$ are exactly the pre-images of
$T_{v_0}\subset Z_{v_0}$ via $\pi_v$. 

\ssn

On the other hand, $T_w$ consists of inertia elements of
$\Gell K$, thus by Theorem~\inertiaelem\ it follows that
$T_w$ is contained in the topological closure of 
$\,\divinrt{\fell K}$. And this last set is known, as
all the divisorial subgroups $Z_v$ together with their
inertia groups $T_v$ are encoded in $\Gell K$, as shown
in Proposition~\chardivsubg. Finally, 
in the ``trouble case'' $\td(K|k)=2$, we deduce the 
following: The image of $\,\divinrt{\fell K}\cap Z_v$ in 
$\Gell{Kv}$ via $\pi_v$ consists of exactly all the 
$v$-residually inertia elements of $\Gell{Kv}$. And finally, 
by general decomposition theory, a $v$-residually inertia 
element $g_v$ generates an inertia group of $\Gell{Kv}$ 
if and only if the pro-cyclic closed subgroup $T_{g_v}$
generated by $g_v$ in $\Gell{Kv}$ is maximal among the 
pro-cyclic subgroups of $\Gell{Kv}$. Summarising we the 
following: 
\proclaim{Fact} \
{\rm In the above context, let $g_v\in\Gell{Kv}$ be a 
given element, and $T_{g_v}$ the pro-cyclic subgroup 
generated by $g_v$. Then $T_{g_v}$ is an inertia 
subgroup of $\Gell{Kv}$ if and only if the following 
holds:

\ssn

1) If $\td(K|k)=2$, then $g_v$ lies in the image of 
$\,\divinrt{\fell K}\cap Z_v$ in $\Gell{Kv}$, and $T_{g_v}$ 
is maximal among the pro-cyclic subgroups of $\Gell{Kv}$.

\ssn

2) If $\td(K|k)>2$, then $T_{g_v}$ is the inertia group
of a divisorial subgroup of $\Gell{Kv}$ as resulting 
from Proposition~\chardivsubg\ from the Local Theory for 
$\Gell{Kv}$. 
}
\endproclaim
\label\vresinertia\theoremnum*

C) {\it More about the case $\td(K|k)=2$\/}

\ssn

We first recall some basic facts concerning Galois
theory over curves and their function fields over
algebraically closed fields as follows. Let $\KK|\kpkp$ 
be a function field in one variable over the
algebraically closed field $\kpkp$. Let $X\to\kpkp$ be
the unique complete normal model of $\KK|\kpkp$, thus 
$X\to\kpkp$ is projective and smooth. The closed 
points $a\in X$ are in one-to-one correspondence with 
the divisorial valuations $v_a$ of $\KK|\kpkp$.

\ssn

Let $\ell\neq{\rm char}$ be some prime number, and
$\Gprm\KK={\rm Gal}(\KK^{\ell, ab}|\KK)$ be the Galois 
group of a maximal abelian pro-$\ell$ Galois extension of 
$\KK$. Let $T_a=Z_a$ be the inertia/decomposition groups 
over each $a\in X(\kpkp)$ in $G$. (N.B., these groups
depend only on $a$, and not on the specific Zariski 
prime divisor of $\KK^{\ell,{\rm ab}}|\kpkp$ used to 
define each particular $T_a=Z_a$.)

\ssn

a) \ {\it Canonical inertia generators\/}

\ssn
Let $g\geq0$ be genus of $X$. It is well known (e.g., using 
for instance generalized Jacobians, or the specialization
theorem of Grothendieck) that there exist elements 
$\rho_1,\sigma_1,\dots,\rho_g, \sigma_g$ in $\Gprm\KK$, 
and ``Abelian inertia elements'' $\tau_a\in T_a$ (all 
$a\in X(\kpkp)$) such that one has: 
$$
\Gprm\KK=G=<\rho_1,\sigma_1,\dots,\rho_g,\sigma_g,
    (\tau_a)_a \>\mid\,{\scriptstyle\prod}_a\tau_a=1>.
$$

In particular, if $\tau'_a\in T_a$ are such that 
$\prod_a\tau'_a=1$ inside $G$, then there exists a unique
$\epsilon\in{\lv Z}_\ell$ such that $\tau'_a=\tau^\epsilon_a$ 
for all $a\in X$.

\ssn

We will say that $(\tau^{\rm ab}_a)_{a\in X}$ is a 
{\it canonical system of Abelian inertia generators\/} in 
$\Gprm\KK$. Such a system is unique up to simultaneous 
multiplication by an element of ${\lv Z}_\ell^\times$.

\bsn

b) \ {\it $\piellab(X)$ and the genus\/}

\msn
In the context above, let $T'_\KK$ be the closed 
subgroup of $\Gprm\KK$ generated by the inertia subgroups 
$T_a$ (all $a$). Then we have a canonical exact sequence 
of the form:
$$
1\to T'_\KK\to\Gprm\KK\to\pi'_1(X)\to1\,.
$$

In particular, $\pi'_1(X)$ is the free Abelian pro-$\ell$
of rank $2g$, where $g$ is the genus of $X$. Thus the genus 
of $X$ is encoded in $\Gellab\KK$ endowed via the inertia
groups $(T_a)_{a\in X}$.

\ssn

c) \ {\it Detecting\/} ${\rm Div}^0(X)$

\ssn
For $a\in X$, recall the commutative diagram form section~1, A),~4): 
$$
\def\normalbaselines{%
                    \baselineskip20pt\lineskip3pt\lineskiplimit3pt}
\matrix{
       \; \KK^\times   &\horr{\rm can}
                         &\hhhom{\Gprm\KK}{{\lv Z}_\ell}   &\cr
     \dwn{v_a}         &           &     \dwn{{\rm res}}      &\cr
       \; v_a\KK       &\horr{\jmath^a}
                         &\hhhom{T_a}{{\lv Z}_\ell}          &\cr }
\leqno{\hbox to25pt{}(a)}
$$

For every $a$, let $\gamma_a\in v_a\KK$ be the unique positive
generator, thus $\gamma_a=v_v\pi_a$ for every uniformizing 
parameter $\pi_a$ at $a$. Then $T_a$ has a {\it unique generator\/} 
$\tau^0_a$ such that $\jmath^a\,(\gamma_a)(\tau^0_a)=1\in{\lv Z}_\ell$. 
In particular, if $x\in\KK^\times$ is arbitrary, and 
$v_a x=\gamma_x=m_{a,x}\,\gamma_a$ for some integer $m_{a,x}$, 
then $\jmath_{v_a}(\gamma_x)(\tau'_a)=m_{a,x}$ inside 
${\lv Z}_\ell$. We will call $\tau^0_v$ the canonical inertia 
generator of $T_a$ (all $a$). 

\ssn

It is a well known fact, that ${\eu T}^0=(\tau^0_a)_a$ is also 
a canonical system of Abelian inertia generators in $\Gprm\KK$, 
and ${\rm Div}^0(X)$ is {\it canonically\/} isomorphic to 
$
{\eu D}^0_{{\eu T}^0}=\{\psi\in\hhmell{T'_\KK}
 \mid \psi(\tau^0_a)\in{\lv Z}\hbox{ (all $a$)}, 
\ \psi(\tau^0_a)=0 \hbox{ (almost all $a$)}\,\}.
$ 
\proclaim{Fact} \
{\rm 
Let ${\eu T}=(\tau_a)_a$ be an arbitrary canonical system of 
Abelian inertia generators, thus $\tau^0_a=\tau_a^\epsilon$ 
(all $a$) for some $\ell$-adic unit $\epsilon$. Further denote 
$$
{\eu D}^0_{\eu T}=\{\psi\in\hhmell{T'_\KK}
 \mid \psi(\tau_a)\in{\lv Z}\hbox{ (all $a$)}, 
\ \psi(\tau_a)=0 \hbox{ (almost all $a$)}\,\}.
$$
Then $\epsilon\cdot{\rm Div}^0(X)=\epsilon\cdot
  {\eu D}^0_{{\eu T}^0}={\eu D}^0_{\eu T}$ inside 
$\hhmell{\Tellab\KK}$. 
}
\endproclaim

\demo{Prof} \ Clear from the discussion above. \enddemo

\proclaim{Fact} \ 
{\rm Now let $K|k$ be a function 
field with $\td(K|k)=2$, and $k$ the algebraic closure of
a finite field. Let $v$ be a Zariski prime divisor of $K$,
and $\Gell{Kv}$ the corresponding residual Galois group. 
Then by Fact~\vresinertia\ above, the inertia groups in
$\Gell{Kv}$ are known. Thus applying the above discussion
in the case $\KK=Kv$ and $X=X_v$, we get the following:

\ssn

1) The inertia groups $(T_a)_{a\in X_v}$ and the canonical 
generating systems of inertia ${\eu T}_v=(\tau_a)^\nix_{a\in X_v}$
can be deduced in a group theoretic way from $\Gell K$.

\ssn

2) The exact sequence 
$
1\to\Tellab{Kv}\to\Gellab{Kv}\to\piellab(X_v)\to1\,,
$
can be deduced in a group theoretic way from $\Gell K$.
Thus also the genus $g_v$ of $X_v$, as $\piellab(X_v)$
has $\Zell$-rank equal to $2g_v$. 

\ssn

3) Each canonical generating systems of inertia 
${\eu T}_v=(\tau_a)^\nix_{a\in X_v}$ defines in a 
canonical way the subgroup 
$$
{\eu D}^0_{{\eu T}_v}=\{\psi\in\hhmell{\Tellab{Kv}}
 \mid \psi(\tau_a)\in{\lv Z}\hbox{ (all $a$)}, 
\ \psi(\tau_a)=0 \hbox{ (almost all $a$)}\,\}
$$
which up to multiplication by an $\ell$-adic unit equals
the image of ${\rm Div}^0(X_v)$ in $\hhmell{\Tellab{Kv}}$ 
defined canonically via Kummer Theory.
}
\endproclaim 
\label\dequalstwo\theoremnum*

\ssn

\vfill\eject

D) \ {\it Geometric families of Zariski prime divisors}

\ssn

First we recall the basic definitions and facts.
Let $\KK|\kpkp$ be a function field over an algebraically
closed base field $\kpkp$.  

\proclaim{Fact/Definition} {\rm 
In the above context, for a quasi-projective normal model 
$X\to\kpkp$ of $\KK|\kpkp$, let ${\eu D}_X$ be the Zariski 
prime divisors of $\KK|\kpkp$ it defines. We will say 
that a subset ${\eu D}\subset{\eu D}_{\KK|\kpkp}$ is a 
{\it geometric set of Zariski prime divisors,\/} if it
satisfies the following equivalent conditions:

\ssn

{(i)} For all projective normal models 
$X\to\kpkp$ of $\KK|\kpkp$ one has: ${\eu D}$ and 
${\eu D}_X$ are almost equal.\footnote{$^2$}{{We say 
that two sets are almost equal, if their symmetric 
difference is finite.}}

\ssn

{(ii)} There exists quasi-projective normal 
models $X'\to\kpkp$ of $\KK|\kpkp$ such that 
the equality ${\eu D}={\eu D}_{X'}$ holds. 
}
\endproclaim

Our aim in this subsection is to recall a criterion for 
describing the geometric sets ${\eu D}$ of divisorial 
valuations of $\KK|\kpkp$, and to show that the criterion
is encoded group theoretically in $\Gell K$ in the case
$K|k$ is a function field with $\td(K|k)>1$ and $k$ the 
algebraic closure of a finite field. See [P4], Section~3,
for more details for general function fields.   

\ssn

Recall that a {\it line\/} on a $\kpkp$-variety is by 
definition an integral $\kpkp$-subvariety ${\eu l}\,\subseteq X$, 
which is a curve of geometric genus equal to $0.$ We 
denote by $X^{\rm line}$ the union of all the lines on $X$. 

\ssn

We will say that a variety $X\to\kpkp$ is \vurl\ if the 
set $X^{\rm line}$ is not dense in $X$. In particular, 
a curve $X$ is \vurl\ if and only if its geometric genus  
$g_X$ is positive.

\ssn
 
Further recall that being \vurl\ is a birational 
notion. Thus we will say that a function field $\KK|\kpkp$ 
with $\td(\KK|\kpkp=d>0$ is {\it \vurl,\/} if $\KK|\kpkp$ 
has models $X\to\kpkp$ which are \vurl. 

\ssn

Suppose that $d>1$. We call a Zariski prime divisor of 
$\KK|\kpkp$ {\it \vurl,\/} if $\KK v\,|\,\kpkp$ is 
\vurl. A Zariski prime divisor $v$ of $\KK|\kpkp$ 
is \vurl\ if and only if there exists a normal model 
$X\to\kpkp$ of $\KK|\kpkp$, and a \vurl\ prime Weil 
divisor $X_1$ of $X$ such that $v=v_{X_1}$. 

\smallskip

Now let ${\eu D}$ be a set of Zariski prime divisors of 
$\KK|\kpkp$. For every finite extension $\KK_i|\KK$, let 
${\eu D}_{i}$ be the prolongation of ${\eu D}$ to $\KK_i$. 
Thus in particular, ${\eu D}_{i}$ is a set of Zariski 
prime divisors of $\KK_i|\kpkp$. As shown in loc.cit., 
one has the following:

\proclaim{Proposition} In the context above, set 
$d=\td(\KK|\kpkp)$. A set ${\eu D}$ of Zariski prime 
divisors of $\KK|\kpkp$ is geometric if and only if the 
following conditions are satisfied: 

\ssn

{\rm(i)} There exists a finite $\ell$-elementary extension
$\KK_0|\KK$ of degree $\leq\ell^d$ such that ${\eu D}_{0}$ 
is almost equal to the set of all \vurl\ prime divisors of 
$\KK_0$.

\ssn

{\rm(ii)} If $\KK_2|\KK$ is any $\ell$-elementary extension 
of degree $\leq\ell^d$, and $\KK_1=\KK_2\KK_0$ is the 
compositum, then ${\eu D}_{1}$ is almost equal to the set 
of all \vurl\ prime divisors of $\KK_1$.
\endproclaim
\label\geomsetdiv\theoremnum*

From this Proposition one deduces the following inductive
procedure on $d=\td(\KK|\kpkp)$ for deciding whether a
given set ${\eu D}$ of Zariski prime divisors is geometric,
respectively whether $\KK|\kpkp$ is \vurl. 

\proclaim{Criterion} 
{\rm In the above context one has:

\ssn

1) Case $d=1$: 

\ssn

${\cal P}_{\rm geom}^{(1)}({\eu D})$: \
A set of Zariski prime divisors ${\eu D}$ of $\KK|\kpkp$ 
is geometric if and only if it is almost equal to the set 
of all Zariski prime divisors ${\eu D}_{\KK|\kpkp}$. 

\ssn

${\cal P}_{\rm v.u.}^{(1)}(\KK|\kpkp)$: \
$\KK|\kpkp$ is \vurl\ if and only if the geometric genus 
$g_X$ of the complete normal model $X\to\kpkp$ of $\KK|\kpkp$  
satisfies $g_X>0$.

\msn

2) Case $d>1$: Then by induction on $d$, we already 
have criteria ${\cal P}_{\rm geom}^{(\delta)}$ and 
${\cal P}_{\rm v.u.}^{(\delta)}$ which assure that sets of 
Zariski prime divisors are geometric, respectively that 
function fields are \vurl\ (all $1\leq\delta<d$). 
We now make the induction step $d=\td(\KK|\kpkp)$ as 
follows:

\ssn

${\cal P}_{\rm geom}^{(d)}({\eu D})$: \ The criterion 
given by Proposition~\geomsetdiv\ is satisfied for the
set of Zariski prime divisors ${\eu D}$ of $\KK|\kpkp$. 

\ssn

Here we remark that given a finite extension $\KK_i|\KK$, 
the assertion ``$\,v_i$ {\it is a \vurl\ Zariski prime divisor 
of\/} $\KK_i|\kpkp\,$'' which is essentially used in loc.cit., is 
actually equivalent to ``$\,\KK_i v_i|\kpkp$ {\it satisfies\/} 
$\,{\cal P}_{\rm v.u.}^{(d-1)}(\KK_i v_i|\kpkp)\,$''. 

\ssn

${\cal P}_{\rm v.u.}^{(d)}(\KK|\kpkp)$: \ There exists
a subset ${\eu D}\subset{\eu D}_{\KK|\kpkp}$ such that  
${\cal P}_{\rm geom}^{(d)}({\eu D})$ holds; and further,
${\cal P}_{\rm v.u.}^{(d-1)}(\KK v\,|\kpkp)$ is true for 
almost all $\KK v\,|\kpkp$ ($v\in{\eu D}$).
}
\endproclaim
\label\criterion\theoremnum*

We finally come to case $K|k$ with $\td(K|k)>1$ and $k$
the algebraic closure of a finite field. We show that
Criterion~\criterion\ can be interpreted in $\Gell K$. 

\ssn

Indeed, by Proposition~\chardivsubg, the Zariski prime 
divisors of $\fell  K$ are in bijection with the divisorial 
subgroups of $\Gell K$ via $v\mapsto Z_v$; thus a 
bijection $\fell{\eu D}\mapsto{\cal Z}_{\fell{\eu D}}$ 
from sets of Zariski prime divisors 
$\fell{\eu D}\subset{\eu D}_{\fell K}$ to sets of 
divisorial subgroups ${\cal Z}$ of $\Gell K$. 
Moreover, a set $\fell{\eu D}$ is the prolongation 
to $\fell K$ of a set of Zariski prime divisors 
${\eu D}\subset{\eu D}_K$ if and only if 
${\cal Z}_{\fell{\eu D}}$ is invariant under 
conjugation in $\Gell K$. If this is the case, 
then every conjugacy class represents an element 
of ${\eu D}$.

\ssn

For ${\cal Z}^\ell$ a set of divisorial subgroups of 
$\Gell K$, which is closed under conjugation in $\Gell K$,
let ${\eu D}_{{\cal Z}^\ell}$ be the corresponding set
of Zariski prime divisors of $\fell K$. For every 
sub-extension $K_i|K$ of $\fell K|K$, let ${\eu D}^i_{\cal Z}$
be the restriction of $\fell{\eu D}_{\cal Z}$ to $K_i$.
Thus ${\eu D}^i_{\cal Z}$ is exactly the prolongation of 
${\eu D}_{\cal Z}$ to $K_i$. 

\ssn

Now suppose that $d=2$. Let $v$ be some Zariski prime 
divisor of $K|k$. Let further $X_v\to k$ be the complete 
smooth model of $Kv|k$. Then by Fact~\dequalstwo, the 
inertia groups $T_a\subset\Gell{Kv}$ (all $a\in X_v$), 
thus $\pi_1(X_v)$ and the genus $g_{X_v}$ of $X_v$ can be 
recovered from $\Gell K$. The same is true if we replace 
$K|k$ by some finite sub-extension $K_i|K$ of $\fell K|K$. 

\ssn

Therefore, by induction on $d=\td(K|k)>0$, we have the 
following Galois translation of the Criterion above.

\proclaim{{\rm Gal}-Criterion} 
{\rm In the above context one has:

\ssn

1) Case $d=2$: 

\ssn

Gal$\,{\cal P}_{\rm geom}^{(2)}({\eu D}_{\cal Z})$ is said 
to be satisfied, if the Galois theoretic translation of 
$\,{\cal P}_{\rm geom}^{(2)}({\eu D}_{\cal Z})$ given above
holds.

\ssn

Gal$\,{\cal P}_{\rm v.u.}^{(2)}(K|k)\,$ is said to be 
satisfied, if for some ${\eu D}_{\cal Z}$ as above one has:
Gal$\,{\cal P}_{\rm geom}^{(2)}({\eu D}_{\cal Z})$ holds, 
and $g_{X_v}>0$ for almost all $v\in{\eu D}_{\cal Z}$. 

\ssn

2) Case $d>1$: 

\ssn

By induction on $d$, we already have the Galois translation 
Gal$\,{\cal P}_{\rm geom}^{(\delta)}$ and 
Gal$\,{\cal P}_{\rm v.u.}^{(\delta)}$ of the criteria 
${\cal P}_{\rm geom}^{(\delta)}$ and 
${\cal P}_{\rm v.u.}^{(\delta)}$ (all $1\leq\delta<d$). 
We now make the induction step $d=\td( K| k)$ as follows:

\ssn

Gal$\,{\cal P}_{\rm geom}^{(d)}({\eu D}_{\cal Z})$ is said 
to be satisfied, if (the Galois translation of) the criterion 
given by Proposition~\geomsetdiv\ is satisfied. 

\ssn

Here we remark that given a finite extension $K_i|K$, the 
assertion ``$\,v_i$ {\it is a \vurl\ Zariski prime divisor of\/} 
$K_i|k,$'' which is essentially used in loc.cit., is actually 
equivalent to Gal$\,{\cal P}_{\rm v.u.}^{(\delta)}(K_i v_i\,|\,k)$

\ssn

Gal$\,{\cal P}_{\rm v.u.}^{(d)}(K|k)$ is said to be satisfied 
if and only if for some set ${\eu D}_{\cal Z}$ as above 
satisfying Gal$\,{\cal P}_{\rm geom}^{(d)}({\eu D}_{\cal Z})$ 
one has: For almost all conjugacy classes of decomposition 
groups $Z_v$, the condition Gal$\,{\cal P}_{\rm v.u.}^{(d-1)}(Kv|k)$ 
is satisfied.
}
\endproclaim
\label\galcrit\theoremnum*
%
%
%
%
%
\section{Abstract pro-$\ell$ Galois theory}
In this section we develop a ``pro-$\ell$ Galois theory'', 
which in some sense has a flavor similar to one of the 
abstract class field theory. The final aim of this theory 
is to provide a machinery for recovering function fields 
(over algebraically closed base fields) from their pro-$\ell$ 
Galois theory in an axiomatic way. The material here is a 
simplified version of the corresponding part from 
\nmnm{Pop}~[P4].

\msn

A) \ {\it Axioms and definitions\/}

\ssn

The context is the following: Let $\ell$ be a fixed 
prime number, and $\calz$ a quotient of ${\lv Z}_\ell$. 

\nonumproclaim{Definition} {\rm
A {\it \abstr\ level $\delta\geq0$ \/} (pro-$\ell$) 
{\it \cycellprefrm\/} ${\cal G}=\big(\ggg, {\eu Z})$ is 
defined by induction on $\delta$ to be an Abelian pro-$\ell$ 
group $\ggg$ endowed with extra structure as follows: 

\msn
Axiom I) \ The level $\delta=0$: 

\ssn

{\it A level $\delta=0$ \abstr\ \cycellprefrm\ is simply\/} 
${\cal G}=(\ggg,\vid)$.

\ssn

We will say that ${\cal G}$ {\it has level $\delta>0$,\/} if 
the following  Axioms II, III, are inductively satisfied. 

\msn
Axiom II) \ Decomposition structure:

\ssn

${\eu Z}=(\zzv)_v$ {\it is a family of closed subgroups of $\ggg$. 
We call $\zzv$ the decomposition group at $v$, and suppose that 
each $\zzv$ is endowed with a subgroup $\ttv$ such that:\/} 

\ssn

i) {\it $\ttv\cong\calz$ as an abstract pro-$\ell$ groups, 
and $\ttv\cap\ttt_w=\{1\}$ for all $v\neq w$.\/}  

\ssn

For every co-finite subset ${\eu U}_i$ of valuations $v$, let 
$\ttt_{{\eu U}_i}$ be the closed subgroup of $\ggg$ generated 
by the $\ttv$ ($v\in{\eu U}_i$). A system $({\eu U}_i)_i$ of 
such ${\eu U}_i$ is called co-final, if every finite set of 
valuations is contained in the complement of some ${\eu U}_i$. 

\ssn

ii) {\it There exist co-final systems $({\eu U}_i)_i$ 
with $\ttv\cap T_{{\eu U}_i}=\{1\}$ (all $i$, and 
$v\not\in{\eu U}_i$).\/}

\msn
Axiom III) \ Induction:

\ssn

{\it Every $\ggv:=\zzv/\ttv$ caries itself the structure of a 
\abstr\ \cycellprefrm\ of level $(\delta-1)$.\/}
}
\endproclaim

{\bf Convention.} 
Let ${\cal G}=(\ggg,{\eu Z})$ be a level 
$\delta\geq0$ \abstr\ \cycellprefrm. In order to have a
uniform notation, we enlarge the index set $v$ (which 
in the case $\delta=0$ is empty) by a new symbol, which 
we denote $v_*$, by setting ${\zzv}_{_*}=\ggg$ and 
${\ttv}_{_*}=\{1\}$. We will say that $v_*$ is the 
``trivial valuation''. Thus $\ggg$ is the decomposition 
group of the trivial valuation, and its ``inertia group'' 
is the trivial group. In particular, the residue Galois 
group at $v_*$ is ${\ggv}_{_*}=\ggg$. 

\proclaim{Definition/Remark}
{\rm 
Let ${\cal G}=(\ggg,{\eu Z})$ be a \abstr\ 
\cycellprefrm\ of some level $\delta\geq0$. Consider 
any $\alpha$ such that $0\leq\alpha\leq\delta$.  

\ssn

1) By induction on $\delta$ it is easy to see that one can 
view ${\cal G}=(\ggg,{\eu Z})$ canonically as a \abstr\
\cycellprefrm\ of level $\alpha$.

\ssn

2) We define inductively the system 
${\cal G}^{(\alpha)}=({\cal G}_{\alpha,i})_i$  of the 
{\it $\alpha$-residual \abstr\ \cycellprefrm s\/} of 
${\cal G}$ as follows: First, if $\alpha=0$, then this 
system consists of ${\cal G}$ only. In general, if 
$0<\alpha\leq\delta$, we recall that every residue 
group $\ggv$ carries the structure of a \abstr\ 
\cycellprefrm\ of level $(\delta-1)$; thus of level 
$(\alpha-1)$ by remark 1) above. Let ${\cal G}_v$ 
be this \abstr\ \cycellprefrm. Then by induction, 
the system of the $(\alpha-1)$-residual \abstr\ 
\cycellprefrm s ${\cal G}^{(\alpha-1)}_v$ of each 
${\cal G}_v$ are defined. We then set 
${\cal G}^{(\alpha)}=\big({\cal G}_v^{(\alpha-1)}\big)_v$. 

\ssn

We remark that the ``correct'' notation for the system 
of the $\alpha$-residual \abstr\ \cycellprefrm s 
$({\cal G}_{\alpha,i})_i$ is to index it by multi-indices 
of length $\alpha$ of the form 
$\bfv=(v_{i_1}\dots v_{i_\alpha})$, where $v_{i_1}$ 
is a valuation of $\ggg$, $v_{i_2}$ is a valuation of 
${\ggv}_{_1}$, etc.. 

\ssn

3) For every multi-index $\bfv$ as above, let 
$\ggg_{\bfv}$ be the profinite group on which the 
$\bfv$-residual \abstr\ \cycellprefrm\ ${\cal G}_{\bfv}$ 
is based. Then we will call $\ggg_{\bfv}$ a 
{\it $\bfv$-residual group\/} of ${\cal G}$. One 
can further elaborate here as follows: Given a 
multi-index $\bfv=(v_{i_1}\dots v_{i_\alpha})$, 
we can define inductively the following: 

\ssn

a) The $\bfv$-decomposition group $Z_{\bfv}$ of $\ggg$, as 
being ---inductively on $\alpha$--- the pre-image of 
$Z_{(v_{i_2}\dots v_{i_\alpha})}\subseteq\ggg_{v_{i_1}}$ 
in ${\zzv}_{_{i_1}}$ via $Z_{v_{i_1}}\to\ggg_{v_{i_1}}$. 

\ssn

b) The inertia group $T_{\bfv}$ of $Z_{\bfv}$, as being 
the kernel of $Z_{\bfv}\to\ggg_{\bfv}$. From the definition 
it follows that $T_{\bfv}\cong\calz^\alpha$.
} 
\endproclaim
\label\deltaresidual\theoremnum*

\proclaim{Definition/Remark} {\rm 
Let ${\cal G}=\big(\ggg,(\zzv)_v\big)$ be a \abstr\
\cycellprefrm\ of some level $\delta'\geq0$.  

\ssn

1) We denote $\whlx_{\cal G}=\hhhom{\ggg}{\calz}$ and 
call it the ($\ell$-adic completion of the) {\it \abstr\ 
pre-field formation defining\/} ${\cal G}$. 

\msn

From now on suppose that $\delta'>0$. 

\msn

2) Let $\ttt\subseteq\ggg$ be the closed subgroup generated 
by all the inertia groups $\ttv$ (all $v$). We set
$\pi_{1,{\cal G}}:=\ggg/\ttt$ and call it the abstract 
fundamental group of ${\cal G}$. One has a canonical exact
sequence
$$
1\to\ttt\to\ggg\to\pi_{1,{\cal G}}:=\ggg/\ttt\to1
$$
Taking $\calz$-Homs, we get an exact sequence of the form
$$
0\to\widehat U_{\cal G}:=\hhhom{\pi_{1,{\cal G}}}{\calz}
    \hor{\rm can}\whlx_{\cal G}:=\hhhom{\ggg}{\calz}
         \hor{\jmath^{\cal G}}\hhhom{\ttt}{\calz}.
$$
We will call $\widehat U_{\cal G}=
        \hhhom{\pi_{1,{\cal G}}}{\calz}$
the {\it unramified part\/} of $\whlx_{\cal G}$. And if no 
confusion is possible, we will identify $\widehat U_{\cal G}$ 
with its image in $\whlx_{\cal G}$. 

\ssn

3) We now have a closer look at the structure of 
$\whlx_{\cal G}$. For an arbitrary $v$ we have 
inclusions $\ttv\hookrightarrow\zzv\hookrightarrow\ggg$. 
Thus we can/will consider/denote restriction maps as 
follows: 
$$
\jmath^v:\whlx_{\cal G} =\hhhom{\ggg}{\calz}
      \horr{\rszv}\hhhom{\zzv}{\calz}
            \horr{{\rm res}_v}\hhhom{\ttv}{\calz}.
$$ 
We set $\widehat U_v=\ker(\jmath^v)$ and call it the 
{\it $v$-units\/} in $\whlx_{\cal G}$. Thus the 
unramified part of $\whlx_{\cal G}$ is exactly 
$\widehat U_{\cal G}=\cap_v\ker(\jmath^v)$. 

\ssn

We further denote 
$\LX_{{\cal G},{\rm fin}}=\{\,x\in\whlx_{\cal G}
         \mid\jmath^v(x)=0\hbox{ for almost all } v\,\}$.
We remark that by Axiom II,~ii), $\LX_{{\cal G},{\rm fin}}$
is dense in $\whlx_{\cal G}$. Indeed, for a co-final 
system $({\eu U}_i)_i$ as at loc.cit., denote 
$\ggg_i=\ggg/\ttt_{{\eu U}_i}$ and $\ttt_i=\ttt/T_{{\eu U}_i}$. 
We have a canonical exact sequence
$$
1\to\ttt_i\to\ggg_i\to\pi_{1,{\cal G}}\to1\,,
$$
and $\ttt_i$ is generated by the images $T_{v,i}$ of 
$\ttv$ in $\ggg_i$ ($v\not\in{\eu U}_i$). Clearly, the 
image of the inflation map
$$
{\rm inf}_i:\hhhom{\ggg_i}{\calz}\to\hhhom{\ggg}{\calz}
$$
is exactly $\Delta_i:=\{x\in\whlx_{\cal G}\mid\jmath^v(x)=0
\hbox{ for all } v\in{\eu U}_i\,\}$. Finally, taking
inductive limits over the co-final system $({\eu U}_i)_i$,
the density assertion follows.

\ssn

A closed submodule $\Delta\subset\whlx_{\cal G}$ 
is said to have {\it finite co-rank,\/} if 
$\Delta\subset\LX_{{\cal G},{\rm fin}}$, and 
$\Delta\hhb1/\,\widehat U_{\cal G}$ is a finite 
$\calz$-module (or equivalently, $\Delta$ is contained 
in $\ker(\jmath^v)$ for almost all $v$). Clearly, the sum 
of two finite co-rank submodules of $\whlx_{\cal G}$ is 
again of finite co-rank. Thus the set of such submodules 
is inductive. And one has:
$$
\LX_{{\cal G},{\rm fin}}=\cup_\Delta\,
              \hbox{(all finite co-rank $\Delta$)}
$$

4) By the discussion above, the family $(\jmath^v)_v$ 
gives rise canonically to a continuous homomorphism 
$\widehat{\oplus_v}\,\jmath^v$ of $\ell$-adically 
complete $\Zell\,$-modules 
$$
\widehat{\oplus_v}\,\jmath^v:
   \whlx_{\cal G}=\hhhom{\ggg}{\calz}\to
     \hhhom{\ttt}{\calz}\hookrightarrow
      \widehat{\oplus_v}\,\hhhom{\ttv}{\calz}\,.
$$ 
We will identify $\hhhom{\ttt}{\calz}$ with its 
image inside $\widehat{\oplus_v}\,\hhhom{\ttv}{\calz}$.
Therefore, $\jmath^{\cal G}=\widehat{\oplus_v}\,\jmath^v$
on $\whlx_{\cal G}$. 

\ssn

We will denote $\widehat{\rm Div}_{\cal G}:= 
   \widehat{\oplus_v}\,\hhhom{\ttv}{\calz}$ 
and call it the ($\ell$-adic completion of the) 
{\it abstract divisor group\/} of ${\cal G}$. We 
will say that the image of $\whlx_{\cal G}$ in 
$\widehat{\rm Div}_{\cal G}$ is the {\it divisorial 
quotient\/} (or the divisorial part) of $\whlx_{\cal G}$.

\ssn

We further set 
$\widehat{\eu Cl}_{\,\cal G}={\rm coker}(\jmath^{\cal G})$, 
and call it the ($\ell$-adic completion of the) {\it abstract 
divisor class group\/} of ${\cal G}$. Therefore, we finally 
have a canonical exact sequence
$$
0\to\widehat U_{\cal G}\>\hookrightarrow\>
    \whlx_{\cal G}\>\hor{\jmath^{\cal G}}\>
        \widehat{\rm Div}_{\cal G}
          \>\hor{\rm can}\>\widehat{\eu Cl}_{\,\cal G}\to0\,.
$$  

5) We say that ${\cal G}$ is {\it complete curve like\/} 
if the following holds: There exist generators $\tau_v$ 
of $\ttv$ such that $\prod_v\tau_v=1$, and this is the 
only pro-relation satisfied by the system of elements 
${\eu T}=(\tau_v)_v$. 

\ssn

Further consider $0\leq\delta<\delta'$. 
We say that ${\cal G}$ is {\it $\delta$-residually 
complete curve like\/} if all the $\delta$-residual 
\abstr\ \cycellprefrm s ${\cal G}_{\bfv}$ are residually 
complete curve like. In particular, ``$0$-residually 
complete curve like'' is the same as ``residually 
complete curve like''.

\ssn

6) Consider the exact sequence
$1\to\ttv\to\zzv\to\ggv:=\zzv/\ttv\to1$ given by Axiom 
II,~i), (all $v$). 
Let ${\rm inf}_v:\hhhom{\ggv}{\calz}\to\hhhom{\zzv}{\calz}$
be the resulting inflation homomorphism. 

\ssn

Since $\ttv=\ker(\zzv\to\ggv)$, it follows that 
$\rszv(\widehat U_v)$ is the image of the inflation 
map ${\rm infl}_v$. Therefore there exists a canonical 
continuous homomorphism, which we call the 
{\it $v$-reduction homomorphism:\/} 
$$
\jmath_v:\widehat U_v\to\hhhom{\ggv}{\calz}=
     \whlx_{{\cal G}_v}
$$
where $\whlx_{{\cal G}_v}$ is the $\ell$-adic 
completion of the abstract $v$-residual field, i.e., 
the one attached to the $v$-residual \cycellprefrm\  
${\cal G}_v$.
 
\ssn

7) Next let $v$ be arbitrary, and ${\cal G}_v$ be the 
corresponding residual \cycellprefrm. To ${\cal G}_v$ 
we have the corresponding exact sequence as defined 
for ${\cal G}$ at point 3) above: 
$$
0\to\widehat U_{{\cal G}_v}
 \hookrightarrow\whlx_{{\cal G}_v}
  \hor{\jmath^{{\cal G}_v}}\widehat{\rm Div}_{{\cal G}_v}
$$   

We will say that ${\cal G}$ is {\it ample,\/} if in the 
notations from point~3) above the following hold:

\ssn

i) The canonical projection
$\whlx_{\cal G}\horr{(\jmath^v)_v}\prod_{v\not\in{\eu U}_i}
                                \hhhom{\ttv}{\calz}$
is surjective. 

\ssn

ii) There exists $v$ such that:

\ssn

\itemitem{a)} $\Delta_i\subseteq\widehat U_v$ and 
      $\jmath_v(\Delta_i)\subseteq\LX_{{\cal G}_v,{\rm fin}}$.

\ssn

\itemitem{b)} $\ker\big(\>\Delta_i\hor{\jmath_v}\whlx_{{\cal G}_v}
  \hor{\jmath^{{\cal G}_v}}\widehat{\rm Div}_{{\cal G}_v}\>\big)
    \subseteq\widehat U_{\cal G}$.

\msn

Next consider $0\leq\delta<\delta'$. We say that 
${\cal G}$ is {\it ample up to level $\delta$,\/} 
if either $\delta=0$, or by induction on $\delta$ the 
following hold: First, if $\delta=1$, then ${\cal G}$ 
is ample in the sense defined above. Second, if $\delta>1$, 
then ${\cal G}_v$ is ample up to level $(\delta-1)$ for 
all $v$.

\ssn

Before going to the next Subsection, let us remark that 
Condition a) and b) for a co-final system $({\eu U}_i)_i$ 
implies the corresponding assertions for every co-finite
subset of valuations $U$, respectively for every finite
co-rank module $\Delta$ (and this fact is obvious).
}
\endproclaim
\label\abstractprefields\theoremnum*

\ssn

B) \ {\it Abstract divisor groups\/}
 
\ssn

{\bf Convention.} In order to avoid too technical 
formulations, we will suppose from now ---if not explicitly 
otherwise stated--- that $\calz=\Zell$. In particular, 
${\lv Z}\subset\zetell$ are subgroups/subrings of 
$\calz$.

\proclaim{Definition} {\rm

\ssn

1) Let $M$ be an arbitrary $\calz$-module. We say that
subsets $M_1,M_2$ of $M$ are {\it$\ell$-adically 
equivalent,\/} if there exists an $\ell$-adic unit 
$\epsilon\in\calz$ such that $M_2=\epsilon\cdot M_1$ 
inside $M$. Correspondingly, given systems $S_1=(x_i)_i$ 
and $S_2=(y_i)_i$ of elements of $M$, we will say 
that $S_1$ and $S_2$ are $\ell$-adically equivalent, 
if there exists an $\ell$-adic unit $\epsilon\in\calz$ 
such that $x_i=\epsilon\, y_i$ (all $i$). 

\ssn

2) Let $M$ be an arbitrary $\ell$-adically complete 
module. We will say that a $\zetell$-submodule 
${\cal M}_\pel\subseteq M$ of $M$ is a 
$\zetell$-{\it lattice\/} in $M$, (for short,
a lattice) if ${\cal M}_\pel$ is a free $\zetell$-module, 
and it is $\ell$-adically dense in $M$, and it satisfies 
the following equivalent conditions:

\ssn

a) $M/\ell={\cal M}_\pel/\ell$

\ssn

b) ${\cal M}_\pel$ has a $\zetell$-basis ${\eu B}$ 
which is $\ell$-adically independent in $M$.

\ssn

c) The condition b) above is satisfied for every 
$\zetell$-basis of ${\cal M}_\pel$.   

\ssn

More general, let $N\subseteq{\cal M}_\pel\subseteq M$ 
be $\zetell$-submodules of $M$ such that $M/N$ is again
$\ell$-adically complete. We will say that ${\cal M}_\pel$ is 
an {\it $N$-lattice\/} in $M$, if $N\subseteq{\cal M}_\pel$,
and ${\cal M}_\pel/N$ is a lattice in $M/N$. 

\ssn

4) Finally, in the context from 3) above, a {\it true lattice\/} 
in $M$ is a free Abelian subgroup ${\cal M}$ of $M$ such that 
${\cal M}_\pel:={\cal M}\otimes{\lv Z}_{(\ell)}$ is 
a lattice in $M$ in the sense of 3) above. And we will say 
that a ${\lv Z}$-submodule ${\cal M}\subseteq M$ containing 
$N$ is a {\it true $N$-lattice\/} in $M$, if ${\cal M}/N$ 
is a true lattice in $M/N$. 
}
\endproclaim

\proclaim{Construction} {\rm
Let ${\cal G}=\big(\ggg,\(\zzv)_v\big)$ be a \abstr\
\cycellprefrm\ which is both {\it ample up to level\/} 
$\delta$ and {\it $\delta$-residually complete curve 
like\/} ($\delta\geq0$ given). Recall the last exact sequence 
from point 3) from Definition/Remark~\abstractprefields:
$$
0\to\widehat U_{\cal G}\>\hookrightarrow\>
    \whlx_{\cal G}\>\hor{\jmath^{\cal G}}\>
        \widehat{\rm Div}_{\cal G}
           \>\hor{\rm can}\>\widehat{\eu Cl}_{\,\cal G}\to0\,.
$$  

The aim of this subsection is to describe the $\ell$-adic 
equivalence class of a lattice $\mardv_\clg$ in 
$\widehat\mardv_\clg$, which ---in the case it exists--- 
will be called the {\it abstract divisor group\/} of 
${\cal G}$. By construction, this will be equivalent to 
giving the equivalence class of a $\widehat U_\clg$-lattice 
in $\whlx_{\cal G}$, which will then be exactly the 
pre-image of $\mardv_\clg$ in $\whlx_{\cal G}$.
}
\endproclaim
\label\construction\theoremnum*

\ssn

{\bf The case $\delta=0$}, i.e., 
          ${\cal G}$ complete curve like.

\ssn

In the notations from Definition/Remark 
\abstractprefields,~5) above, let ${\eu T}=(\tau_v)_v$ 
be the system of generators of the groups $\ttv$ as 
there. Let as call such a system a {\it distinguished 
system of inertia generators.\/} We remark that any 
two distinguished system of inertia generators are 
{\it strictly $\ell$-adically equivalent.\/} Indeed, 
if $\tau'_v\in\ttv$ is another generator of $\ttv$, then 
$\tau'_v=\tau^{-\epsilon_v}_v$ for some $\ell$-adic units 
$\epsilon_v\in\calz$. Therefore, if 
${\eu T}'=(\hhb1\tau'_v\,)^\nix_v$ does 
also satisfy condition ii) from Definition/Remark 
\abstractprefields,~5), then we have also
$\prod_v\,\tau^{-\epsilon_v}_v=1$. By the 
uniqueness of the relation $\prod_v\tau_v=1$,
it follows that $\epsilon_v=\epsilon$ for some 
fixed $\ell$-adic unit $\epsilon\in\calz$, as claimed.

\ssn

Now let ${\eu T}=(\tau_v)_v$ be a distinguished 
system generators of $\ttt$. Further let 
${\cal F}_{\eu T}$ be the Abelian pro-$\ell$ free group 
on the system ${\eu T}$ (written multiplicatively). 
Then one has a canonical exact sequence of pro-$\ell$ 
groups
$$
1\to\,\tau^{\calz}\,\to\,{\cal F}_{\eu T}\,
                                  \to\,\ttt\to1\,,
$$
where $\tau=\prod_v\,\tau_v$ in ${\cal F}_{\eu T}$ is 
the pro-$\ell$ sum of the generators $\tau_v$ (all 
$v$). Taking $\ell$-adically continuous Homs we get an 
exact sequence
$$
0\to\hhhom{\ttt}{\calz}\to
       \hhhom{{\cal F}_{\eu T}}{\calz}
              \to\hhhom{\tau^{\calz}}{\calz}\to 0.
$$
Remark that $\hhhom{{\cal F}_{\eu T}}{\calz}\,\cong\,
     \widehat{\rm Div}_{\cal G}$
in a canonical way, and the last homomorphism is simply 
the summation: $\varphi\mapsto\big(\tau\mapsto
    {\textstyle\sum}_v\,\varphi(\tau_v)\big)$.
Thus $\hhhom{\ttab}{\calz}$ consists of all the 
homomorphisms $\varphi\in\hhhom{{\cal F}}{\calz}$ with 
trivial ``trace''. 

\ssn

Consider the system 
${\eu B}=(\varphi_v)_v$ of all the functionals 
$\varphi_v\in\hhhom{{\cal F}_{\eu T}}{\calz}$ 
defined by $\varphi_v(\tau_w)=1$ if $v=w$, and $0$, 
otherwise (all $v,w$). We denote by 
$$
{\rm Div}_{{\eu T},\pel}=\;<{\eu B}>_{(\ell)}\;
       \subset\hhhom{{\cal F}}{\calz}
$$ 
the $\zetell$-submodule of $\hhhom{{\cal F}}{\calz}$ 
generated by ${\eu B}$. Then ${\eu B}$ 
{\it is an $\ell$-adic basis\/} of 
$\hhhom{{\cal F}}{\calz}$, i.e., ${\rm Div}_{{\eu T},\pel}$ 
is $\ell$-adically dense in $\hhhom{{\cal F}}{\calz}$, 
and there are no non-trivial $\ell$-adic relations 
between the elements of ${\eu B}$. We will say that 
${\eu B}=(\varphi_v)_v$ is the ``dual basis'' to ${\eu T}$. 
We next set
$$
{\rm Div}^0_{{\eu T},\ell}:=\{\,{\textstyle\sum}_v\, 
  a_v\,\varphi_v\in{\rm Div}_{{\eu T},\pel}\mid
    {\textstyle\sum}_v\,a_v=0\,\} =
      {\rm Div}_{{\eu T},\pel}\,\cap\,\hhhom{\ttt}{\calz}.
$$

Clearly, if $\varphi_{v_0}$ is fixed, then the 
system $e_v=\varphi_v-\varphi_{v_0}$ (all $v\neq v_0$) 
is an $\ell$-adically independent $\zetell$-basis of 
${\rm Div}^0_{{\eu T},\pel}$.

\ssn

Thus finally, ${\rm Div}_{{\eu T},\pel}$ is a lattice 
in $\hhhom{\ttt}{\calz}$, and ${\rm Div}^0_{{\eu T},\pel}$ 
is a lattice in $\hhhom{\ttt}{\calz}$.  

\ssn

The dependence of ${\rm Div}_{{\eu T},\pel}$ 
on ${\eu T}=(\tau_v)_v$ is as follows. 
Let ${\eu T}'=(\tau'_v)_v={\eu T}^\epsilon$ be 
another distinguished system of inertia generators. 
If ${\eu B}'=(\varphi'_v)_v$ is 
the dual basis to ${\eu T}'$, then 
$\epsilon\cdot{\eu B}'={\eu B}$. Thus ${\eu B}$ and 
${\eu B}'$ are $\ell$-adically equivalent, and we have: 
${\rm Div}_{{\eu T},\pel}
       =\epsilon\cdot{\rm Div}_{{\eu T}',\pel}$ 
and 
${\rm Div}^0_{{\eu T},\pel}=
     \epsilon\cdot{\rm Div}^0_{{\eu T}',\pel}$. 

\ssn

Therefore, all the subgroups of 
$\hhhom{\ttab}{\calz}$ the form ${\rm Div}_{{\eu T},\pel}$, 
respectively of the form ${\rm Div}^0_{{\eu T},\pel}$, 
are $\ell$-adically equivalent (all distinguished 
${\eu T}$). Hence the $\ell$-adic equivalence classes of 
${\rm Div}_{{\eu T},\pel}$ and ${\rm Div}^0_{{\eu T},\pel}$ 
do not depend on ${\eu T}$, but only on ${\cal G}$. 

\proclaim{Fact} {\rm  
In the above context, denote by $\lxlxzl{\eu T}$ 
the pre-image of $\mardv^0_{{\eu T},\pel}$ in 
$\whlx_{\cal G}$. Further consider all the finite 
co-rank submodules $\Delta$ of $\whlx_{\cal G}$ 
containing $\widehat U_{\cal G}$. Then the following 
assertions are equivalent:

\ssn

{\rm (i)} $\lxlxzl{\eu T}$ is a 
$\widehat U_{\cal G}$-lattice in $\whlx_{\cal G}$.

\ssn

{\rm (ii)} $\Delta\cap\lxlxzl{\eu T}$ is a 
$\widehat U_{\cal G}$-lattice in $\Delta$ (all 
$\Delta$ as above).

\ssn
Moreover, if (i), (i), are satisfied, then
$\jmath^v(\lxlxzl{\eu T})={\lv Z}_\pel\varphi_v$ 
(all $v$). 
}
\endproclaim
\label\curvecase\theoremnum*

\demo{Proof}
Clear.
\enddemo

\proclaim{Definition} {\rm
In the context of Fact above, suppose that the equivalent
conditions (i), (ii) are satisfied. Then we define 
$\mardv_\clg$ to be any of the lattices
$\mardv_{{\eu T},\pel}\subset\widehat\mardv_\clg$,
and call it an {\it abstract divisor group\/} of 
${\cal G}$.

\ssn

We further say that $\mardv^0_{{\eu T},\pel}$ is the 
divisor group of degree $0$ in $\mardv_{{\eu T},\pel}$.
And remark that any two abstract divisor groups
$\mardv_{\cal G}$ and $\mardv'_{\cal G}$ are 
$\ell$-equiv\-alent latices in $\widehat\mardv_{\cal G}$,
and the same is true for $\mardv^0_{\cal G}$ and
${\mardv'}^0_{\cal G}$.
}
\endproclaim

{\bf The case:} $\delta>0$. 

\ssn

We begin by mimicking the construction from the case
$\delta=0$, and by induction on $\delta$ conclude the
construction. Thus let ${\eu T}=(\tau_v)_v$ be any system 
generators for the inertia groups $\ttv$ (all $v$). Further 
let ${\cal F}_{\eu T}$ be the Abelian pro-$\ell$ free group 
on the system ${\eu T}$ (written multiplicatively). Then 
$\ttt$ is a quotient ${\cal F}_{\eu T}\,\to\,\ttt\to1$. Thus 
taking $\ell$-adic Homs we get an exact sequence
$$
0\to\hhhom{\ttt}{\calz}\to
       \hhhom{{\cal F}_{\eu T}}{\calz}\,,
$$              
and remark that $\hhhom{{\cal F}_{\eu T}}{\calz}\,\cong\,
     \widehat{\rm Div}_{\cal G}$
in a canonical way. Next consider the system 
${\eu B}=(\varphi_v)_v$ of all the functionals 
$\varphi_v\in\hhhom{{\cal F}_{\eu T}}{\calz}$ 
defined by $\varphi_v(\tau_w)=1$ if $v=w$, and $0$, 
otherwise (all $v,w$). We denote by 
$$
{\rm Div}_{{\eu T},\pel}=\;<{\eu B}>_{(\ell)}\;
       \subset\hhhom{{\cal F}}{\calz}
$$ 
the $\zetell$-submodule of $\hhhom{{\cal F}}{\calz}$ 
generated by ${\eu B}$. Then ${\eu B}$ {\it is an 
$\ell$-adic basis\/} of $\hhhom{{\cal F}}{\calz}$, 
i.e., ${\rm Div}_{{\eu T},\pel}$ is $\ell$-adically 
dense in $\hhhom{{\cal F}}{\calz}$, and there are 
no non-trivial $\ell$-adic relations between the 
elements of ${\eu B}$. We will call 
${\eu B}=(\varphi_v)_v$ the ``dual basis'' to 
${\eu T}$. Thus finally, ${\rm Div}_{{\eu T},\pel}$ 
is a lattice in $\hhhom{\ttt}{\calz}$.  

\ssn

Now let ${\eu T}'=(\tau'_v)_v$ be another system
of inertia generators. And suppose that for some 
$\epsilon\in\calz$ we have ${\eu T}'={\eu T}^\epsilon$. 
If ${\eu B}=(\varphi'_v)_v$ is the dual basis to 
${\eu T}'$, then $\epsilon\,\varphi'_v=\varphi_v$ 
inside $\hhhom{\ttt}{\calz}$. Thus 
$\epsilon\cdot{\eu B}'={\eu B}$. In other words, 
${\eu B}$ and ${\eu B}'$ are $\ell$-adically 
equivalent. And we have: 
${\rm Div}_{{\eu T},\pel}
       =\epsilon\cdot{\rm Div}_{{\eu T}',\pel}$. 

\ssn

Finally, we fix notations as follows: 

\ssn

For a system of inertia generators ${\eu T}$ as above, and 
the corresponding lattice ${\rm Div}_{{\eu T},\pel}$ in 
$\widehat{\rm Div}_{\cal G}$, let $\lxlxzl{\eu T}$ be its 
pre-image in $\whlx_{\cal G}$. Next, for every $v$ we 
denote by ${\cal G}_v$ the corresponding $v$-residual 
\abstr\ \cycellprefrm, etc.. In particular, 
${\cal G}_v$ is both $(\delta-1)$-residually 
complete curve like, and ample up to level
$(\delta-1)$. Now let $\Delta$ be a finite co-rank 
$\calz$-submodule of $\whlx_{\cal G}$ 
such that $\Delta\cap\widehat U_{\cal G}=1$. For 
every $v$ such that $\Delta\subset\widehat U_v$, 
we set $\Delta_v=\jmath_v(\Delta)$. Since ${\cal G}$ 
is {\it ample up to level\/} $\delta$, there exists 
some $v$ such that:  

\msn

$(*)$ \ $\Delta\subset\widehat U_v$, and $\Delta_v$ has 
finite co-rank in $\whlx_{{\cal G}_v}$, and 
$\Delta_v\cap\widehat U_{{\cal G}_v}=1$. In particular,
$\Delta$ and $\Delta_v$ have the same $\calz$-rank.

\proclaim{Fact} {\rm 
For every $v$ and the corresponding ${\cal G}_v$,  
let a $\widehat U_{{\cal G}_v}$-lattice $\lxlxzl{v}$ 
in $\whlx_{{\cal G}_v}$ be given. 
Then up to $\ell$-adic equivalence, there exits at 
most one $\widehat U_{\cal G}$-lattice $\lxzl$ in 
$\whlx_{\cal G}$ such that for every finite 
co-rank $\calz$-module $\Delta\subset\whlx_{\cal G}$ 
with $\Delta\cap\widehat U_{\cal G}=1$ and $v$ as at 
$(*)$ above the following hold: 

\ssn

{\rm i)} $\LX_\Delta:=\Delta\cap\lxzl$, 
$\LX_{\Delta_v}:=\Delta_v\cap\lxlxzl{v}$ are lattices in 
$\Delta$, respectively $\Delta_v$.

\ssn

{\rm ii)} The lattices $\jmath_v(\LX_\Delta)$ and 
$\LX_{\Delta_v}$ are $\ell$-adically equivalent in $\Delta_v$. 

\ssn
Moreover, if the $\widehat U_{\cal G}$-lattice $\lxzl$
exists, then its $\ell$-adic equivalence class does 
depend only on the $\ell$-adic equivalence classes of 
the $\widehat U_{{\cal G}_v}$-lattices $\lxlxzl{v}$
(all $v$).
}
\endproclaim
\label\gencase\theoremnum*

\demo{Proof} 
Let $\lxzl,\lxzl'$ be $\widehat U_{\cal G}$-latices in 
$\whlx_{\cal G}$ satisfying i), ii), above. For $\Delta$ 
as in Fact above, set $\LX'_\Delta=\Delta\cap\lxzl'$. 
Then by hypothesis ii) it follows that both latices 
$\jmath_v(\LX_\Delta)$ and $\jmath_v(\LX'_\Delta)$
are $\ell$-adically equivalent to the lattice
$\LX_{\Delta_v}$ inside $\Delta_v$. Thus they are 
equivalent. After replacing $\lxzl'$ by some multiple
$\epsilon\cdot\lxzl'$ with $\epsilon\in\Zell^\times$, we 
can suppose that $\LX_\Delta=\LX'_\Delta$. As $\Delta$ 
was arbitrary, it finally follows that $\lxzl=\lxzl'$. 
\enddemo 

\proclaim{Definition} 
{\rm 
Let ${\cal G}$ be a \abstr\ \cycellprefrm\ which is 
both $\delta$-residual complete curve like and ample up 
to level $\delta>0$. We define an abstract divisor group 
$\mardv_{\cal G}\subset\widehat\mardv_{\cal G}$ ---if it 
exists--- inductively as follows: First, suppose that 
for all $v$, an abstract divisor group $\mardv_{{\cal G}_v}$
for ${\cal G}_v$ does exist. In particular, we have 
by definition, see hypothesis i) below: 
$\mardv_{{\cal G}_v}$ is a lattice in 
$\widehat\mardv_{{\cal G}_v}$, and its pre-image 
$\lxlxzl{v}$ via $\jmath^{{\cal G}_v}$ is a 
$\widehat U_{\clg_v}$-lattice in $\whlx_{\clg_v}$. 

\ssn

Then an {\it abstract divisor group\/} of ${\cal G}$ 
---if it exists--- is any lattice of the form 
$\mardv_\clg={\rm Div}_{{\eu T},\pel}$ 
in $\widehat\mardv_\clg$ such that its pre-image 
$\lxzl$ in $\whlx_{\cal G}$ satisfies the 
following conditions: 

\ssn

j) $\lxzl $ is a $\widehat U_\clg$-lattice in 
$\whlx_\clg$.

\ssn

jj) $\lxzl$ satisfies the conditions i), ii) from 
Fact~\gencase. 
} 
\endproclaim
\label\abstrdiv\theoremnum*
\proclaim{Remark/Fact}
{\rm 
Let ${\cal G}$ be a \abstr\ \cycellprefrm\ which is 
both $\delta$-residual complete curve like and ample up 
to level $\delta>0$. Suppose that abstract divisor 
groups $\mardv_\clg={\rm Div}_{{\eu T},\pel}$ for 
${\cal G}$ do exist. For such an abstract divisor group 
$\mardv_\clg$, let $\lxzl $ be its pre-image in $\whlx_\clg$. 
Then for all $v$ one has:
$$
\jmath^v(\lxzl )={\lv Z}_\pel\,\varphi_v\,.
$$
Indeed, by condition i) of the ampleness, see 
Definition/Remark \abstractprefields,~7), it follows that 
$\jmath^v(\whlx_\clg)=\Zell\,\varphi_v$. Further, 
since $\lxzl$ is $\ell$-adically dense in $\whlx_\clg$, 
it follows that $\jmath^v(\lxzl)$ is dense in 
$\Zell\,\varphi_v$. Thus the assertion. 

\ssn

2) In particular, the $\widehat U_\clg$-lattice $\lxzl$ 
determines $\mardv_\clg$, as being the additive subgroup
$$
\mardv_\clg=\oplus_v\,\jmath^v(\lxzl )
$$
Therefore, giving an abstract divisor group $\mardv_\clg$,
is equivalent to giving a $\widehat U_\clg$-lattice $\lxzl$ 
in $\whlx_\clg$ such that:

\ssn

j) $\lxzl$ satisfies the conditions i), ii) from 
Fact~\gencase\ with respect to the pre-images $\lxlxzl v$
of some abstract divisor groups $\mardv_{\clg_v}$ (all $v$). 

\ssn

jj) $\jmath^v(\lxzl)\cong{\lv Z}_\pel$ (all $v$).

\msn

3) Finally, for an abstract divisor group $\mardv_\clg$ for 
${\cal G}$, and its pre-image $\lxzl$ in $\whlx_\clg$, we set 
$\,{\eu Cl}\,_{\lxzl}=\mardv_\clg/\jmath^{\cal G}(\lxzl)$, 
and call it the abstract ideal class group of $\lxzl$. 
Thus one has a commutative diagram of the form

\ssn

$$
\matrixbaselines
\matrix{
0\to\>\widehat U_{\cal G}
    &\hookrightarrow
       &\lxzl
          &\hor{\jmath^{\cal G}}
             &{\rm Div}_{\lxzl}
                &\hor{\rm can}&{\eu Cl}\,_{\lxzl}\>\to0\cr
 \hhb{20}\dwn{}&&\dwn{}\hhb5
            &&\dwn{}\hhb5&&\dwn{\jmath_{\eu Cl}}\hhb{32}\cr
 0\to\>\widehat U_{\cal G}
    &\hookrightarrow
       &\whlx_{\cal G}
          &\hor{\jmath^{\cal G}}
             &\widehat{\rm Div}_{\cal G}
                &\hor{\rm can}&\widehat{\eu Cl}_{\,\cal G}\>
                                                      \to0\cr
}
\eqno{\indent(*)}
$$
where the first three vertical morphisms are the canonical
inclusions, and the last one is the $\ell$-adic 
completion homomorphism.
}
\endproclaim
\label\abstractfd\theoremnum*
\proclaim{Proposition}
Let ${\cal G}$ be a \abstr\ \cycellprefrm\ which is 
both $\delta$-residual complete curve like and ample up 
to level $\delta>0$. Then any two abstract divisor groups 
$\mardv_\clg$ and $\mardv'_\clg$ for ${\cal G}$ are 
$\ell$-{\it adically equivalent\/} as lattices
in $\widehat\mardv_\clg$. Or equivalently, their
pre-images $\lxzl$ and $\lxzl'$ in $\whlx_{\cal G}$
are $\ell$-adically equivalent $\widehat U_\clg$-lattices
in $\whlx_{\cal G}$.
\endproclaim
\label\abstractdiv\theoremnum*
\demo{Proof}
We prove this assertion by induction on 
$\delta$. For $\delta=0$, the uniqueness is
already shown, see Fact~\curvecase, and Definition
in case $\delta=0$. Now suppose that $\delta>0$. Let 
$\mardv_{{\cal G}_v}$ and $\mardv'_{{\cal G}_v}$ be abstract
divisor groups for ${\cal G}$ used for the definition
of $\mardv_\clg$, respectively $\mardv'_\clg$ (all $v$).
By the induction hypothesis, $\mardv_{{\cal G}_v}$ 
and $\mardv'_{{\cal G}_v}$ Thus 
their pre-images $\lxlxzl v$ and $\lxlxzl v'$ in 
$\widehat{\cal L}_{{\cal G},v}$ are $\ell$-adically 
equivalent $\widehat U_{{\cal G},v}$-lattices. 
Therefore, by Fact~\gencase, the lattices $\lxzl$ 
and $\lxzl'$ (which are the pre-images of $\mardv_\clg$ 
respectively $\mardv'_\clg$ in $\whlx_{\cal G}$) 
are $\ell$-adically equivalent. Finally, use 2) above
to conclude.
\enddemo
\proclaim{Definition}
{\rm
A level $\delta>0$ {\it divisorial \cycellprefrm\/} is
by definition every \abstr\ \cycellprefrm\ ${\cal G}$ 
which is both $(\delta-1)$-residually complete curve like 
and ample up to level $(\delta-1)$ and has abstract 
divisor groups $\mardv_\clg$. If this is the case, we will 
denote by $\lxzl$ the pre-image of $\mardv_\clg$ in 
$\widehat{\cal L}_{\cal G}$, and call it a {\it divisorial 
$\widehat U_{\cal G}$-lattice\/} in $\widehat{\cal L}_{\cal G}$.
}
\endproclaim

C) {\it Example: Geometric \cycellprefrm\/}

\ssn
Let $\ell$ be a prime number. Let $\KK|\kpkp$ be a field 
extension with $\kpkp$ an algebraically closed field, 
${\rm char}(\kpkp)\neq\ell$. Consider a maximal Abelian 
pro-$\ell$ field extension $\prm\KK=\KK^{\ell{\rm ab}}$
of $\KK$, and set $\ggg={\rm Gal}(\prm\KK|\KK)$. Then we get:

\proclaim{Fact} {\rm
${\cal G}=\big(\ggg, \>\vid\big)$ is a level $\delta=0$ 
\abstr\ \cycellprefrm.}
\endproclaim

Nevertheless, in the case $\KK|\kpkp$ is a {\it function 
field\/} with $d=\td(\KK|\kpkp)>0$, we can refine the pre-Galois 
formation above by making it into a \abstr\ \cycellprefrm\ 
of level $\delta$ for every $0\leq\delta\leq d$ as follows: 
First, let us endow every function field $\KK|\kpkp$ as 
above with a geometric set ${\eu D}={\eu D}_X$ of Zariski 
prime divisors of $\KK|\kpkp$. Here $X\to\kpkp$ is a 
quasi-projective normal model of $\KK|\kpkp$, and 
${\eu D}_X$ is its set of Weil prime divisors. Second,
since $\prm\KK|\KK$ is Abelian, if $v'$ is a prolongation
of $v\in{\eu D}$ to $\prm K$, then $T_{v'}\subset Z_{v'}$
do not depent on $v'$, but on $v$.

\ssn

Then one can endow $\ggg$ with the family $(\zzv)_v$
of the decomposition groups of all the $v\in{\eu D}$.
Clearly, if $\ttv\subset\zzv$ is the inertia group at
$v$, then $\ttv\cong\Zell$ is pro-$\ell$ cyclic. 
Further, if $v\neq w$, then $\ttv\cap\ttw=\{1\}$. Finally, 
let $U_i$ be a basis of (small enough) affine open subsets
in $X$, and for every $i$, set ${\eu U}_i=U_i\cap{\eu D}$.
Then by general facts about geometric fundamental groups,
it follows that $({\eu U}_i)_i$ is a co-final system 
in the sense of Axiom II),~ii), from Subsection~A).
Thus finally, we have a \abstr\ \cycellprefrm\ supported 
by $\ggg$ as follows:
$$
{\cal G}_{\eu D}:=\big(\ggg, (\zzv)_v\big)=:{\cal G}_X
$$
Moreover, for every $v\in{\eu D}$ we have: The 
residue field $\KK v$ is again a function field over 
$\kpkp$ with $\td(\KK v\,|\,\kpkp)=(d-1)$. Further,
by general decomposition theory, the residue field 
extension $\prm\KK v\,|\,\KK v$ is a maximal 
Abelian pro-$\ell$ extension of $\KK v$. 

\ssn

Thus if $d>1$, each such function field $\KK v|\kpkp$ 
comes equipped with a geometric set of Zariski prime 
divisors ${\eu D}_v={\eu D}_{X_v}$, where $X_v$ a 
quasi-projective normal model of $\KK v\,|\,\kpkp$. 

\ssn

Now consider any $\delta$ with $0\leq\delta\leq d$. Then  
by induction on the transcendence degree over $\kpkp$, 
we can suppose that each $\ggv$ endowed with the set of 
decomposition groups defined by ${\eu D}_v$ is a \abstr\
\cycellprefrm\ ${\cal G}_{{\eu D}_v}$ of level $(\delta-1)$. 
Thus we have the following: 

\proclaim{Proposition/Definition} Let
${\cal G}_{\eu D}=\big(\ggg,(\zzv)_v\big)$ be as constructed 
above. Then ${\cal G}_{\eu D}$ is in a canonical way a 
level $\delta$ \abstr\ \cycellprefrm, for every $\delta$ 
satisfying $0\leq\delta\leq\td(\KK|\kpkp)$. 

\ssn

A \abstr\ \cycellprefrm\ of the form ${\cal G}_{\eu D}$ 
will be called a {\rm geometric \cycellprefrm} of level 
$\delta$. 
\endproclaim
\label\theexample\theoremnum*

\proclaim{Remarks} {\rm
Let ${\cal G}_X=\big(\ggg,(\zzv)_v\big)$ be a geometric
\cycellprefrm\ as constructed/defined above. Then 
$\widehat\KK=\hhhom{\Gprm\KK}{\calz}=\whlx_{{\cal G}_X}$
by Kummer Theory. In order to compute 
$\widehat U_{{\cal G}_X}$, $\widehat{\rm Div}_{{\cal G}_X}$, 
and $\widehat{\eu Cl}_{{\cal G}_X}$, we do the following:
First, let ${\cal H}(X)$ denote the group of principal 
divisors on $X$, and consider the canonical exact sequence 
$0\to {\cal H}(X)\hor{}{\rm Div}(X)\hor{pr}
                                  {\eu Cl}\,(X)\to 0\,.$  
\ssn

1) Passing to $\ell$-adic completions, we get an exact 
sequence of $\ell$-adic complete groups of the form:
$0\to{}{\cal T}_{\ell,X}\hor{}\widehat{\cal H}(X)\hor{}
   \widehat{\rm Div}(X)\hor{}\widehat{\eu Cl}\,(X)\to 0$, 
where ${\cal T}_{\ell,X}=\plim{}\,{}_{\ell^n}{\eu Cl}\,(X)$ 
(with multiplication as homomorphisms), and the last 
three objects the corresponding $\ell$-adic completions. 

\ssn

2) On the other hand one has
${\rm Div}(X)=\oplus_{{v}\in{\eu D}_X}\,v\KK$. 
For every $v$, consider the commutative diagram from 
subsection A),~3). We get a commutative diagram of the 
following form:
$$
\matrixbaselines
\matrix{
&\widehat\KK 
 &\hhb{-5}\shto&\hhb{-20}\widehat{\rm Div}(X)=\widehat{\oplus_v}\,v\KK
                      &\shto\hhb5&\widehat{\eu Cl}\,(X)&\to\;0&\cr
&\dwn{\hat\delta}&&\hhb{35}\dwn{\oplus\theta^v}&
                                      &\dwn{}&&\cr
&\hhhom{\Gell\KK}{\Zell}
&\shorrr{\jmath^{{\cal G}_X}}&\widehat{\oplus_v}\,
 \hhhom{\ttv}{\Zell}&\;\shor{\rm can}&
                                \widehat{\cal P}_{X}&\to\;0&\cr
}
\leqno{}
$$
where the vertical maps are isomorphisms, and 
$\widehat{\cal P}_{X}$ is simply the quotient of the 
middle group by the second. We remark/recall that 
${\oplus}\theta^v$ is defined as follows: For every 
$v$, let $\gamma_v=1\cdot v$ be the unique positive 
generator of $v\KK={\lv Z}v$. Then there exists a unique 
generator of the arithmetical inertia $\tau_v\in{\eu T}_v$ 
such that $\jmath_v(\gamma_v)\big(\tau_v\big)=1$. 
Hence, the commutativity of the diagram above follows 
directly from the definition of the homomorphisms.

\ssn

3) From this we deduce: 
$\widehat{\rm Div}_{{\cal G}_X}=\widehat{\rm Div}(X)
    \ \hbox{ and } \  
       \widehat{\eu Cl}\,_{{\cal G}_X}=\widehat{\eu Cl}(X)$.

\ssn

4) Concerning $\widehat U_{{\cal G}_X}$, we recall that 
$\widehat U_{{\cal G}_X}=\hhhom{\pi^{\ell,{\rm ab}}_1(X)}{\Zell}$,
see Defin\-ition/Remark 4.2,~2). Let $\ttt\subset\ggg$ be 
the closed subgroup generated by the decomposition groups 
of all the divisorial valuations of $\KK$. We set 
$\pi^{\ell,\rm ab}_{1,\KK}=\ggg/\ttt$, and call it the 
{\it fundamental group to $\KK$,\/} as it is a birational 
invariant of $\KK$. It equals the fundamental group of any 
complete regular model $X_0\to\kpkp$ of $\KK|\kpkp$, if 
any such do exist...

\ssn

We will say that a normal model $X\to\kpkp$ 
for $\KK|\kpkp$ is {\it regular complete like\/} if 
$\pi_1^{\ell,\rm ab}(X)=\pi^{\ell,\rm ab}_{1,\KK}$.   
From the structure theorem for (the Abelian pro-$\ell$ 
quotient of the) fundamental groups of a normal curve 
it follows that such a curve is regular complete like 
if and only if it is a complete one.  
}
\endproclaim
\label\moreexample\theoremnum*
\proclaim{Proposition} 
Let $\KK|\kpkp$ and a geometric \cycellprefrm\ ${\cal G}_X$ 
be as introduced/defined in Proposition above. Consider some 
$\delta\leq\td(\KK|\kpkp)$, and let us view ${\cal G}_X$ 
as a \abstr\ \cycellprefrm\ of level $\delta$. We denote
by ${\cal G}_{X_{\bfv}}$ the residual \cycellprefrm s of
${\cal G}_X$ with $\bfv$ multi-indices of length $\leq\delta$. 
In particular, $X_{\bfv}\to\kpkp$ are quasi-projective normal 
varieties. Then one has:

\ssn

{\rm(1)} ${\cal G}_X$ is $\delta$-residually complete curve 
like if (and only if) $\delta=\td(\KK|\kpkp)-1$, and all the 
$\delta$-residual varieties $X_{\bfv}$ are complete normal 
curves. 

\ssn

{\rm(2)} Suppose that $\delta<\td(\KK|\kpkp)$. Further 
suppose that $X_{\bfv}$ is regular complete like (all 
$\bfv$ as above). Then ${\cal G}_X$ is ample up to level 
$\delta$. 

 \ssn

{\rm(3)} Suppose that ${\cal G}_X$ is $\delta$-residually 
curve like, and ample up to level $\delta$. Then 
${\rm Div}(X)_\pel:={\rm Div}(X)\otimes{\lv Z}_\pel$ inside 
$\widehat{\rm Div}(X)$ is an abstract divisor group of 
${\cal G}_X$. We will denote by $\lxlxzl{X}$ its pre-image 
in $\widehat\KK$.

\ssn

Thus for every geometric \cycellprefrm\ ${\cal G}_X$ 
as above, ${\rm Div}(X)_\pel$ is an abstract divisor group
of ${\cal G}_X$, and   $\lxlxzl{X}$ is a divisorial 
$\widehat U_{\cal G}$-lattice of ${\cal G}_X$.

\endproclaim
\label\geometric\theoremnum*

\demo{Proof} To (1): Clear. It is nevertheless 
more/quite difficult to prove the ``only if'' part of~(1), 
which we will not directly use, thus omit the proof here... 

\ssn

To (2): We make induction on $d=\td(\KK|\kpkp)$. In the 
notations from Definition/Remarks \abstractprefields,~3), 
let $\Delta\subset\widehat\KK$ be a co-finite rank 
submodule such that $\Delta\cap\widehat U_{{\cal G}_X}=0$. 
Then there exists an open affine subset $X'\subset X$ 
such that for all $v\in{\eu D}_{X'}$ one has: 
$\Delta \subset\widehat U_v$. In particular, the canonical
projection $\pi_1(X')\to\pi_1(X)$ gives rise to an 
embedding $\widehat U(X)\hookrightarrow\widehat U(X')$, and 
$\Delta\subset\widehat U(X')$. Using \nmnm{de Jongs}'s
alterations, and the inclusion-norm maps, we can suppose
that $X$ is a smooth $\kpkp$-variety. Finally, let 
$X_1,\dots,X_n$ be the finitely many Weil prime divisors 
in $X\backslash X'$, and $S:=\{x_1,\dots,x_n\}$ their
generic points. Thus the inertia groups at these finitely 
many points generate 
$T_{X',X}:=\ker\big(\pi_1(X')\to\pi_1(X)\big)$. 

\ssn

Now let $\imath:X\to{\lv P}^N_\kpkp$ be some $\kpkp$-embedding. 
For a general hyper-plane $H\subset{\lv P}^N_\kpkp$, we set 
$Y=H\cap X$, $Y'=H\cap X'$, $Y_i=X_i\cap H$ (all $i$). Then
by Bertini, each $Y_i$ is a prime divisor of $Y$. Let 
$T:=\{y_1,\dots,y_n\}$ is the set of their generic points
of the $Y_i$ (all $i$). Then by Bertini we have:

\ssn

i) The canonical projections $\pi_1(Y)\to\pi_1(X)$ and  
$\pi_1(Y')\to\pi_1(X')$ are surjective.

\ssn

ii) Set $T_{Y',Y}:=\ker\big(\pi_1(Y')\to\pi_1(Y)\big)$. Then
$T_{Y',Y}$ is generated by the inertia groups at all the $y_i$ 
(all $i$). 

\ssn

iii) Finally, $T_{Y',Y}$ is mapped surjectively onto $T_{X',X}$ 
under the projection $\pi_1(Y')\to\pi_1(X')$ from i) above.

\ssn

This is now exactly the translation of the fact that for
the Zariski prime divisor $v$ of $\KK$ defined by the
Weil prime divisor $Y=H\cap X$ of $X$, the assertion from
Definition/Remarks \abstractprefields,~3), holds for 
$\Delta$ at $v$.

\msn

To (3): It follows immediately from (1) and (2) above.
\enddemo

\ssn

D) {\it The case $k$ is the algebraic closure
               of a finite field\/}

\ssn

In this subsection we will discuss the case $K|k$ 
is a function field with $\td(K|k)>1$ and $k$ an 
algebraic closure of a finite field. Let 
$\Gell K={\rm Gal}(\fell K|K)$ be as usually defined. 
The aim of this subsection is to show that the system of 
all the geometric \cycellprefrm s ${\cal G}_X$ are actually 
group theoretically encoded in $\Gell K$. Moreover, the 
extra information concerning such a \abstr\ \cycellprefrm, 
e.g., $\delta$-residually complete curve like, and/or 
ample up to level $\delta$, is also encoded in the 
$\Gell K$. First some general preparation as follows: 
\proclaim{Fact} 
{\rm
Let $\KK|\kpkp$ be some field extension with $\kpkp$ 
algebraically closed, ${\rm char}(\KK)\neq\ell$. Let
$\fell\KK$ and $\Gell\KK$ as usually define. Denote by 
$$
\Gell\KK\to\abab{G_\KK}=:G
$$
the canonical projection. For every subgroup $(\cdot)^\ell$
of $\Gell\KK$ we will denote by $\abab{(\cdot)}$ the image 
of $(\cdot)$ in $\abab{G_\KK}$. By general decomposition 
theory we have:

\ssn

1) For every valuation $v^\ell$ of $\KK^\ell$, 
and its inertia, respectively decomposition groups 
$T^\ell_v\subset Z^\ell_v$ in $\Gell\KK$ one has: 
$\abab{T_v}$ and $\abab{Z_v}$ are exactly the inertia, 
respectively decomposition groups of $v^\ell$ in 
$\abab{G_\KK}$. 

\ssn

Second, the groups of the form $\abab{T_v}$ and 
$\abab{Z_v}$ are exactly all the inertia, respectively 
decomposition groups in $\abab{G_\KK}$.

\ssn

2) Further remark that $\abab{T_v}$ and $\abab{Z_v}$ do 
depend only on the restriction $v$ of $v^\ell$ on $\KK$, 
and not on $v^\ell$ itself (which is one of the prolongations 
of $v$ to $\KK^\ell$). 

\ssn

3) \ $T^\ell_v\to\abab{T_v}$ is an isomorphism, and 
second, $\abab{G_{Kv}}$ is canonically isomorphic to 
$\abab{Z_v}/\abab{T_v}$.

\ssn

4) The Theorem of F.~K.\ Schmidt from \nmnm{Pop} [P], 
Proposition~1.3, holds for the Galois extension $\abab\KK|\KK$. 
In particular, if two valuations $v$ and $w$ equal their cores 
respectively, then $\abab{Z_v}\cap\abab{Z_w}\neq\{1\}$ if and 
only if $v=w$. 
}
\endproclaim
\label\generalab\theoremnum*
\proclaim{Fact} \
{\rm Let $K|k$ be a function field with $k$ the algebraic
closure of a finite field. {\it In order to simplify notations\/} 
let us denote $G:=\abab{G_K}$, and for every valuation
$v^\ell$ and $\ttv^\ell\subset\zzv^\ell$ as above, set
$\ttv=\abab\ttv$ and $\zzv=\abab\ttv$. Then we have the 
following:

\ssn

1) $d=\td(K|k)$ is encoded in $\Gell K$, by Section 2,~A).

\ssn

2) The divisorial inertia groups $\ttv$ and the divisorial
decomposition groups $\zzv$ in $G$ are known, by Section~1.

\ssn

The same is true for all the defectless inertia elements,
by Section 1), Theorem~1.11.    

\ssn

For given such subgroups $\zzv$ and $Z_w$ we have: 
$\zzv\cap Z_w$ is non-trivial if and only if $v=w$.

\ssn

3) If $\pi^\ell_{1,K}=\Gell K/T^\ell_K$ is the fundamental
group to $K$, then  $\pi_{1,K}:=\pi^{\ell,{\rm ab}}$ is also 
encoded in $\Gell K$, by 2) above. 

\ssn

4) Furthermore, the geometric sets of divisorial decomposition
groups, say ${\cal Z}=\{\zzv\}_v$ of $G$ are encoded in 
$\Gell K\to G$, see Section 2,~C).

\ssn

5) Let ${\eu D}={\eu D}_X$ be a geometric set of Zariski 
prime divisors, say defined by a set ${\cal Z}=\{\zzv\}_v$ 
of decomposition groups. Then the fact that $X$ is complete
regular like is encoded in $\Gell K$ and the set 
${\cal Z}=\{\zzv\}_v$.    
}
\endproclaim

Therefore we finally have the following: 

\proclaim{Proposition}
Let $K|k$ be a function field with $\td(K|k)>1$, and $k$ 
the algebraic closure of a finite field. Then in the 
notations from above we have the following: The geometric
\cycellprefrm s ${\cal G}_{\eu D}$, say ${\eu D}={\eu D}_X$, 
of some level $\delta\leq\td(K|k)$ are encoded in $\Gell K$ 
along the lines indicated above.

\ssn

Moreover, the fact that such a geometric \cycellprefrm\ 
${\cal G}_{\eu D}$ is  $\delta$-residual\-ly complete curve 
like and/or that all residual geometric \cycellprefrm s
${\cal G}_{\bfv}$ are regular complete like (all $\bfv$)
is also encoded in the group theoretic information carried
by ${\cal G}$. 

\ssn

Let ${\cal G}_{\eu D}$ be a geometric \cycellprefrm\  
which is $\delta$-residually complete curve like and 
residually geometric complete like up to level $\delta$.
Then one has: $\delta=\td(K|k)-1$, and viewing ${\cal G}_{\eu D}$ 
as a \abstr\ \cycellprefrm, we have: ${\rm Div}(X)_\pel$ is 
a representative for its divisorial lattices, and its pre-image 
$\lxlxzl{X}$ in $\widehat K$ is a representative for the 
divisorial $\widehat U_X$-lattices in $\widehat K$. 
\endproclaim
\label\charprop\theoremnum*
\demo{Proof}
Let ${\cal G}=\big(\ggg,(\zzv)_v\big)$ be a \abstr\
\cycellprefrm\ of level $\delta$ on $\ggg=\abab{G_K}$.
Then  ${\cal G}$ is a geometric \cycellprefrm\ of level 
$\delta\leq\td(K|k)=:d$ if an only if by induction the 
following hold:

\ssn

i) $\zzv$ is a divisorial decomposition group as given 
by the Fact above  for all ``abstract valuations'' $v$,
and $\ttv$ is the inertia subgroup of $\zzv$.

\ssn

ii) The set ${\eu D}$ of all the valuations $v$ of ${\cal G}$
is a geometric set of Zariski prime divisors.

\ssn

(iii) By induction on $d$, each residual \abstr\ \cycellprefrm\
${\cal G}_v$ on $G_v=\Gell{Kv}$ is a geometric \cycellprefrm\ 
of level $\delta-1$. Remark that in the limit case when $d=2$, 
all the inertia groups $T_a\subset\Gell{Kv}$ are encoded in 
$\Gell K$, see Fact~\dequalstwo, thus ``known''.

\msn

Now suppose that ${\cal G}=\big(\Gell K,(Z_v)_v\big)$ 
a geometric \cycellprefrm\ of level $\delta\leq\td(K|k)$ 
is given. Then ${\cal G}$ is {\it $\delta$-residually 
complete curve like\/} if and only if $\delta=d-1$; and 
further for all the $\delta$-residual \cycellprefrm s 
${\cal G}_{\bfv}=\big(G_{\bfv},(Z_{\bfv})_{\bfv}\big)$ 
one has: $(Z_{\bfv})_{\bfv}$ is exactly the family of 
all the inertia groups of $\ggg_{\bfv}$. Remark that by 
induction on the length of the multi-index $\bfv'$, and 
Fact~\dequalstwo, we know all the geometric inertia elements 
in $\ggg_{\bfv'}$. And that for a multi-index $\bfv$ of 
length $\delta$, a subgroup of $Z_{\bfv}$ is inertia group 
if and only if it is maximal among the pro-cyclic subgroups 
generated by inertia elements of $\ggg_\bfv$.

\msn

Finally, let ${\cal G}=\big(\Gell K,(Z_v)_v\big)$ be a 
geometric \cycellprefrm\ of level $\delta<d$ be given.
Let ${\cal G}_{\bfv}=\big(G_{\bfv},(Z_{\bfv})_{\bfv}\big)$ 
be the $\bfv$-residual geometric \cycellprefrm\ to some 
multi-index $\bfv$ of length $\leq\delta$. Denote by
$\pi_{1,K\!\bfv}$ the fundamental group attached to the 
residual function field $K\!\bfv|k$. Then by the Fact 
above, the fact that the normal model $X_\bfv\to k$ of 
$K\!\bfv|k$ is complete regular like if and only if  
$\pi_{1,K\!\bfv}=\pi_{1,{\cal G}_\bfv}$. Thus ${\cal G}$ 
is {\it ample up to level\/} $\delta$  if and only if 
$\pi_{1,K\!\bfv}=\pi_{1,{\cal G}_\bfv}$ (all multi-indeces 
$\bfv$ of length $\leq\delta$).
\enddemo

\vfill\eject

\section{Concluding the Proof of the 
                       Theorem {\rm (Introduction)}}

\ssn

A) \ {\it Detecting rational projections\/}

\ssn
 
First some definitions: Let $\KK|\kpkp$ be an arbitrary
function field, say with $\kpkp$ algebraically closed. 
For every non-constant function $x\in\KK$, let 
$\KK_x$ be the relative algebraic closure of $\kpkp(x)$ 
in $\KK$. Then $\KK_x|\kpkp$ is a function field in one 
variable. For a given Galois extension $\prm\KK|\KK$, 
the relative algebraic closure of $\KK_x$ in $\prm\KK$
will be denoted by $\prm{\KK_x}|\KK_x$. And the inclusion 
$\prm{\KK_x}\hookrightarrow\prm\KK$ gives rise to a 
surjection projection
$$
\prm p_x:\Gprm\KK\to\Gprm{\KK_x}
$$
which we then call the $1$-{\it dimensional projection 
defined by or attached to $x\in\KK$\/} (for the Galois
group $\Gprm\KK$). We will say that a $1$-dimensional 
projection $\prm p_x$ is a {\it rational projection,\/} 
if $K_x$ is a rational function field. 

\ssn

Our aim in this section is to show that in the case 
$K|k$ with $\td(K|k)>1$ and $k$ an algebraic closure of 
a finite field, the rational projections of $\abab{G_K}$ 
are group theoretically encoded in $\Gell K$.

\ssn

The strategy is as follows, compare with \nmnm{Bogomolov}
[Bo], Lemma~4.2. [But note also that in loc.cit.\ one 
gives/has a recipe for detecting the projections $p_x$ 
{\it under the hypothesis that $\jmath(K^\times)_\pel$ 
inside $\widehat K$ is known...\/} whereas we cannot 
work under this hypothesis: We are actually 
{\it looking for detecting $\jmath(K^\times)_\pel$ 
inside $\widehat K$.\/}] First, for $n=\ell^e$ a power 
of $\ell\neq{\rm char}$, consider the cup product:
$$
\psi_n:\KK^\times\!\!/n\times\KK^\times\!\!/n=
 \hhone{\Gell\KK}{{\lv Z}/n}\times\hhone{\Gell\KK}{{\lv Z}/n}
  \to{\rm H}^2(\Gell\KK,{\lv Z}/n)={}_n{\rm Br}\,(\KK).
$$
Clearly, the whole above picture is of pure group 
theoretical nature. Now \nmnm{Bogomolov}'s idea from 
loc.cit.\ is to use the following fact (which is clear
by Tsen's Theorem, as $\kpkp$ is algebraically closed):

\msn

{\it Let $x,y\in\KK$. Then $\KK_x=\KK_y$
if and only if $\psi_n(x,y)=0$ for all $n$.\/}

\msn

This means that if we know $\jmath_\KK(\KK^\times)$ inside
$\widehat\KK$, then we know $\jmath_\KK(\KK_x^\times)$ 
too. Now Kummer theory gives an $\ell$-adic perfect 
pairing $<,\!>:\widehat\KK\times\abab{G_\KK}\to\Zell$.
Thus the kernel of the projection 
$p_x:\abab{G_\KK}\to\abab{G_{\KK_x}}$, is nothing but
the orthogonal complement of $\jmath_\KK(\KK_x^\times)$:
$$
{\cal N}_x:=\ker(p_x)=\{\sigma\in\Gell\KK\mid<x',\sigma>=0
\hbox{ for all } x'\in\jmath_\KK(\KK_x^\times)\,\}.
$$
And this gives a description of all the kernels ${\cal N}_x$
of projections $p_x$. 

\msn

Now come back to the situation $K|k$ is a function
field with $\td(K|k)>1$, and $k$ the algebraic closure
of a finite field.

\ssn

Let ${\cal G}_X$ be a geometric \cycellprefrm, which is
complete curve like and ample up to level $\delta=\td(K|k)-1$.
Further let $\mardv_{X,\pel}$ be an abstract divisor group
for ${\cal G}_X$ and $\lxlxzl{X}$ its pre-image in $\widehat K$.
Then after multiplying $\mardv_{X,\pel}$ by an 
$\ell$-adic unit $\epsilon$, we can suppose that
$\mardv(X)_\pel=\mardv_{X,\pel}$. In particular, the 
divisorial $\widehat U_X$-lattice $\lxzl$ contains 
$K^\times_\pel:=\jmath_K(K^\times)\otimes{\lv Z}_\pel$
inside $\widehat K$. We now have the following
\nonumproclaim{Lemma 1} In the above context, 
$\lxlxzl{X}\,/\,\big(\widehat U_X\cdot K^\times_\pel\big)$ 
is a torsion group.
\endproclaim
\demo{Proof} Since $\lxlxzl{X}$ is the pre-image of
$\mardv_{X,\pel}$ in $\widehat K$, the contention of the 
lemma is equivalent to the following: 
$\jmath^{{\cal G}_X}(\lxlxzl{X})/{\rm div}(K^\times)_\pel$ 
is a torsion group. We have a commutative diagram of the form

$$
\matrixbaselines
\matrix{
 0\>\to\>&{\rm div}(K^\times)_\pel
          &\hookrightarrow
             &{\rm Div}(X)_\pel
                &\hor{\rm can}&{\eu Cl}(X)_\pel  &\to\>0\cr
\hhb{20}&\dwn{\rm incl}\hhb5 
          &  &\dwn{\rm id}\hhb5
                &             &\dwn{\rm can}\hhb5&\cr
 0\>\to\>&\jmath^{{\cal G}_X}(\lxlxzl{X})
          &\hookrightarrow
             &{\rm Div}(X)_\pel
                &\hor{\rm can}&{\eu Cl}_{X,\pel} &\to\>0\cr
}
\eqno{\indent(*)}
$$

\msn
where the first vertical map is an inclusion. So by the 
Snake lemma it follows that 
$\jmath^{{\cal G}_X}(\lxlxzl{X})/{\rm div}(K^\times)_\pel\cong
\ker\big({\eu Cl}(X)_\pel\to{\eu Cl}_{X,\pel}\big)$.
Thus we have to show that this last group is a torsion
one.

\ssn

Now since $k$ is the algebraic closure of a finite field, 
one has an exact sequence of the form:
$$
0\to\hbox{(divisible group)}_\pel\to
             {\eu Cl}(X)_\pel\to{\rm NS}(X)_\pel\to 0
\leqno(*)
$$
where ${\rm NS}(X)$ is a finitely generated Abelian group.
(Actually, if $X$ is a smooth model, the assertions above 
are all well known. If not, then consider a smooth 
alteration $\tilde X$ of $X$, and conclude by using the 
inclusion-norm map.) 

\ssn

And remark that the ${\lv Q}$-rank of ${\rm NS}(X)_\pel$ 
is the same as its ${\lv Q}_\ell$-rank. And further, that
this last rank is the same as the ${\lv Q}_\ell\,$-rank of 
$\widehat{\eu Cl}_{X,\pel}$. Thus finally equal to the 
${\lv Q}$-rank of ${\eu Cl}_{X,\pel}$. Therefore, 
$\ker\big({\eu Cl}(X)_\pel\to{\eu Cl}_{X,\pel}\big)$
is a torsion group, as contented. 
\enddemo
\nonumproclaim{Corollary} 
For every $\tilde x\in\lxlxzl{X}$ there exists some 
$u\in\widehat U_X$ and $x\in\kpel$ such that $x^m=u\cdot x$ 
for some integer $m>0$. 

\ssn

In other words we have the following: $\lxlxzl X$ is exactly
the relative divisible hull of $\widehat U_X\cdot\kpel$ 
in $\widehat K$, i.e., one has:
$$
\lxlxzl X=\{\tilde x\in\widehat K\mid \exists\, m
 \hbox{\rm\ such that }\tilde x^m\in\widehat U_X\cdot\kpel\,\}
$$
\endproclaim
\nonumproclaim{Notations/Remark}
{\rm
\ssn

We set $\kxzk'_{X,\pel}=\widehat U_X\cdot\kpel$. Further, 
for $x\in\kpel$ we denote:

\ssn

1) \ $\kxzk'_{x,\pel}=\widehat U_X\cdot\jmath_K(K_x)_\pel$. 

\ssn

2) \ $\lxlxzl{x}'$ the divisible hull of $\kxzk'_{x,\pel}$
in $\lxlxzl X$.   

\ssn
 
For $\tilde x\in\lxlxzl X$ let us consider some 
$u\in\widehat U_X$ and $x\in\jmath_K(K^\times)_\pel$
such that $\tilde x^m=u\cdot x$ for some $m>0$. Remark 
that since $\jmath_K(K_x)_\pel\cap\widehat U_X=\{1\}$,
the presentation above is unique in the following sense: 
If $\tilde x^m=u\cdot x$ and $\tilde x^{m'}=u'\cdot x'$, 
then $u^{m'}={u'}^m$ and $x^{m'}={x'}^m$. 

\ssn

For $\tilde x$ and $\tilde x^m=u\cdot x$ as above but 
with $x\in\jmath_K(K^\times)$, we denote:

\ssn

3) $\kxzk'_{\tilde x,\pel}:=\kxzk'_{x,\pel}$, and 
$\lxlxzl{\tilde x}':=\lxlxzl{x}'$.

\ssn

Remark that both $\kxzk'_{\tilde x,\pel}$ and 
$\lxlxzl{\tilde x}'$ depend only on $\tilde x$,
and not the representation $\tilde x^m=u\cdot x$. 

\ssn

N.B., for all $\tilde x$ there exists presentations as 
above with $x\in\jmath_K(K^\times)$.
}
\endproclaim

A next step is now to give a coarser version of the 
procedure indicated above for detecting $K_x$, but which
needs less information. Let namely $\Tell K$ be as usual
the subgroup of $\Gell K$ generated by all the inertia
subgroups of all Zariski prime divisors. And for $n=\ell^e$ 
consider 
$$
\Psi_n:K^\times/n\times K^\times/n\hor{\psi_n}
{}_n{\rm Br}\,(K)\hor{\rm res}{\rm H^2}(\Tell K,{\lv Z}/n)
$$
\nonumproclaim{Lemma 2} In the above notations for all 
$\tilde x\in\lxlxzl{X}$ one has:
$$
\tilde x^\perp:=\{\,\tilde y\in\lxlxzl X\mid
              \Psi_n(\tilde x,\tilde y)=0
                 \hbox{ for all $n=\ell^e$}\,\}=\lxlxzl{x}'
$$
\demo{Proof} For $\tilde y\in\tilde x^\perp$, 
consider presentations $\tilde x^m=u\cdot x$ and 
$\tilde y^{m'}=u'\cdot y$ as usual. Then we have:
$0=\Psi_n(\tilde x^m,\tilde y^{m'})=\Psi_n(x,y)$, as
$\widehat U_X$ is contained in both the right 
and the left kernel of $\Psi_n$ for all $n$. Thus 
$\Psi_n(x,y)=0$ for all $n$. As at the beginning of 
this Subsection, it follows that $K_x=K_y$. Equivalently 
we have $\kxzk'_{x,\pel}=\kxzk'_{y,\pel}$, thus also
$\lxlxzl{x}'=\lxlxzl{y}'$. Hence by the definition of 
$\lxlxzl{x}'$ it follows that $\tilde y\in\lxlxzl{x}'$.
Conversely, we show that $\lxlxzl{x}'$ is contained
in $\tilde x^\perp$. Indeed, each $\tilde y\in\lxlxzl{x}'$,
is of the form $\tilde y^{m'}=u'\cdot y$ for some 
integer $m'>0$, and some $u'\in\widehat U_X$, 
$y\in\kxzk'_{x,\pel}$. Thus $0=\psi_n(x,y)=
\psi_n(u\,x,u'\,y)=m\,m'\psi_n(\tilde x,\tilde y)$ 
for all $n$. Thus $\psi_n(\tilde x,\tilde y)=0$ for all
$n$, and finally $\tilde y\in\tilde x^\perp$.
\enddemo

\ssn

Thus Lemma~2 above shows that the following first two sets
are equal, and in a bijective correspondences with the third
one:
$$
\{\,\tilde x^\perp\mid\tilde x\in\lxlxzl{X}
                 \backslash\,\widehat U_X\,\}=
\{\,\lxlxzl{x}'\mid x\in\kxzk'_{X,\pel}
                 \backslash\,\widehat U_X\,\}\,\cong\,
\{\,\kxzk'_{x,\pel}\mid x\in K\backslash k\,\}
$$

\ssn

Now let some $\lxlxzl{\tilde x}'$ be given. By the
above correspondence, there exists a unique subfield 
$K_x$ of $K$ such that $\lxlxzl{\tilde x}'=\lxlxzl{x}'$,
or equivalently, that $\lxlxzl{\tilde x}'$ is the 
relative divisible hull of $\widehat U_X\cdot K_{x,\pel}$.
To simplify notations, let us set
$$
G:=\abab G_K,\quad G_x:=\abab G_{K_x}, \quad {\rm and} 
\quad p_x:G\to G_x
$$ 
the resulting canonical projection. And let
$\whph_x:\widehat K_x\to\widehat K$ be the embedding 
defined by $p_x$ via Kummer theory. By $\ell$-adic
Kummer theory, the knowledge of $p_x$ is equivalent
to the knowledge of $\whph_x$. We consider the unique 
complete smooth model $X_x\to k$ of $K_x|k$, and the 
corresponding geometric \cycellprefrm\ ${\cal G}_x$, 
etc.. We set $\kxzk_{x,\pel}=\widehat U_x\cdot K_{x,\pel}$ 
inside $\widehat K_x$. Then we have: 
$$
\kxzk_{x,\pel}:=\whph_x(\widehat U_x\cdot K_{x,\pel})
   =\whph_x(\widehat K_x)\cap\kxzk'_{x,\pel}
       =\whph_x(\widehat K_x)\cap\lxlxzl{x}'
$$ 

We next indicate how to identify $\kxzk_{x,\pel}$ inside 
$\lxlxzl{x}'$. Equivalently, by the discussion at the 
beginning of this Subsection, this gives then a recipe 
for detecting the $1$-dimensional projection $p_x$. 

\ssn
 
First, let $v$ be a Zariski prime divisor of $K$, and
$\jmath^v:\widehat K\to\Zell$ as introduced/defined in
Section~3, Definition/Remark 3.2,~3). Since $\jmath^v$ 
is trivial on $\widehat U_X$, it follows that $\jmath^v$
is trivial on $\lxlxzl{\tilde x}$ if and only if $v$
is trivial on $K_x$. 

\ssn

Next let $\jmath_v:\widehat U_v\to\widehat{Kv}$ be the
reduction map attached to a Zariski prime divisor $v$, 
see Section~3, Definition/Remark 3.2,~6). Then the 
restriction of $\jmath_v$ to the subgroup $K^\times_\pel$ 
coincides with the reduction map defined/introduced in 
Section 1,~A), Definition/Remark 2.4,~3). In particular, 
since $\td(K_x|k)=1$, one has: $\jmath^v$ is not trivial 
on $K_x$ if an only if $\jmath_v$ is trivial on $K_x$.

\nonumproclaim{Lemma 3} In the above context we have:
$$
\kxzk_{x,\pel}=\{\,\tilde y\in\lxlxzl{x}'\mid
    \jmath_v(y)=1\hbox{ for all } v \hbox{ with }
        \jmath^v(\lxlxzl{x}')\neq0\,\}
$$
\endproclaim 

\demo{Proof} The direct inclusion is clear by 
the discussion above. In order to prove the converse,  
we view the regular field extension $K|K_x$ as a 
function field of positive transcendence degree
$\td(K|K_x)=\td(K|k)-1$ over the base field $K_x$. 
Let $K''|K$ be a finite extension, and $K''_x$ the 
relative algebraic closure of $K_x$ in $K''$. Set 
$K':=KK''_x$ inside $K''$. Since $K''_x$ is Hilbertian, 
there exists ``many'' Zariski prime divisor $v$ of $K''$ 
with the following properties: 

\ssn

i) $v$ is not trivial on $K''_x$, thus on $K_x$.

\ssn

ii) $[K''v: K'v]=[K'':K']$

\ssn

One applies the fact above to Abelian extensions of $K$ 
which are linearly disjoint from $\fell K_x$. 
\enddemo
%
%
%
%
\proclaim{Proposition}
In the above notations, consider the 
perfect $\ell$-adic pairing given by Kummer theory
$<\,,>:\widehat K\times\abab{G_K}\to\Zell$.

{\rm(1)} The kernels of $1$-dimensional projections 
$\abab p_x:G\to G_x$ are exactly the closed subgroups 
${\cal N}_x$ of $G$ of the form: 
$$
{\cal N}_x=\{\sigma\in G\mid\>
      <\tilde y\,,\sigma>\>\>=0\hbox{ for all } 
         \tilde y\in\kxzk_{x,\pel}\,\}
$$

Let $p_x:G\to G/{\cal N}_x=:G_x$ be a $1$-dimensional 
projection. Then one has:

\ssn

{\rm(2)} A procyclic subgroup $T_{v_x}$ of $G_x$ is an 
inertia group if and only if there exits an inertia
group $T_v$ of $G$ such that $p_x(T_v)$ has finite 
index in $T_{v_x}$. 

\ssn

{\rm3)} Finally, $p_x:G\to G_x$ is a rational projection 
if and only if the inertia groups  $T_{v_x}$ (all $v_x$) 
generate $G_x$. 

\ssn

Clearly, all the facts above are group theoretically 
encoded in $\Gell K$.  
\endproclaim
\label\ratproj\theoremnum*

\demo{Proof} Clear by the observations above and the 
characterization of the $1$-dimensional projections via
$K_\pel$ mentioned above.
\enddemo

\ssn

B) \ {\it The proof of Theorem (Introduction)\/} 
 
\ssn

By Proposition~\charprop\ it follows that the geometric
\cycellprefrm s over $K|k$ are encoded in $\Gell K$. And
it is clear that the group theoretic description of the
geometric \cycellprefrm s is preserved under isomorphisms.
Therefore, every isomorphism
$$
\Phi:\Gell K\to \Gell L
$$
defines a bijection from the set of all geometric
\cycellprefrm s over $K|k$ onto the set of all geometric
\cycellprefrm s over $L|l$. Further, if ${\cal G}_Y$ 
corresponds to ${\cal G}_X$ under this bijection, then 
$\Phi$ gives rise to an isomorphism of divisorial
\cycellprefrm s
$$
\Phi:{\cal G}_X\to{\cal G}_Y\,.
$$

Now suppose that $X$ is complete regular like.
Since this fact is encoded group theoretically in 
$\Gell K$ and the projection $\Gell K\to G:=\abab{G_K}$,
it follows that the corresponding $Y$ is complete regular 
like too. Consider the canonical commutative diagram
$$
\matrixbaselines
\matrix{
0\to\hhb6\widehat U_Y
    &\shto&\lxlxzl{Y} 
             &\shto&\mardv(Y)_\pel
                      &\shto&{\eu Cl}_{Y,\pel}
                                              \,\to0    \cr
\hhb{25}\dwn{\hat\psi}&&\dwn{\hat\psi}\hhb5
                           &&\dwn{}\hhb5&&\dwn{}\hhb{30}\cr 
0\to\hhb6\widehat U_X
    &\shto&\lxlxzl{X} 
             &\shto&\mardv(X)_{\pel}
                      &\shto&{\eu Cl}_{X,\pel}
                                              \,\to0    \cr
}
\leqno{\indent(*)_\pel}
$$
where the vertical map $\hat\psi$ is such a multiple 
$\hat\psi=\epsilon\,\whph$ of the Kummer homomorphism 
$\whph:\widehat L\to\widehat K$ that one has:
$\psi\big(\lxlxzl Y\big)=\lxlxzl X$. 

\ssn

We next remark that by Proposition~\ratproj\ above, 
$\Phi$ also defines a bijection from the set of all 
the rational projections of ${\cal G}_X$ onto the set 
of all the rational projections of ${\cal G}_Y$. 

\ssn

On the other hand, $K_{(\ell)}$ is encoded in $\Gell K$
as follows: The multiplicative group $K^\times$ is 
generated by the set of all the $K^\times_x$ with 
$K_x$ a rational function field. In other words, if
$(p_x)_x$ is the family of all the rational projections
as above, then $K_\pel\subset\widehat K$ is generated 
by all the images $\kxzk_{x,\pel}$. Therefore, in 
order to detect the $\ell$-adic equivalence class 
of the lattice $K_\pel$ inside $\widehat K$, one can
proceed as follows: Choose a divisorial lattice 
$\lxlxzl K$ in $\widehat K$. For every rational 
projection $p_x$, consider the unique divisorial 
lattice $\lxlxzl{K_x}$ in $\widehat K_x$ whose image 
$\kxzk_{x,\pel}$ under the completion morphism 
$\whph_x:\widehat K_x\to\widehat K$ is contained in 
$\lxlxzl K$. Then the subgroup $\lxlxzl{K}^0$ of 
$\lxlxzl K$ which is generated by all the $\kxzk_{x,\pel}$
is the unique lattice in $\lxlxzl K$ which is equivalent 
to $K_\pel$. 

\ssn

Since the description of $\lxlxzl{K}^0$ above is given 
in pure group theoretic terms, it is clear that the 
Kummer homomorphism $\whph_K:\widehat L\to\widehat K$ 
maps $L_\pel$ isomorphically onto a lattice $\lxlxzl{K}^0$ 
in $\widehat K$ which is $\ell$-adically equivalent to 
$K_\pel$. Thus after multiplying $\whph$ by a properly
chosen $\ell$-adic unit (unique modulo rational
$\ell$-adic units), we can suppose that $\whph$ maps
$L_\pel$ isomorphically on $K_\pel$.

\ssn

We finally apply \nmnm{Pop} [P4], Theorem~5.11, 
and deduce that there exist finite pure inseparable 
extensions $K_0|K$ and $L_0|L$ such that  
$\whph_{K_0}:\widehat{L_0}\to\widehat{K_0}$ is the 
$\ell$-adic completion of a unique field isomorphism 
$\imath_0:K_0\to L_0$. 

\ssn

Moreover, let $K'|K$ be a finite (Galois) 
extension of $K$ inside $\fell K$, and via the 
isomorphism $\Phi$ the corresponding finite (Galois) 
extension $L'|L$. Then using the functoriality of 
Kummer theory, we have a commutative diagram of the 
form

$$
\def\normalbaselines{
       \baselineskip20pt\lineskip3pt\lineskiplimit3pt}
\matrix{\widehat L&\horr{\whph_{K}}&\widehat K&\cr
                      \dwn{\rm can}&&\dwn{\rm can}&\cr
        \widehat{L'}&\horr{\whph_{K'}}& \widehat{K'}&\cr}
$$
Suppose $K'|K$ is 
Galois, and $\Phi':{\rm Gal}(K'|K)\to{\rm Gal}(L'|L)$
is the isomorphism induced by $\Phi$. Then for all 
$\sigma\in{\rm Gal}(K'|K)$, and $\tau=\Phi(\sigma)$ we 
have: 
$$
\sigma\circ\whph_{K'}=\whph_{K'}\circ\tau
  \quad{\rm i.e.,}\quad 
    \Phi(\sigma)=\whph_{K'}^{-1}\circ\sigma\circ\whph_{K'}\,.
$$ 

So if $\whph_K$ maps $L_\pel$ into $K_\pel$, then 
$\whph_{K'}$ automatically maps $L'_\pel$ into $K'_\pel$,
i.e., it is not necessary to ``re-norm'' $\whph_{K'}$ in
order to have $\whph_{K'}(L'_\pel)=K'_\pel$. Thus the
above commutative diagram induces a corresponding one 
with $K_\pel$ in stead of $\widehat K$, etc.. 

\msn
 
{\it Conclusion:\/} 

\ssn

Now taking limits over all finite extensions $K'|K_0$ inside 
$\fell K_0=\fell K K_0$, and the corresponding finite 
extensions $L'|L_0$ inside $\fell L_0$, we finally get an 
isomorphism $\phi:\fell{L^{\rm i}}\to\fell{K^{\rm i}}$
defining $\Phi$, i.e., $\Phi(g)=\phi^{-1}\,g\,\phi$ for 
all $g\in\Gell{K}$. 

\ssn

We still have to prove that any two 
automorphisms $\phi'$ and $\phi''$ both defining 
$\Phi$ differ by a power of Frobenius. Indeed, setting 
$\phi:=\phi''\circ{\phi'}^{-1}$ we get: $\phi$ is 
an automorphism of $\fell{L^{\rm i}}$ which maps 
$L^{\rm i}$ onto itself, and induces the identity on 
$\Gell{L^{\rm i}}$. 
We claim that such an automorphism $\phi$ is a power of
Frobenius. Indeed, let $v$ be an arbitrary Zariski prime
divisor of $\fell{L^{\rm i}}$. Then $w:=v\circ\phi$ is
also a Zariski prime divisor of $\fell{L^{\rm i}}$. Moreover,
by the usual formalism, we have $Z_w=\phi^{-1}Z_v\,\phi$ 
inside $\Gell{L^{\rm i}}$. Since the conjugation by $\phi$
is the identity on $\Gell{L^{\rm i}}$, it follows that $Z_w=Z_v$. 
Thus by Proposition \propdecomgr,~(1), it follows that $v$ and 
$w$ are equivalent valuations on $\fell{L^{\rm i}}$. In 
particular, $v(x)>0$ if and only if $w(x)>0$ (all 
$x\in\fell{L^{\rm i}}$). Let ${\rm char}(k)=p$.  We claim 
that $y:=\phi(x)$ is some $p$-power of $x$. First, if $y$ 
and $x$ are algebraically independent, then there exists a 
Zariski prime divisor $v$ of $\fell{L^{\rm i}}$ such that 
$v(x)=1$ and $v(y)=0$. But $v(y)=v\circ\phi(x)=w(x)$,
contradiction! Therefore, $y=\phi(x)$ is always algebraic
integer over $k[x]$, and vice-versa: $x$ is algebraic 
integer over $k[y]$. Now let $f(X,Y)\in k[X,Y]$ be the 
minimal polynomial polynomial relation between $x$ and 
$y$ over $k$. We claim that $f(X,Y)$ is pure inseparable
in $Y$ (and thus by symmetry, also in $X$). Let namely
$y_a=:y_1,\dots,y_r$ be the distinct roots of $f(x,Y)$
over $k(x)$. The for a ``general $a\in k$ we have: If
$x_a=x+a$, then $y_a=:\phi(x_a)=y+\phi(a)$. Thus $x_a$ 
and $y_a$ have the same zeros in ${\eu D}^1_{\fell K}$.
And clearly, the minimal polynomial polynomial relation 
satisfied by $x_a$ and $y_a$ over $k$ is 
$f_a(X,Y)=f\big(X-a,Y-\phi(a)\big)$. Therefore, for
a general $a\in k$, the polynomial $f_a(0,Y)$ has $s$ 
distinct roots $b_1,\dots, b_s$. And the specializations
$(x,y)\mapsto(a,b_i)$ give rise to $s$ different places 
$v_i$ of $\fell{k(x)}$ all of which are zeros of $x_a$. 
But then these places must also be zeros of $y_a$. 
Thus $0=y_a(b_i)=\phi(a)-b_i$ implies $b_i=\phi(a)$,
thus $s=1$. Thus finally $f(X,Y)$ is purely inseparable 
in $Y$, and by symmetry also in $X$. 

\ssn

Thus finally we have: $\phi(x)=a_x\,x^{n_x}$, where 
$a_x\in k$ and $n_x=p^{e_x}$ depend on $x$. We next 
show that $a_x=1$ and that $n_x=p^e$ do not depend on 
$x$. Indeed, since $\phi$ is a field automorphism, 
we have $\phi(x+y)=\phi(x)+\phi(y)$ for all $x,y$ 
from $\fell{L^{\rm i}}$. Or equivalently:
$a_{x+y}(x+y)^{n_{x+y}}=a_x\,x^{n_x}+a_yy^{n_y}$.
Considering several arbitrary $y\in\fell{L^{\rm i}}$ such 
that $x$ and $y$ are algebraically independent over $k$, we 
therefore must have $n_{x+y}=n_y=n_x=p^e$ and $a_{x+y}=a_y=a_x=a$,
thus independent of $x$ and $y$. Since $a_1=1$, we have 
$a=1$. Finally, choosing a transcendence basis ${\cal T}$ 
of $L$ over $k$ (say, which contains $x$), we see that 
the restriction of $\phi$ to $k({\cal T})$ is ${\rm Frob}^e$.
Therefore, $\phi={\rm Frob}^e$ on $\fell{L^{\rm i}}$.

\ssn

The Theorem (Introduction) is proved.
%
%
%
\AuthorRefNames[F--H--V]

\references

[A] 
\name{Artin, E.}, Geometric algebra, Interscience Publishers, Inc.,
      New York 1957.

[Bo]
\name{Bogomolov, F.\ A.}, {\it On two conjectures in birational
     algebraic geometry,\/} in Algebraic Geometry and Analytic
     Geometry, ICM-90 Satellite Conference Proceedings, ed
     A.~Fujiki et all, Springer Verlag Tokyo 1991. 

[B--T\hhb{1}1]
\name{Bogomolov, F.\ A.}, and \name{Tschinkel, Y.} {\it Commuting 
     elements in Galois groups of function fields,\/} in: Motives,
     Polylogs and Hodge Theory, International Press 2002, 75--120. 

[B--T\hhb{1}2]
\name{Bogomolov, F.\ A.}, and \name{Tschinkel, Y.} {\it 
    Reconstruction of function fields,\/} Manu\-script, 
    March 2003.

[BOU]
\name{Bourbaki, N.}, Alg\`ebre commutative, Hermann Paris 1964.

[D1]
\bibline, {\it Le groupe fondamental de la droite projective
     moins trois points,\/} in: Galois groups over ${\bf Q},$ Math.\
     Sci.\ Res.\ Inst.\ Publ. {\bf16}, 79--297, Springer 1989.

[D2]
\name{Deligne, P.} {Letter to the Author, June 1995.}

[Ef]
\name{Efrat, I,}, {\it Construction of valuations from K-theory,\/} 
    Mathematical Research Letters 6 (1999), 335-344.

[E--K]
\name{Engler, A.\ J.} and \name{Koenigsmann, J.}, {\it Abelian subgroups 
    of pro-$p$ Galois groups,\/} Trans.\ AMS {\bf350} (1998), no. 6, 
    2473--2485.

[F]
\name{Faltings, G.}, {\it Endlichkeitss\"atze f\"ur abelsche
Variet\"aten \"uber Zahlk\"orpern,\/} Invent.\ Math. {\bf 73}
    (1983), 349--366.

[G1]
\name{Grothendieck, A.}, Letter to Faltings, June 1983.

[G2]
\bibline, Esquisse d'un programme, 1984.

[Ha]
\name{Harbater, D.} {\it Fundamental groups and embedding problems
    in characteristic $p$,\/} in Recent developments in the inverse
    Galois problem, eds M.~Fried et al, Contemp.\ Math.\ Series
    {\bf186} (1995).

[Hr]
\name{Hartshorne, R.} Algebraic Geometry, GTM 52, Springer 
     Verlag 1993. 

[H]
\name{Hironaka, H.}, {\it Resolution of singularities of an
     algebraic variety over a field of characteristic zero,\/}
     Ann.~of Math. {\bf 79} (1964), 109--203; 205--326.

[I]
\name{Ihara, Y.}, {\it On Galois representations arising from towers 
    of covers of\/} ${\scriptstyle{\lv P}^1\backslash\{0,1,\infty\}}$,
    Invent.\ math. {\bf 86} (1986), 427--459.

[Ik]
\name{Ikeda, M.}, {\it Completeness of the absolute Galois group
     of the rational number field,\/} J.~reine angew.\ Math.
     {\bf291} (1977), 1--22.

[J]
\name{de Jong, A.\ J.}, {\it Families of curves and alterations,\/}
      Annales de l'institute Fourier, {\bf47} (1997), pp. 599--621.

[K--L]
\name{Katz, N.}, and \name{Lang, S.} {\it Finiteness theorems in 
     geometric class field theory,\/} with an Appendix by K.~Ribet,
     Enseign.\ Math.\ {\bf 27} (1981), 285--319.

[K--K] 
\name{Knaf, H.} and \name{Kuhlmann, F.-V.} {\it Abhyankar places 
     admit local uniformization in any characteristic,\/} Manuscript,
     see \ {http://mathsci.usask.ca/~fvk/Fvkprepr.html}

[Ko]
\name{Koenigsmann, J.}, {\it From $p$-rigid elements to valuations
     (with a Galois characterisation of $p$-adic fields),\/} 
     J.\ reine angew.\ Math.,  465  (1995), 165--182.

[Km]
\name{Komatsu, K.}, {\it A remark to a Neukirch's conjecture,\/}
     Proc.\ Japan Acad. {\bf50} (1974), 253--255.

[Mo]
\name{Mochizuki, Sh.}, {\it The local pro-p Grothendieck conjecture 
     for hyperbolic cur\-ves,\/} Invent.\ Math.\ {\bf 138} (1999), 
     319-423.

[N]
\name{Nagata, M.}, {\it A theorem on valuation rings and its
     applications,\/} Nagoya Math.~J. {\bf 29} (1967), 85--91.

[Na]
\name{Nakamura, H.}, {\it Galois rigidity of the \'etale fundamental
     groups of punctured projective lines,\/} J.\ reine angew.\ Math.
     {\bf411} (1990) 205--216.

[N1]
\name{Neukirch, J.}, {\it \"Uber eine algebraische Kennzeichnung der
     Hen\-sel\-k\"or\-per,\/} J.\ reine angew.\ Math. {\bf231} (1968),
     75--81.

[N2]
\bibline, {\it Kennzeichnung der $p$-adischen und endlichen algebraischen
     Zahl\-k\"or\-per,\/} Inventiones math.\ {\bf6} (1969), 269--314.

[N3]
\bibline, {\it Kennzeichnung der endlich-algebraischen Zahl\-k\"or\-per
     durch die Galois\-grup\-pe der maximal aufl\"os\-baren
     Erweiterungen,\/} J.\ f\"ur Math. {\bf238} (1969), 135--147.

[O]
\name{Oda, T.}, {\it A note on ramification of the Galois
     representation of the fundamental group of an algebraic
     curve I,\/} J.\ Number Theory (1990) 225--228.

[Pa]
\name{Parshin, A. N.}, {\it Finiteness Theorems and Hyperbolic 
     Manifolds,\/} in: The Grothen\-dieck Festschrift III, ed 
     P.~Cartier et all, PM Series Vol 88, Birkh\"auser Boston Basel 
     Berlin 1990.

[P1]
\name{Pop, F.}, {\it On Grothendieck's conjecture of birational 
     anabelian geometry,\/} Ann.\ of Math.\ {\bf 138} (1994), 145--182.

[P2]
\bibline, {\it On Grothendieck's conjecture of birational 
     anabelian geometry II,\/} Heidel\-berg--Mannheim Preprint series 
     Arithmetik~II, ${\rm N}^0\,{\bf 16},$ Heidelberg 1995.

[P3]
\bibline, {\it Alterations and birational anabelian geometry,\/} 
      in: Resolution of singularities, Birkh\"auser PM Series, 
      Vol.\hhb2{\bf 181}, p.\hhb{2}519--532; eds Herwig Hauser et 
      al, Birkh\"auser Verlag, Basel 2000.

[P4]
\name{Pop, F.}, {\it The birational anabelian conjecture
      revisited,\/} Manuscript, Oct 2002.

[Po]
\name{Pop, F.}, {\it \'Etale Galois covers of affine, smooth
     curves,\/} Invent.\ Math. {\bf120} (1995), 555--578.

[R]
\name{Roquette, P.}, {\it Zur Theorie der Konstantenreduktion
     alge\-brai\-scher Mannig\-fal\-tig\-keiten,\/} J.\ reine
     angew.\ Math. {\bf200} (1958), 1--44.

[S]
\name{Serre, J.-P.}, {Cohomologie Galoisienne,} LNM 5,
     Springer 1965.

[Sp]
\name{Spiess, M.}, {\it An arithmetic proof of Pop's Theorem
     concerning Galois groups of function fields over number
     fields,\/} J.\ reine angew.\ Math.\ {\bf478} (1996), 107--126.
[T]
\name{Tamagawa, A.}, {\it The Grothendieck conjecture for affine 
      curves,\/} Compositio Math.\ {\bf 109} (1997), 135-194.

[U1]
\name{Uchida, K.}, {\it Isomorphisms of Galois groups of algebraic
     function fields,\/} Ann.\ of Math. {\bf106} (1977), 589--598.

[U2]
\bibline, {\it Isomorphisms of Galois groups of solvably closed Galois
     extensions,\/} T\^oho\-ku Math.\ J. {\bf31} (1979), 359--362.

[U3]
\bibline, {\it Homomorphisms of Galois groups of solvably closed Galois
     extensions,\/} J.~Math.\ Soc.\ Japan {\bf33}, No.4, 1981.

[W]
\name{Ware, R.}, {\it Valuation Rings and rigid Elements in Fields,\/}
     Can.\ J.\ Math. {\bf 33} (1981), 1338--1355.

[Z--S] 
\name{Zariski, O.}, and \name{Samuel, P.} Commutative Algebra,
      Vol II, Springer-Verlag, New York, 1975.
\endreferences
\bye